\documentclass[12pt,leqno]{article}
\usepackage{amssymb}
\usepackage[mathscr]{eucal}
\usepackage{amsmath,amssymb,latexsym,theorem,bbm}
\usepackage{color}
\usepackage{appendix, url}

\newcommand{\DS}{\displaystyle}

\newcommand{\CC}{\mathbb{C}}

\newcommand{\NN}{\mathbb{N}}

\newcommand{\RR}{\mathbb{R}}

\newcommand{\ZZ}{\mathbb{Z}}

\newcommand{\bA}{{\boldsymbol{A}}}

\newcommand{\ta}{\widetilde{a}}
\newcommand{\tA}{\widetilde{A}}

\newcommand{\hA}{\widehat{A}}
\newcommand{\hcA}{\widehat{\cA}}

\newcommand{\bB}{{\boldsymbol{B}}}
\newcommand{\hB}{\widehat{B}}

\newcommand{\tB}{\widetilde{B}}

\newcommand{\be}{{\boldsymbol{e}}}

\newcommand{\bI}{{\boldsymbol{I}}}

\newcommand{\bN}{{\boldsymbol{N}}}

\newcommand{\bP}{{\boldsymbol{P}}}

\newcommand{\bQ}{{\boldsymbol{Q}}}

\newcommand{\bR}{{\boldsymbol{R}}}

\newcommand{\bS}{{\boldsymbol{S}}}

\newcommand{\bV}{{\boldsymbol{V}}}

\newcommand{\bx}{{\boldsymbol{x}}}
\newcommand{\bX}{{\boldsymbol{X}}}
\newcommand{\by}{{\boldsymbol{y}}}
\newcommand{\bY}{{\boldsymbol{Y}}}

\newcommand{\bZ}{{\boldsymbol{Z}}}
\newcommand{\bU}{{\boldsymbol{U}}}

\newcommand{\bkappa}{{\boldsymbol{\kappa}}}

\newcommand{\Beta}{{\boldsymbol{\eta}}}

\newcommand{\btheta}{{\boldsymbol{\theta}}}

\newcommand{\bzeta}{{\boldsymbol{\zeta}}}
\newcommand{\bzero}{{\boldsymbol{0}}}

\newcommand{\cA}{{\mathcal A}}
\newcommand{\cB}{{\mathcal B}}
\newcommand{\cC}{{\mathcal C}}
\newcommand{\cD}{{\mathcal D}}

\newcommand{\cI}{\mathcal{I}}

\newcommand{\cG}{{\mathcal G}}
\newcommand{\cF}{{\mathcal F}}

\newcommand{\cN}{{\mathcal N}}

\newcommand{\tbS}{\widetilde{\bS}}
\newcommand{\tbN}{\widetilde{\bN}}

\newcommand{\cc}{\mathrm{c}}
\newcommand{\dd}{\mathrm{d}}
\newcommand{\ee}{\mathrm{e}}
\newcommand{\ii}{\mathrm{i}}

\newcommand{\EE}{\operatorname{\mathbb{E}}}

\newcommand{\PP}{{\operatorname{\mathbb{P}}}}
\newcommand{\QQ}{{\operatorname{\mathbb{Q}}}}

\newcommand{\var}{\operatorname{Var}}

\newcommand{\diag}{\operatorname{diag}}

\DeclareMathOperator*{\argmin}{arg\,min}

\newcommand{\hvarrho}{\widehat{\varrho}}

\newcommand{\tM}{\widetilde{M}}
\newcommand{\tN}{\widetilde{N}}

\newcommand{\tY}{\widetilde{Y}}

\newcommand{\tZ}{\widetilde{Z}}

\newcommand{\vare}{\varepsilon}

\renewcommand{\mid}{\,|\,}

\renewcommand{\leq}{\leqslant}
\renewcommand{\geq}{\geqslant}

\newcommand{\stoch}{\stackrel{\PP}{\longrightarrow}}
\newcommand{\stochG}{\stackrel{\PP_G}{\longrightarrow}}

\newcommand{\stochGeta}{\stackrel{\PP_{G\cap\{\exists\Beta^{-1}\}}}{\longrightarrow}}

\newcommand{\distrP}{\stackrel{\cD(\PP)}{\longrightarrow}}
\newcommand{\distrQ}{\stackrel{\cD(\QQ)}{\longrightarrow}}

\newcommand{\distre}{\stackrel{\cD}{=}}

\newcommand{\meanP}{\stackrel{L_1(\PP)}{\longrightarrow}}

\newcommand{\qmeanP}{\stackrel{L_2(\PP)}{\longrightarrow}}

\newcommand{\asP}{\stackrel{{\text{$\PP$-a.s.}}}{\longrightarrow}}

\newcommand{\aseP}{\stackrel{{\text{$\PP$-a.s.}}}{=}}

\newcommand{\bbone}{\mathbbm{1}}

\newcommand{\proofend}{\hfill\mbox{$\Box$}}

\numberwithin{equation}{section}

\theoremstyle{change} \theorembodyfont{\em}
\newtheorem{Lem}{Lemma.}[section]
\newtheorem{Thm}[Lem]{Theorem.}
\newtheorem{Pro}[Lem]{Proposition.}
\newtheorem{Cor}[Lem]{Corollary.}
\newtheorem{Def}[Lem]{Definition.}

\theorembodyfont{\rm}
\newtheorem{Rem}[Lem]{Remark.}

\begin{document}

\begin{center}
 {\bfseries\Large
   Stable convergence of conditional least squares\\[0.5mm]
   estimators for supercritical continuous state and
   continuous time branching processes with immigration}

\vspace*{3mm}

 {\sc\large
  M\'aty\'as $\text{Barczy}^{*}$
  }

\end{center}

\vskip0.2cm

\noindent
 * HUN-REN--SZTE Analysis and Applications Research Group,
     Bolyai Institute, University of Szeged,
     Aradi v\'ertan\'uk tere 1, H--6720 Szeged, Hungary.

\noindent E--mail: barczy@math.u-szeged.hu


\vskip0.5cm


{\renewcommand{\thefootnote}{}
\vspace*{0.2cm}
\footnote{ M\'aty\'as Barczy was supported by the project TKP2021-NVA-09.
Project no.\ TKP2021-NVA-09 has been implemented with the support
 provided by the Ministry of Innovation and Technology of Hungary from the National Research, Development and Innovation Fund,
 financed under the TKP2021-NVA funding scheme.}}

\vspace*{-10mm}

\begin{abstract}
We prove stable convergence of conditional least squares estimators of drift parameters
 for supercritical continuous state and continuous time branching processes with
 immigration based on discrete time observations.
\end{abstract}

\section{Introduction}
\label{section_intro}

Statistical inference for branching processes has a long tradition and history, see, e.g.,
 the books of Guttorp \cite{Gut} and  Gonz\'alez et al.\ \cite[Section 5]{GonPueYan}.
In this paper, we study the asymptotic behavior of conditional least squares (CLS) estimators
 of drift parameters for supercritical continuous state and continuous time branching processes with
 immigration (CBI processes) based on discrete time (low frequency) observations.
According to our knowledge, results on stable (mixing) convergence of CLS estimators
 for parameters of CBI processes are not available in the literature,
 all the existing results state convergence in distribution of the appropriately normalized CLS estimators in question.
For supercritical discrete time Galton-Watson branching processes, H\"ausler and Luschgy \cite{HauLus} proved
 stable convergence of the CLS estimator of the offspring mean under non-extinction (that we recall below),
which served us as a motivation for investigating the problem for CBI processes.

For supercritical (discrete time) Galton-Watson branching processes (without immigration), H\"ausler and Luschgy \cite[Corollaries 10.2, 10.4 and 10.6]{HauLus}
 established stable and mixing convergences of the appropriately scaled conditional moment estimator (also called Lotka-Nagaev estimator),
 CLS estimator and weighted CLS estimator (also called Harris estimator) of the offspring mean under non-extinction.
More precisely, given a supercritical Galton-Watson branching process \ $(Y_n)_{n\geq 0}$, \
 the Lotka-Nagaev estimator \ $\widehat \alpha_n^{\mathrm{(LN)}}$, \ the CLS estimator \ $\widehat \alpha_n^{\mathrm{(CLS)}}$ \ and
 the Harris estimator \ $\widehat \alpha_n^{\mathrm{(H)}}$ \ of the offspring mean \ $\alpha$ \
 based on the observations \ $Y_0,Y_1,\ldots,Y_n$ \ (where \ $n\geq 1$) \ take the following forms:
 \begin{align*}
  \widehat \alpha_n^{\mathrm{(LN)}}:=\frac{Y_n}{Y_{n-1}},\qquad\;
  \widehat \alpha_n^{\mathrm{(CLS)}}:=\frac{\sum_{k=1}^n Y_{k-1}Y_k}{\sum_{k=1}^n Y_{k-1}^2},\qquad\;
  \widehat \alpha_n^{\mathrm{(H)}}:=\frac{\sum_{k=1}^n Y_k}{\sum_{k=1}^n Y_{k-1}},
 \end{align*}
 provided that \ $Y_{n-1}\geq 1$ \ for \ $\widehat \alpha_n^{\mathrm{(LN)}}$, \ and \ $Y_0\geq 1$ \
 for \ $\widehat \alpha_n^{\mathrm{(CLS)}}$ \ and \ $\widehat \alpha_n^{\mathrm{(H)}}$ \ (yielding that
 \ $\sum_{k=1}^n Y_{k-1}^2 \geq 1$ \ and \ $\sum_{k=1}^n Y_{k-1} \geq 1$), \ respectively.
Supposing that the initial value \ $Y_0$ \ is positive integer-valued, \ $Y_0$ \ and the offspring distribution have finite fourth moments, the variance
 \ $\beta^2$ \ of the offspring distribution is positive and that the offspring mean \ $\alpha>1$ \ (supercritical case),
  H\"ausler and Luschgy \cite[Corollary 10.4]{HauLus} proved that
 \[
   \frac{(\alpha^3-1)^{1/2}}{\alpha^2-1}\alpha^{n/2}
      (\widehat \alpha_n^{\mathrm{(CLS)}} - \alpha)
  \to \beta M_\infty^{-1/2} N
  \qquad \text{$\cF_\infty^Y$-stably under \ $\PP_{\{M_\infty > 0\}}$ \ as \ $n \to \infty$,}
 \]
 where \ $\cF_\infty^Y:=\sigma\Big(\bigcup_{n=0}^\infty \sigma(Y_0,Y_1,\ldots,Y_n)\Big)$; \
 \ $M_\infty$ \ is an \ $\cF_\infty^Y$-measurable random variable satisfying \ $\EE(M_\infty^2)<\infty$ \ and
  \ $\alpha^{-n} Y_n\asP M_\infty$ \ as \ $n\to\infty$; \ $N$ \ is a standard normally distributed random
 variable \ $\PP$-independent of \ $\cF_\infty^Y$; \ $\PP_{\{M_\infty >0\}}$ \ denotes the conditional
 probability measure given \ $\{M_\infty >0\}$; \ $\PP(\lim_{n\to\infty} Y_n = \infty) = \PP( M_\infty >0) >0$
  \ (i.e., the probability of non-extinction is positive); and for the notion of stable convergence, see, e.g., Definition \ref{Def_HL_stable_conv}.
H\"ausler and Luschgy \cite[Corollary 10.4]{HauLus} also established \ $\cF_\infty^Y$-mixing convergence of
 \ $\widehat \alpha_n^{\mathrm{(CLS)}}$ \ using a random scaling.

Limit theorems stating stable (mixing) convergence instead of convergence in distribution are important not only from
 theoretical point of view.
These types of limit theorems have such statistical applications where limit theorems stating only
 convergence in distribution cannot be directly used.
Such a nice application for the description of the asymptotic behaviour of the above mentioned Harris estimator of the offspring mean of
 a supercritical Galton--Watson branching process is explained on pages 3 and 4 in H\"ausler and Luschgy \cite{HauLus}.
In the heart of this application there is a generalization of Slutsky's lemma in a way that convergence in distribution is replaced by stable convergence
 and in return the limit of the stochastically convergent sequence can be also random not only a deterministic constant
 (see, e.g., Theorem 3.18 in H\"ausler and Luschgy \cite{HauLus} or Theorem \ref{Thm_HL_Thm3_18}).

Under some moment conditions (given below), a CBI process can be represented as a pathwise unique strong solution of the
 stochastic differential equation (SDE)
 \begin{align}\label{SDE_atirasa_dim1}
  \begin{split}
   X_t
   &=X_0
     + \int_0^t (a + B X_s) \, \dd s
     + \int_0^t \sqrt{2 c \max \{0, X_s\}} \, \dd W_s \\
   &\quad
      + \int_0^t \int_0^\infty \int_0^\infty
         z \bbone_{\{u\leq X_{s-}\}} \, \tN(\dd s, \dd z, \dd u)
      + \int_0^t \int_0^\infty r \, M(\dd s, \dd r)
  \end{split}
 \end{align}
 for \ $t \in [0, \infty)$, \ where \ $X_0$ \ is a nonnegative random variable with \ $\EE(X_0)<\infty$,
 \ $a, c \in [0, \infty)$, \ $B \in \RR$, \  $(W_t)_{t\geq0}$ \ is a standard Wiener process,
 \ $N$ \ and \ $M$ \ are Poisson random measures on
 \ $(0, \infty)^3$ \ and on \ $(0, \infty)^2$ \ with intensity measures
 \ $\dd s \, \mu(\dd z) \, \dd u$ \ and \ $\dd s \, \nu(\dd r)$,
 \ respectively,
 \ $\tN(\dd s, \dd z, \dd u)
    := N(\dd s, \dd z, \dd u) - \dd s \, \mu(\dd z) \, \dd u$ \ is the compensated Poisson random measure
    corresponding to \ $N$,
 \ the branching jump measure $\mu$ and the immigration jump measure
 $\nu$ are Borel measures on $(0,\infty)$ (furnished with the Borel $\sigma$-algebra) that
 satisfy the moment conditions $\int_0^\infty (z\wedge z^2)\,\mu(\dd z)<\infty$ and $\int_0^\infty z\,\nu(\dd z)<\infty$,
 respectively, and $X_0$, $(W_t)_{t\geq0}$, $N$ and $M$ are independent, see Dawson and Li
 \cite[Theorems 5.1 and 5.2]{DawLi1} or Barczy et al.\ \cite[Sections 4 and 5]{BarLiPap2}.
For more details on CBI processes, see Section \ref{section_CBI}.
A CBI process $(X_t)_{t\geq 0}$ given as the pathwise unique strong solution of the SDE \eqref{SDE_atirasa_dim1}
 is called subcritical, critical or supercritical if $B < 0$, \ $B = 0$ or $B > 0$, see, e.g., Huang et al.\ \cite[page 1105]{HuaMaZhu},
 and hence one can call \ $B$ \ as the criticality parameter.

Below, we recall some known results on the asymptotic behavior of CLS estimators for CBI processes,
 all these results state convergence in distribution of the appropriately normalized CLS estimators in question.

Huang et al.\ \cite{HuaMaZhu} studied the weighted CLS estimator of \ $(B, a)$ \ based on
 discrete time observations \ $(X_{k\Delta t})_{k\in\{1,\ldots,n\}}$, \ $n \in \{1, 2, \ldots\}$, \ where \ $\Delta t \in (0, \infty)$ \ is fixed,
 supposing that the initial value \ $X_0$ \ is a nonnegative constant.
Under second order moment conditions on the branching and immigration mechanisms,
 in the supercritical case (i.e., when $B>0$) they showed that the weighted CLS estimator of \ $B$ \
 is strongly consistent and is asymptotically normal using an appropriate random scaling;
 the weighted CLS estimator of \ $a$ \ is not weakly consistent, but asymptotically normal
 using an appropriate random scaling.
For more details, see Section \ref{section_HuangMaZhu}.
Huang et al.\ \cite{HuaMaZhu} also described the asymptotic behaviour of the weighted CLS
 estimator of \ $(B,a)$ \ in the subcritical and critical cases.

Overbeck and Ryd\'en \cite{OveRyd} investigated the CLS and weighted CLS estimators
 for a Cox--Ingersoll--Ross process in the subcritical case under the assumption that the initial value
 is distributed according to the unique stationary distribution.
A Cox--Ingersoll--Ross process is a CBI process of diffusion type (without jump part), i.e., when \ $\mu = 0$ \ and
 \ $\nu = 0$ \ in \eqref{SDE_atirasa_dim1}.
Based on discrete time observations \ $(X_k)_{k\in\{0,1,\ldots,n\}}$,
 \ $n \in \{1, 2, \ldots\}$, \  Overbeck and Ryd\'en \cite{OveRyd} derived a formula for the CLS estimator of
 \ $(B, a, c)$ \ and proved its asymptotic normality in the subcritical
 case.
Li and Ma \cite{LiMa} described the asymptotic behaviour of the CLS and weighted CLS estimators of \ $(B, a)$ \
 of an $\alpha$-stable CIR process given by the SDE \eqref{SDE_atirasa_dim1} with $c = 0$, $\nu=0$ and $\mu(\dd z) = z^{-1-\alpha}\,\dd z$ (where $\alpha\in(1,2)$)
 based on discrete time observations in the subcritical case under the assumption that $X_0$ is distributed
 according to the unique stationary distribution.

We turn back to the CBI process \ $(X_t)_{t\geq 0}$ \ given by the pathwise unique strong solution of the SDE \eqref{SDE_atirasa_dim1}.
Assuming that \ $X_0 = 0$, $c$, $\mu$ \ and \ $\nu$ \ are known, under eighth order moment conditions on the branching and immigration mechanisms,
 Barczy et al.\ \cite[Theorem 3.1]{BarKorPap} considered the (non-weighted) CLS estimator of \ $(B, A)$ \ based on discrete time observations
 \ $X_1,\ldots,X_n$, $n\in\{1,2,\ldots\}$, \ where
 \begin{align}\label{help1_BarKorPap}
  A := a + \int_0^\infty r \, \nu(\dd r).
 \end{align}
In the critical case (i.e., when $B=0$) they described its asymptotic behavior as the number of observations tends to infinity,
 provided that $A > 0$.
The limit distribution turns out to be non-normal except $C = 0$, where
 \begin{align}\label{help_C}
  C := 2 c + \int_0^\infty z^2 \, \mu(\dd z).
 \end{align}
In the critical case, if \ $C=0$, \ then the branching process in question is a pure immigration process,
 and the CLS estimator of \ $(B, A)$ \ is asymptotically normal.

Now we turn to summarize our new results on stable convergence of the appropriately normalized CLS estimators
 of (transformed) drift parameters for a CBI process.
In what follows, we suppose that \ $c$, \ $\mu$ \ and \ $\nu$ \ are known,
 and we consider the (non-weighted) CLS estimator of \ $(B, A)$ \ based on discrete time
 observations \ $(X_k)_{k\in\{0,1,\ldots,n\}}$, \ $n \in \{1, 2, \ldots\}$.
\ In the supercritical case (i.e., when \ $B > 0$) \ in case of a nontrivial immigration mechanism
 (equivalently, when \ $A > 0$) \ under second order moment conditions on the initial law and on the branching and immigration mechanisms,
 we describe the asymptotic behavior of the CLS estimator in question as the number of observations \ $n \to \infty$, \
 by proving stable convergence, see Theorem \ref{main}.
The limit distribution is mixed normal except the case \ $C = 0$.
In case of \ $C=0$, \ the appropriately scaled CLS estimator of \ $B$ \ converges stably to a random variable that can be written
 as an almost surely convergent infinite sum of some random variables (see \eqref{htBhtb1}),
 and the appropriately scaled CLS estimator of \ $A$ \ converges stably to a normally distributed random variable (see \eqref{htBhtb2}).
Our results immediately yield convergence in distribution of the appropriately normalized CLS estimators in question,
 since stable convergence implies convergence in distribution.
We again call the attention that in all the papers cited previously the authors proved convergence in distribution of
 (weighted) CLS estimators, and, according to our knowledge, results
 on stable convergence of the CLS estimator of drift parameters of CBI processes are not available in the literature.
Nonetheless, in Section \ref{section_HuangMaZhu}, we make a detailed comparison of our results on the asymptotic behaviour
 of the CLS estimator of $(B,A)$ and those of Huang et al.\ \cite{HuaMaZhu} on the asymptotic behaviour of
 the weighted CLS estimator of \ $(B,a)$.
\ We also compare our proof technique and that of Huang et al.\ \cite{HuaMaZhu}.

Next, we summarize our proof method.
The first and main step is to establish stable limit theorems for two martingales
 associated to the CBI process \ $(X_t)_{t\geq 0}$, \ see Theorems \ref{main_Ad} and \ref{main1_Ad}
 according to the cases $C>0$ and $C=0$.
It is worth noting here that there are several examples of martingales in the literature in connection with a number of branching processes,
 see, e.g., Athreya and Ney \cite[Chapter VI, Section 4]{AthNey} or Smorodina and Yarovaya \cite{SmoYar}.
The limit theorems in Theorems \ref{main_Ad} and \ref{main1_Ad} involve deterministic scaling, in case of $C>0$ the limit distribution is mixed normal,
 and in case of $C=0$, for one of the martingales is normal, and for the other martingale is given as the law
 of an absolutely convergent infinite sum of some random variables.
The proofs of Theorems \ref{main_Ad} and \ref{main1_Ad} are based on a general
 multidimensional stable limit theorem due to Barczy and Pap \cite[Theorem 1.4]{BarPap1}
 (see also Theorem \ref{MSLTES}), and a martingale central limit theorem involving mixing convergence
 due to H\"ausler and Luschgy \cite[Corollary 6.26]{HauLus} (see also Theorem \ref{Thm_HL_Cor6_26}).
Then, using Theorems \ref{main_Ad} and \ref{main1_Ad}, a suitable decomposition of the CLS estimator
 of \ $(\varrho , \cA) := (\ee^B, A \int_0^1 \ee^{Bs} \, \dd s)$, \ and Theorem 3.18 in H\"ausler and Luschgy \cite{HauLus}
 (see also Theorem \ref{Thm_HL_Thm3_18}, a continuous mapping theorem for stable convergence), we derive stable convergence of the CLS estimator of
 \ $(\varrho,\cA)$ \ in case of $C>0$, and those of the CLS estimator of \ $\varrho$ \ and the CLS estimator of \ $\cA$ \ in case of \ $C=0$, \ see Theorem \ref{main_rb}.
Finally, using a more general version of the continuous mapping theorem for stable convergence (see Lemma \ref{Lem_Kallenberg_stable})
 and Theorem \ref{main_rb}, we prove our main result Theorem \ref{main} on the stable convergence of the CLS estimator
 of \ $(B,A)$.

The paper is structured as follows.
In Section \ref{section_CBI} we recall some preliminaries on CBI processes
 such as the notions of branching and immigration mechanisms, and
 almost sure and \ $L_2$-asymptotic behaviour of supercritical CBI processes,
 see Theorems \ref{convX}, \ref{main_Ad0} and \ref{Thm_L2}.
In Section \ref{section_CBI_2} one can find a precise introduction of the  CLS estimator of \ $(\varrho,\cA)$ \
 based on the observations \ $X_0,X_1,\ldots,X_n$, \ and that of \ $(B,A)$, \ together with a result
 on their existence and uniqueness in the supercritical case (that we are interested in), see Lemma \ref{LEMMA_CLSE_exist_discrete}.
In the rest of Section \ref{section_CBI_2}, we present our results on stable convergence for the above mentioned CLS estimators
 and for the two associated martingales, see Theorems \ref{main}, \ref{main_rb}, \ref{main_Ad} and \ref{main1_Ad}.
Section \ref{section_HuangMaZhu} is devoted to a comparison of our results on the CLS estimator of \ $(B,A)$ \ and those of Huang et al.\ \cite{HuaMaZhu} on the weighted CLS estimator of \ $(B,a)$.
Sections \ref{section_proof_Theorem_main_Ad0}--\ref{section_proof_main1_Ad} contain the proof of our results,
 namely the proofs of Lemma \ref{LEMMA_CLSE_exist_discrete}, Theorems \ref{main_Ad0}, \ref{Thm_strong_consistency},
 \ref{main}, \ref{main_rb}, \ref{main_Ad} and \ref{main1_Ad}.
We close the paper with four appendices.
In Appendix \ref{App_1}, we recall the notions of stable and mixing convergences together with some known results
 that are used in the proofs.
We also recall a multidimensional stable limit theorem due to Barczy and Pap \cite{BarPap1}, which plays a crucial role
 in our proofs, see Theorem \ref{MSLTES}.
In Appendix \ref{section_SDE_moments}, we collect and (re)establish some results on SDE representation and moments of CBI processes.
In Appendix \ref{section_Lenglart}, we recall the Lenglart's inequality, and in Appendix \ref{CMT}, we present
 a version of the continuous mapping theorem for stable convergence.

Finally, we collect our notations.
Let \ $\ZZ_+$, \ $\NN$, \ $\RR$, \ $\RR_+$, \ $\RR_{++}$, \ $\CC$, \ $\CC_-$ \ and \ $\CC_{--}$ \ denote the set
 of non-negative integers, positive integers, real numbers, non-negative real
 numbers, positive real numbers, complex numbers, complex numbers with non-positive real parts and complex numbers with negative real parts,
 respectively.
The Borel \ $\sigma$-algebra on a subset \ $U$ \ of \ $\RR^d$ \ is denoted by \ $\cB(U)$, \ where \ $d\in\NN$,
 \ and recall that \ $\cB(U)=U\cap\cB(\RR^d)$.
\ By a Borel measure on a Borel set \ $S\in\cB(\RR^d)$, \ we mean a measure on \ $(S,\cB(S))$.
\ Further, let \ $\log^+(x):= \log(x)\bbone_{\{x\geq 1\}} + 0\cdot \bbone_{\{ 0\leq x < 1\}}$ \ for \ $x\in\RR_+$.
\ For \ $x , y \in \RR$, \ we will use the notations
 \ $x \land y := \min \{x, y\} $ \ and \ $x^+:= \max \{0, x\}$.
The integer part of \ $x\in\RR$ \ is denoted by \ $\lfloor x\rfloor$.
\ Throughout this paper, we make the conventions \ $\int_a^b := \int_{(a,b]}$
 \ and \ $\int_a^\infty := \int_{(a,\infty)}$ \ for any \ $a, b \in \RR$ \ with
 \ $a < b$.
\ By \ $\|\bx\|$ \ and \ $\|\bA\|$, \ we denote the Euclidean norm of a vector
 \ $\bx \in \RR^d$ \ and the induced matrix norm of a matrix
 \ $\bA \in \RR^{d\times d}$, \ respectively.
By \ $\langle \bx,\by\rangle$, \ we denote the Euclidean inner product of vectors \ $\bx,\by\in\RR^d$.
The null vector and the null matrix will be denoted by \ $\bzero$.
\ Moreover, \ $\bI_d \in \RR^{d\times d}$ \ denotes the $d\times d$ identity matrix,
 and \ $\be_1$, \ldots, $\be_d$ \ denotes the natural bases in \ $\RR^d$.
For $a_1,\ldots,a_d\in\RR$, we denote by $\diag(a_1,\ldots,a_d)$, the $d\times d$ diagonal matrix
 with diagonal entries $a_1,\ldots,a_d$.
For a symmetric and positive semidefinite matrix \ $\bA \in \RR^{d\times d}$, \
 its unique symmetric, positive semidefinite square root is denoted by \ $\bA^{1/2}$.
\ By \ $r(\bA)$, \ we denote the spectral radius of \ $\bA \in \RR^{d\times d}$.
\ By \ $C^2_\cc(\RR_+,\RR)$ \ we denote the set of twice continuously
 differentiable real-valued functions on \ $\RR_+$ \ with compact support.
Convergence almost surely, in probability, in \ $L_1$, \ in \ $L_2$ \ and
 in distribution under a probability measure \ $\PP$ \ will be denoted by \ $\asP$,
 \ $\stoch$, \ $\meanP$, \ $\qmeanP$ \ and \ $\distrP$, \ respectively.
For an event \ $A$ \ with \ $\PP(A) > 0$,
 \ let \ $\PP_A(\cdot) := \PP(\cdot\mid A) = \PP(\cdot \cap A) / \PP(A)$ \ denote
 the conditional probability measure given \ $A$.
\ Let \ $\EE_\PP$ \ denote the expectation under a probability measure \ $\PP$
 (in case of no confusion, we simple use the notation \ $\EE$).
\ Almost sure equality under a probability measure \ $\PP$ \ and equality in distribution
 will be denoted by \ $\aseP$ \ and \ $\distre$, \ respectively.
Every random variable will be defined on a suitable probability space \ $(\Omega,\cF,\PP)$.
\ For a random variable \ $\xi:\Omega\to\RR^d$, \ the distribution of \ $\xi$ \ on \ $(\RR^d,\cB(\RR^d))$ \ is denoted by \ $\PP^\xi$.
\ If \ $\bV \in \RR^{d\times d}$ \ is symmetric and positive semidefinite, then
 \ $\cN_d(\bzero, \bV)$ \ denotes the $d$-dimensional normal distribution with
 mean vector \ $\bzero\in\RR^d$ \ and covariance matrix \ $\bV$.
\ In case of \ $d=1$, \ instead of \ $\cN_1$ \ we simply write \ $\cN$ \ or, in some cases, \ $\cN(0,1)$.
\ The notions of stable and mixing convergences, some of their important properties used in the present paper
 and a multidimensional stable limit theorem are recalled in Appendix \ref{App_1}.

\section{Preliminaries on CBI processes}
\label{section_CBI}

CBI processes constitute an important class of Markov processes taking values in the nonnegative half line $\RR_+$.
These processes can be well-used for modelling the evolution of large populations with small individuals.
The study of continuous state and continuous time branching processes (without immigration) was initiated by Feller
 \cite{Fel}, who proved that a diffusion process may arise in a limit theorem for (discrete time)
 Galton-Watson branching processes.
CBI processes cover the situation where immigrants may come to the population from outer sources.
For a brief introduction to CBI processes, see Li \cite{Li3}.
CBI processes can be also viewed as special measure-valued branching Markov processes,
 for a modern book in this field, see Li \cite{Li}.
In this section, we recall the definition and some known results for CBI processes such as the notions of branching and immigration mechanisms,
 and almost sure and $L_2$-asymptotic behaviour of supercritical CBI processes.

\begin{Def}\label{Def_admissible}
A tuple \ $(c, a, b, \nu, \mu)$ \ is called a set of admissible
 parameters if \ $c, a \in \RR_+$, \ $b \in \RR$, \ and \ $\nu$ \ and
 \ $\mu$ \ are Borel measures on \ $\RR_{++}$ \ satisfying
 \ $\int_0^\infty (1 \land r) \, \nu(\dd r) < \infty$ \ and
 \ $\int_0^\infty (z \land z^2) \, \mu(\dd z)  < \infty$.
\end{Def}

Note that the measures \ $\nu$ \ and \ $\mu$ \ in Definition \ref{Def_admissible} are \ $\sigma$-finite,
 since the function \ $\RR_{++}\ni r\mapsto 1 \land r$ \ is strictly positive and measurable with a finite integral with respect to \ $\nu$,
 \ and the function \ $\RR_{++}\ni z \mapsto z \land z^2$ \ is strictly positive and measurable with a finite integral with respect to \ $\mu$,
 \ see, e.g., Kallenberg \cite[Lemma 1.4]{K2}.

\begin{Thm}\label{CBI_exists}
Let \ $(c, a, b, \nu, \mu)$ \ be a set of admissible parameters.
Then there exists a unique conservative transition semigroup \ $(P_t)_{t\in\RR_+}$ \ acting on
 the Banach space (endowed with the supremum norm) of real-valued bounded
 Borel measurable functions on the state space \ $\RR_+$ \ such that
 \begin{align*}
  \int_0^\infty \ee^{uy} P_t(x, \dd y)
  = \exp\biggl\{x \psi(t, u) + \int_0^t F(\psi(s, u)) \, \dd s\biggr\} ,
  \qquad x \in \RR_+, \quad u \in \CC_- , \quad t \in \RR_+ ,
 \end{align*}
 where, for any \ $u \in \CC_-$, \ the continuously differentiable
 function \ $\RR_+ \ni t \mapsto \psi(t, u) \in \CC_{--}$
 \ is the unique locally bounded solution to the differential
 equation
 \begin{equation*}
  \partial_t \psi(t, u) = R(\psi(t, u)) , \qquad
   \psi(0, u) = u ,
 \end{equation*}
 with
 \[
   R(u)
   := c u^2 +  b u
      + \int_0^\infty
         \bigl( \ee^{uz} - 1
                - u (1 \land z) \bigr)
         \, \mu(\dd z) , \qquad
   u \in \CC_- ,
 \]
 and
 \[
   F(u)
   := a u
      + \int_0^\infty
         \bigl(\ee^{ur} - 1\bigr)
         \, \nu(\dd r) , \qquad
   u \in \CC_- .
 \]
Further, the function \ $\RR_+\times \CC_{-}\ni (t,u)\mapsto \psi(t,u)$ \ is continuous.
\end{Thm}

Theorem \ref{CBI_exists} is a special case of Theorem 2.7  in Duffie et al.\ \cite{DufFilSch}
 with \ $m = 1$, \ $n = 0$ \ and zero killing rate.
The unique existence of a locally bounded solution to the differential equation
 in Theorem \ref{CBI_exists} is proved by Li \cite[page 48]{Li}.
The continuity of the function \ $\RR_+\times \CC_{-}\ni (t,u)\mapsto \psi(t,u)$ \
 follows from Proposition 6.4 in Duffie et al.\ \cite{DufFilSch}.

\begin{Def}\label{Def_CBI}
A conservative Markov process with state space \ $\RR_+$ \ and with transition semigroup
 \ $(P_t)_{t\in\RR_+}$ \ given in Theorem \ref{CBI_exists} is called a
 CBI process with parameters \ $(c, a, b, \nu, \mu)$.
\ The function \ $\CC_- \ni u \mapsto R(u) \in \CC$
 \ is called its branching mechanism, and the function
 \ $\CC_- \ni u \mapsto F(u) \in \CC_{-}$ \ is called its
 immigration mechanism.
A CBI process with parameters \ $(c, a, b, \nu, \mu)$ \ is called a CB process
 (a continuous state and continuous time branching process without immigration)
 if \ $a=0$ \ and \ $\nu=0$.
\end{Def}

Let \ $(X_t)_{t\in\RR_+}$ \ be a CBI process with parameters
 \ $(c, a, b, \nu, \mu)$ \ such that \ $\EE(X_0)<\infty$ \ and the moment condition
 \begin{equation}\label{moment_condition_m_nu}
  \int_1^\infty r \, \nu(\dd r) < \infty
 \end{equation}
 hold.
By Dawson and Li \cite[Theorems 5.1 and 5.2]{DawLi1} (or Barczy et al.\ \cite[Theorem 4.6]{BarLiPap2}),
 \ $(X_t)_{t\in\RR_+}$ \ can be represented as a pathwise unique strong solution of the SDE \eqref{SDE_atirasa_dim1},
 where
 \begin{equation}\label{tBbeta}
  B := b + \int_1^\infty (z - 1) \, \mu(\dd z).
 \end{equation}
Note that \ $B \in \RR$ \ due to \ $\int_0^\infty (z \land z^2) \, \mu(\dd z)  < \infty$
 and $0\leq z-1\leq z$ for $z\geq 1$.
\ Further, by formula (3.4) in Barczy et al. \cite{BarLiPap2} (see also formula (79) in Li \cite{Li3}), we get
 \begin{equation}\label{EXcond}
  \EE(X_t \mid X_0 = x)
  = \ee^{Bt} x + A \int_0^t \ee^{Bu} \, \dd u ,
  \qquad x \in \RR_+ , \quad t \in \RR_+ ,
 \end{equation}
 where $A$ is given in \eqref{help1_BarKorPap}.
Note that \ $A \in \RR_+$ \ due to \eqref{moment_condition_m_nu} and the moment condition on $\nu$ in Definition \ref{Def_admissible}.
Further, \ $\EE(X_t \mid X_0=x)$, $x\in\RR_+$, \ does not depend on the parameter \ $c$.
\ One can give probabilistic interpretations of the modified parameters
 \ $B$ \ and \ $A$, \ namely,
 \ $\ee^B = \EE(Y_1 \mid Y_0 = 1)$ \
 \ and \ $A = \EE(Z_1 \mid Z_0 = 0)$, \ where
 \ $(Y_t)_{t\in\RR_+}$ \ and \ $(Z_t)_{t\in\RR_+}$ \ are CBI
 processes with parameters \ $(c, 0, b, 0, \mu)$ \ and
 \ $(0, a, 0, \nu, 0)$, \ respectively, see formula
 \eqref{EXcond}.
The processes \ $(Y_t)_{t\in\RR_+}$ \ and \ $(Z_t)_{t\in\RR_+}$ \ can be
 considered as \emph{pure branching} (without immigration) and
 \emph{pure immigration} (without branching) processes, respectively.
Consequently, \ $\ee^B$ \ and \ $A$ \ may be called the branching
 and immigration mean, respectively.
Note that the branching mechanism depends only on the parameters \ $c$,
 $b$ \ and \ $\mu$, \ while the immigration mechanism depends only on the
 parameters \ $a$ \ and \ $\nu$.

Under the condition \eqref{moment_condition_m_nu} together with $\EE(X_0)<\infty$,
 \ by the help of the modified parameters \ $B$ \ and \ $A$, \ the SDE \eqref{SDE_atirasa_dim1} can be rewritten as
 \begin{align}\label{SDE_atirasa_dim1_mod}
  \begin{split}
   X_t
   &=X_0
     + \int_0^t (A + B X_s) \, \dd s
     + \int_0^t \sqrt{2 c X_s^+} \, \dd W_s \\
   &\quad
      + \int_0^t \int_0^\infty \int_0^\infty
         z \bbone_{\{u\leq X_{s-}\}} \, \tN(\dd s, \dd z, \dd u)
      + \int_0^t \int_0^\infty r \, \tM(\dd s, \dd r)
  \end{split}
 \end{align}
 for \ $t \in \RR_+$, \ where
 \ $\tM(\dd s, \dd r) := M(\dd s, \dd r) - \dd s \, \nu(\dd r)$.
\ Note that the SDE \eqref{SDE_atirasa_dim1_mod} does not contain integrals
 with respect to non-compensated Poisson random measures.

Next, we recall a convergence result for supercritical CBI processes, see Barczy et al.\ \cite[Theorem 3.1]{BarPalPap1}.

\begin{Thm}\label{convX}
Let \ $(X_t)_{t\in\RR_+}$ \ be a supercritical CBI process with
 parameters \ $(c, a, b, \nu, \mu)$ \ such that \ $\EE(X_0) < \infty$ \ and the
 moment condition \eqref{moment_condition_m_nu} hold.
Then there exists a non-negative random variable \ $w_{X_0}$ \ with \ $\EE(w_{X_0}) < \infty$ \ such that
 \begin{gather}\label{Conv_X_as}
  \ee^{-Bt} X_t
  \asP w_{X_0} \qquad
  \text{as \ $t \to \infty$.}
 \end{gather}
Further, if \ $\int_1^\infty z\log(z)\,\mu(\dd z)=\infty$, \ then \ $\PP(w_{X_0}=0)=1$.
\end{Thm}

Next, we formulate a consequence of Theorem \ref{convX}.

\begin{Thm}\label{main_Ad0}
Under the assumptions of Theorem \ref{convX}, for each \ $\ell \in \NN$, \ we have
 \begin{align}\label{help_main_Ad0_1}
   \ee^{-\ell Bn} \sum_{k=1}^n X_{k-1}^\ell \asP \frac{w_{X_0}^\ell}{\ee^{\ell B}-1}  \qquad
   \text{as \ $n \to \infty$,}
 \end{align}
 and
 \begin{align}\label{help_main_Ad0_2}
   \ee^{-2Bn} \sum_{k=1}^n X_{k-1}X_k \asP \frac{\ee^B}{\ee^{2B}-1}w_{X_0}^2 \qquad
   \text{as \ $n \to \infty$,}
 \end{align}
 where the random variable $w_{X_0}$ is given in Theorem \ref{convX}.
\end{Thm}

For a proof of Theorem \ref{main_Ad0}, see Section \ref{section_proof_Theorem_main_Ad0}.
Both Theorems \ref{convX} and \ref{main_Ad0} play important roles in the proofs.

Note that in case of \ $\PP(w_{X_0} = 0) = 0$, \ the scaling factor \ $\ee^{-Bt}$ \ is correct in Theorem \ref{convX}
 in the sense that the limit random variable $w_{X_0}$  is not zero almost surely.
In the next statement, a set of sufficient conditions is recalled for \ $\PP(w_{X_0} = 0) = 0$, \ see Barczy et al.\ \cite[Theorem 3.2]{BarPalPap2}.

\begin{Thm}\label{w=0}
Let \ $(X_t)_{t\in\RR_+}$ \ be a supercritical CBI process with
 parameters \ $(c, a, b, \nu, \mu)$ \ such that \ $\EE(X_0) < \infty$ \ and the
 moment condition
 \begin{equation}\label{moment_condition_2}
  \int_1^\infty z^2 \, \mu(\dd z) + \int_1^\infty r^2 \, \nu(\dd r) < \infty
 \end{equation}
 hold.
If \ $A > 0$, \ i.e., $a\ne 0$ \ or \ $\nu\ne 0$ \ (in other words $(X_t)_{t\in\RR_+}$ \ is not a pure branching process),
 then \ $\PP(w_{X_0} = 0) = 0$, \ where the random variable $w_{X_0}$ is given in Theorem \ref{convX}.
\end{Thm}

Note that if condition \eqref{moment_condition_2} holds, then we also have \ $\int_1^\infty z\log(z)\,\mu(\dd z)<\infty$.

Barczy et al.\ \cite[formula (3.57)]{BarPalPap1} showed that, under second order moment conditions,
 convergence in $L_2(\PP)$ also holds in \eqref{Conv_X_as}, that we recall below.

\begin{Thm}\label{Thm_L2}
Let \ $(X_t)_{t\in\RR_+}$ \ be a supercritical CBI process with parameters \ $(c, a, b, \nu, \mu)$ \
 such that \ $\EE(X_0^2) < \infty$ \ and the moment condition \eqref{moment_condition_2} hold.
Then \ $\EE(X_t^2)<\infty$, \ $t\in\RR_+$, \ $\EE(w_{X_0}^2)<\infty$, \ and \ $\ee^{-Bt} X_t \qmeanP w_{X_0}$ \ as \ $t\to\infty$.
\end{Thm}

\section{CLS estimators for supercritical CBI processes}
\label{section_CBI_2}

Let \ $(X_t)_{t\in\RR_+}$ \ be a CBI process with parameters
 \ $(c, a, b, \nu, \mu)$ \ such that \ $\EE(X_0)<\infty$ \ and the moment condition
 \eqref{moment_condition_m_nu} hold.
Recall that \ $(X_t)_{t\in\RR_+}$ \ is called subcritical, critical or supercritical
 if \ $B < 0$, \ $B = 0$ \ or \ $B > 0$, \ where \ $B$ \ is defined in \eqref{tBbeta}.

First, we derive CLS estimator of \ $(B,A)$, \ provided that \ $c$, \ $\nu$ \ and \ $\mu$ \ are known, where
 \ $A$ \ is given in \eqref{help1_BarKorPap}.
Taking into account the forms of \ $A$ \ and \ $B$, \ it equivalently means that we find CLS estimator of \ $(b,a)$ \ provided that
 \ $c$, \ $\nu$ \ and \ $\mu$ \ are known.
In order to derive CLS estimator of \ $(B,A)$, \ first we derive CLS estimator of the modified parameters \ $\ee^B$ \
 and \ $A\int_0^1 \ee^{Bs}\,\dd s$, \ and then using it, we get the CLS estimator of \ $(B,A)$.
\ We follow this procedure, since the given below minimization problem \eqref{min_problem} corresponding to the CLS estimator of the modified
 parameters in question is a kind of linear one in the modified parameters, but not in the parameters \ $B$ \ and \ $A$.

For each \ $k \in \ZZ_+$, \ let \ $\cF_k^X := \sigma(X_0, X_1 , \dots, X_k)$,
 \ and let \ $\cF_\infty^X := \sigma(\bigcup_{k=0}^\infty \cF_k^X)$.
We suppose that the probability space \ $(\Omega,\cF_\infty^X ,\PP)$ \ is complete,
 i.e., for all \ $S \in \cF_\infty^X$ \ with \ $\PP(S) = 0$ \ and for all subsets \ $U\subset S$,
 \ we have \ $U\in\cF_\infty^X$.

Since \ $(X_t)_{t\in\RR_+}$ \ is a time-homogeneous Markov process, by
 \eqref{EXcond}, we have
 \begin{equation*}
  \EE(X_k \mid \cF_{k-1}^X) = \EE(X_k \mid X_{k-1})
  = \varrho X_{k-1} + \cA ,
  \qquad k \in \NN ,
 \end{equation*}
 where
 \begin{equation}\label{ttBbeta}
  \varrho := \ee^B \in \RR_{++} , \qquad
  \cA := A \int_0^1 \ee^{Bs} \, \dd s \in \RR_+ .
  \end{equation}
Note that \ $\cA = \EE(X_1 \mid X_0 = 0)$, \ see \eqref{EXcond}.
Note also that \ $\cA$ \ depends both on the branching and immigration
 mechanisms, although \ $A$ \ depends only on the immigration mechanism.
Let us introduce the sequence
 \begin{equation}\label{Mk}
  M_k
  := X_k - \EE(X_k \mid \cF_{k-1}^X)
   = X_k - \varrho X_{k-1} - \cA ,
  \qquad k \in \NN ,
 \end{equation}
 of martingale differences with respect to the filtration
 \ $(\cF_k^X)_{k \in \ZZ_+}$.
\ By \eqref{Mk}, the process \ $(X_k)_{k \in \ZZ_+}$ \ satisfies the recursion
 \begin{equation}\label{regr}
  X_k = \varrho X_{k-1} + \cA + M_k ,
  \qquad k \in \NN .
 \end{equation}
Supposing that \ $c$, \ $\nu$ \ and \ $\mu$ \ are known, for each \ $n \in \NN$, \ a CLS estimator \ $(\hvarrho_n, \hcA_n)$ \ of
 \ $(\varrho, \cA)$ \ based on the observations \ $X_0, X_1, \ldots, X_n$ \ can be
 obtained by minimizing the sum of squares
 \begin{align}\label{min_problem}
   \sum_{k=1}^n (X_k - \varrho X_{k-1} - \cA)^2
 \end{align}
 with respect to \ $(\varrho, \cA)$ \ over \ $\RR^2$, \ and one can easily check that it has the form
 \begin{equation}\label{CLSErb}
  \begin{bmatrix}
   \hvarrho_n \\[1mm]
   \hcA_n
  \end{bmatrix}
  := \frac{1}{n \sum\limits_{k=1}^n X_{k-1}^2
              - \left( \sum\limits_{k=1}^n X_{k-1} \right)^2}
     \begin{bmatrix}
      n \sum\limits_{k=1}^n X_k X_{k-1}
      - \sum\limits_{k=1}^n X_k \sum\limits_{k=1}^n X_{k-1} \\[3mm]
      \sum\limits_{k=1}^n X_k \sum\limits_{k=1}^n X_{k-1}^2
      - \sum\limits_{k=1}^n X_k X_{k-1} \sum\limits_{k=1}^n X_{k-1}
     \end{bmatrix}
 \end{equation}
 on the set
 \begin{gather*}
  H_n:=\Biggl\{\omega \in \Omega
               : n \sum_{k=1}^n X_{k-1}^2(\omega)
                 - \left( \sum_{k=1}^n X_{k-1}(\omega) \right)^2 > 0\Biggr\} .
 \end{gather*}
One can realize that the particular form of the CLS estimator \eqref{CLSErb} also appears as a CLS estimator of some parameters
 for other branching models, see, e.g., Wei and Winnicki \cite[formulas (1.4), (1.5)]{WW}.

In the sequel we investigate the supercritical case, i.e., the case $B>0$.

\begin{Lem}\label{LEMMA_CLSE_exist_discrete}
Let $(X_t)_{t\in\RR_+}$ be a supercritical CBI process with parameters
 $(c, a, b, \nu, \mu)$ such that $\EE(X_0) < \infty$ and the moment condition \eqref{moment_condition_2} hold.
Suppose that $A>0$, where $A$ is given in \eqref{help1_BarKorPap}, i.e., $(X_t)_{t\in\RR_+}$ is not a pure branching process.
Then we have $\PP(w_{X_0}>0)=1$ and the convergence $\PP(H_n) \to 1$ as $n \to \infty$ holds,
 where the random variable $w_{X_0}$ is given in Theorem \ref{convX}.
Hence the probability of the existence of a unique CLS estimator $(\hvarrho_n, \hcA_n)$ of $(\varrho,\cA)$
 converges to 1 as $n \to \infty$, and this CLS estimator has the form given in \eqref{CLSErb}
 on the event $H_n$.
\end{Lem}

The proof of Lemma \ref{LEMMA_CLSE_exist_discrete} can be found in Section \ref{section_proof_LEMMA_CLSE_exist_discrete}.

In what follows, we suppose that $(X_t)_{t\in\RR_+}$ is a supercritical CBI process with parameters
 \ $(c,a,b,\nu,\mu)$ \ such that \ $\EE(X_0) < \infty$, \ the moment condition
 \eqref{moment_condition_2} and \ $A>0$ \ (i.e., \ $(X_t)_{t\in\RR_+}$ \ is not a pure branching process) hold.

Next, we introduce a natural estimator of \ $(B,A)$ \ based on the observations \ $X_0,X_1,\ldots,X_n$. \
Let us define the function
 \ $h : \RR^2 \to \RR_{++} \times \RR$ \ given by
 \[
   h(x, y)
   := \biggl( \ee^x, y \int_0^1 \ee^{xs} \, \dd s \biggr), \qquad (x, y) \in \RR^2 .
 \]
Note that \ $h$ \ is bijective having inverse \ $h^{-1}:\RR_{++}\times \RR\to\RR^2$ \ given by
 \[
   h^{-1}(u, v)
   = \left(\log(u),
           \frac{v}{\int_0^1 u^s \, \dd s} \right), \qquad (u, v) \in \RR_{++} \times \RR .
 \]
Our forthcoming Theorem \ref{Thm_strong_consistency} implies that the CLS estimator
 \ $\hvarrho_n$ \ of \ $\varrho=\ee^{B}>1$ \ is strongly and hence weakly consistent,
 and so, for all \ $\varepsilon,\delta\in(0,1)$ \ there exists an \ $n_0\in\NN$ \ such that
 \ $\PP(\vert \hvarrho_n - \varrho\vert\leq \vare) > 1-\delta$ \ for each $n\geq n_0$, \
 yielding that
 \begin{align}\label{help_est_existence_1}
   \PP(\hvarrho_n \in\RR_{++})\geq \PP(\hvarrho_n \geq  -\vare + \varrho) \geq \PP(\vert \hvarrho_n - \varrho\vert\leq \vare) > 1-\delta
    \qquad \text{for each \ $n\geq n_0$.}
 \end{align}
Consequently, for all \ $\delta\in(0,1)$ \ there exists an \ $n_0\in\NN$ such that
 \ $\PP\big( \hvarrho_n\in\RR_{++}\big)  > 1-\delta$ \ for each \ $n\geq n_0$,
 \ and hence for each \ $n\geq n_0$, \ the event
 \[
     \left\{  (\hvarrho_n, \hcA_n)
      = \argmin_{(\varrho,\cA)\in\RR_{++} \times \RR}\;
       \sum_{k=1}^n (X_k - \varrho X_{k-1} - \cA)^2 \right\}
 \]
 has a probability greater than $1-\delta$.
Thus, motivated by that \ $(B,A)=h^{-1}(\varrho,\cA)$, \ one can introduce
 a natural estimator of \ $(B, A)$ \ based on the observations \ $X_0, X_1, \ldots, X_n$ \
 by applying the inverse of \ $h$ \ to the CLS estimator of
 \ $(\varrho, \cA)$, \ that is,
 \begin{align}\label{help_B_A_est_definition}
   (\hB_n, \hA_n)
   := h^{-1}(\hvarrho_n, \hcA_n)
    = \left(\log(\hvarrho_n),
             \frac{\hcA_n}{\int_0^1 (\hvarrho_n)^s \, \dd s} \right),
   \qquad \ n \in \NN ,
  \end{align}
 on the set
 \ $\{\omega \in \Omega : \hvarrho_n(\omega) \in \RR_{++}\}$.
\ Due to \eqref{help_est_existence_1} and Lemma \ref{LEMMA_CLSE_exist_discrete}, the probability of the existence of the estimator
 \ $(\hB_n, \hA_n)$ \ converges to \ $1$ \ as \ $n\to\infty$.
\ Further, using that \ $h$ \ is bijective, we also obtain that
 \begin{equation*}
      (\hB_n, \hA_n)
      = \argmin_{(B,A)\in\RR^2}\;
        \sum_{k=1}^n
        \left( X_k - \ee^B X_{k-1}
                - A \int_0^1 \ee^{Bs} \, \dd s\right)^2
 \end{equation*}
 on the set \ $\{\omega\in\Omega : \hvarrho_n(\omega)\in\RR_{++} \}$.
\ In view of this, in what follows, we will simply call $(\hB_n, \hA_n)$ the CLS estimator of $(B,A)$
 based on the observations $X_0,X_1,\ldots,X_n$.
\ We would like to stress the point that the probability of the existence of
 the CLS estimator \ $\bigl(\hB_n, \hA_n\bigr)$ \ in general is not \ $1$, \
 it only converges to \ $1$ \ as \ $n\to\infty$.
\ However, this does not cause any problem in limit theorems with convergence in distribution
 or with stable convergence following from the next observations.
Namely, if \ $\xi_n$, \ $\eta_n$, \ $n\in\NN$, \ and \ $\xi$ \ are random variables on a probability space
 \ $(\Omega,\cF,\PP)$ \ such that \ $\xi_n\distrP \xi$ \ as \ $n\to\infty$ \
 and \ $\lim_{n\to\infty}\PP(\xi_n=\eta_n)=1$, \ then \ $\eta_n\distrP \xi$ \ as \ $n\to\infty$.
\ Indeed, for all $\vare>0$, we have $\PP(\vert \eta_n-\xi_n\vert \geq \vare) \leq \PP(\eta_n\ne \xi_n)\to 0$
 \ as \ $n\to\infty$, i.e., $\eta_n-\xi_n\stoch 0$ \ as \ $n\to\infty$, \ and then one can use
 Slutsky's lemma.
Similarly, if \ $\xi_n$, \ $\eta_n$, \ $n\in\NN$, \ and \ $\xi$ \ are random variables on a probability space
 \ $(\Omega,\cF,\PP)$ \ and \ $\cG\subset \cF$ \ is a sub-$\sigma$-field such that \ $\xi_n$ \ converges
 \ $\cG$-stably to \ $\xi$ \ as \ $n\to\infty$ \ and \ $\lim_{n\to\infty}\PP(\xi_n=\eta_n)=1$, \ then
 \ $\eta_n$ \ converges \ $\cG$-stably to \ $\xi$ \ as \ $n\to\infty$.
\ Indeed, as we have seen earlier, \ $\eta_n-\xi_n\stoch 0$ \ as \ $n\to\infty$, \ and then one can use part (a) of Theorem 3.18
 in H\"ausler and Luschgy \cite{HauLus} (see also part (a) of Theorem \ref{Thm_HL_Thm3_18}).

\begin{Thm}\label{main}
Let \ $(X_t)_{t\in\RR_+}$ \ be a supercritical CBI process with parameters
 \ $(c, a, b, \nu, \mu)$ \ such that \ $\EE(X_0^2) < \infty$, \ the
 moment condition \eqref{moment_condition_2} and $A>0$ (i.e., $(X_t)_{t\in\RR_+}$ is not a pure branching process) hold.
Then the probability of the existence of the CLS estimator \ $(\hB_n, \hA_n)$ \ of \ $(B,A)$ \
 converges to 1 as \ $n \to \infty$, \ and
 \begin{equation}\label{htBhtb}
  \begin{bmatrix}
   \ee^{Bn/2} (\hB_n - B) \\
   n \ee^{-Bn/2} (\hA_n - A)
  \end{bmatrix}
  \to \bS^{1/2} \bN
  \qquad \text{$\cF_\infty^X$-stably as \ $n \to \infty$,}
 \end{equation}
 where \ $\bN$ \ is a 2-dimensional random vector $\PP$-independent of \ $\cF_\infty^X$
  \ such that \ $\bN \distre \cN_2(\bzero, C \bI_2)$, \ where \ $C$ \ is given in \eqref{help_C}
  and the random matrix \ $\bS$, \ defined by
 \[
   \bS
   := \begin{bmatrix}
       \frac{(\ee^B-1)(\ee^{2B}-1)^2}{B\ee^B(\ee^{3B}-1)} w_{X_0}^{-1}
        & - \frac{\ee^B(\ee^B-1)}{\ee^{3B}-1} \\[1mm]
       - \frac{\ee^B(\ee^B-1)}{\ee^{3B}-1}
        & \frac{B\ee^{2B}}{(\ee^B-1)(\ee^{3B}-1)} w_{X_0}
      \end{bmatrix},
 \]
 is \ $\PP$-independent of \ $\bN$, \ where the random variable \ $w_{X_0}$ \ is given in Theorem \ref{convX}.

If, in addition, \ $C = 0$\ (i.e., \ $c = 0$ \ and \ $\mu = 0$), \ then
 \begin{gather}\label{htBhtb1}
  \ee^{Bn} (\hB_n - B)\to \frac{\ee^{2B}-1}{\ee^{2B}} w_{X_0}^{-1}
      \sum_{j=0}^\infty \ee^{-Bj} Z_j
  \qquad \text{$\cF_\infty^X$-stably as \ $n \to \infty$,} \\
  n^{1/2} (\hA_n - A)
  \to N_1
  \qquad \text{$\cF_\infty^X$-mixing as \ $n \to \infty$,} \label{htBhtb2}
 \end{gather}
 where $N_1$ is a random variable $\PP$-independent of $\cF_\infty^X$ such that
 \[
    N_1 \distre \cN\left(0, \frac{B(\ee^{2B}-1)}{2(\ee^B-1)^2} \int_0^\infty r^2 \, \nu(\dd r)\right),
 \]
 and \ $(Z_j)_{j\in\ZZ_+}$ \  are \ $\PP$-independent and identically distributed random variables being \ $\PP$-independent  of
 \ $\cF_\infty^X$ \ such that \ $Z_1$ \ has a characteristic function
 \[
         \EE(\ee^{\ii\theta Z_1})
          = \exp\biggl\{\int_0^1 \int_0^\infty
                    \Big(\ee^{\ii\theta r \ee^{Bu} }- 1  - \ii \theta r \ee^{Bu} \Big)
                    \, \dd u \, \nu(\dd r)\biggr\},
                    \qquad \theta\in\RR.
 \]
In particular, we have $Z_1 \distre M_1$.
Further, the series in \eqref{htBhtb1} is absolutely convergent $\PP$-almost surely.
\end{Thm}

The proof of Theorem \ref{main} can be found in Section \ref{section_proof_main}.

\begin{Rem}
In Theorem \ref{main} the random variable \ $\bN$ \ and the random matrix \ $\bS$ \ are indeed \ $\PP$-independent,
 since \ $\bN$ \ is \ $\PP$-independent of \ $\cF_\infty^X$ \ and \ $\bS$ \ is \ $\cF_\infty^X$-measurable
 (following from the \ $\cF_\infty^X$-measurability of \ $w_{X_0}$, \ which is explained in Step 1 of the proof of Theorem \ref{main_Ad}).
The $\cF_\infty^X$-measurability of $w_{X_0}$ also implies that in \eqref{htBhtb1} the sequence of random variables
 \ $(Z_j)_{j\in\ZZ_+}$ \ and the random variable \ $w_{X_0}$ \ are \ $\PP$-independent.
The random variable \ $w_{X_0}$ \ is infinitely divisible, since \ $\ee^{-Bt}X_t\asP w_{X_0}$ \ as \ $t\to\infty$ \ (see Theorem \ref{convX}),
 \ $\ee^{-Bt}X_t$ \ is infinitely divisible for all \ $t\in\RR_{++}$ \
 (following, e.g., from Li \cite[Section 3.3]{Li}) and it is known that if a sequence of infinitely divisible random variables converges
  in distribution, then the limit random variable is infinitely divisible as well.
Note also that the laws of the limit random variables in (3.8) and (3.9) may depend on the law of the initial value \ $X_0$,
 \ since the law of the random variable $w_{X_0}$ \ (which appears in the limit random variables in question)
 may depend on the law of \ $X_0$. \ This phenomenon usually happens for limit laws of CLS estimators for
 supercritical branching or affine models using deterministic scalings.
\proofend
\end{Rem}

In the next two remarks, we discuss the case \ $C=0$ \ in Theorem \ref{main}.

\begin{Rem}\label{REMARK1}
In case of \ $C = 0$, the convergence \eqref{htBhtb} implies that \ $\ee^{Bn/2} (\hB_n - B) \distrP 0$ \
 and \ $n \ee^{-Bn/2} (\hA_n - A) \distrP 0$ \ as \ $n \to \infty$, \ since stable convergence yields
 convergence in distribution.
Hence in case of \ $C=0$ the scaling factors in \eqref{htBhtb} are not the good ones in the sense that
 the limit distribution in \eqref{htBhtb} is the 2-dimensional zero vector.
This motivates establishing the convergences \eqref{htBhtb1} and \eqref{htBhtb2}, which show that
 in case of \ $C=0$ \ and \ $\nu\ne 0$, \
 the scaling factors \ $\ee^{Bn}$ \ and \ $n^{1/2}$ \ are
 the good ones for \ $\hB_n - B$ \ and \ $\hA_n - A$, \ respectively,
 in the sense that the corresponding limit distributions are not identically zero.
Indeed, if \ $\nu\ne 0$, \ then \ $\int_0^\infty r^2\,\nu(\dd r)\ne 0$ \ yielding that the variance
 of the normally distributed random variable \ $N_1$ \ in \eqref{htBhtb2} is not \ $0$.
\ Further, in case of \ $\nu\ne 0$, \ we check that \ $Z_1$ \ is not identically $0$.
\ If $\nu \ne 0$, then there exists a $G\in\cB((0,\infty))$ such that $\nu(G)>0$, and then
 for all \ $\theta\in\RR\setminus\{0\}$, we have
 \[
   \int_0^1 \int_0^\infty
    \Big( \cos(\theta r \ee^{Bu}) - 1 \Big) \, \dd u \, \nu(\dd r)
      \leq
     \int_G \left(\int_0^1 \Big( \cos(\theta r \ee^{Bu}) - 1 \Big) \, \dd u \right) \nu(\dd r) ,
 \]
 where the right hand side is strictly negative, since for all $\theta\in\RR\setminus\{0\}$ and $r\in G$, the set
 \[
  \big\{ u\in[0,1] : \cos(\theta r \ee^{Bu}) \ne 1 \big\}
     = \Big\{ u\in[0,1] : \frac{1}{2\pi}\theta r \ee^{Bu} \notin \ZZ \Big\}
 \]
 has Lebesgue measure $1$.
This yields that the absolute value of \ $\EE(\ee^{\ii\theta Z_1})$ \ is strictly less than $1$ for all $\theta\in\RR\setminus\{0\}$
 implying that $Z_1$ cannot be the identically $0$ random variable.
Consequently, the absolute value of \ $\EE(\ee^{\ii\theta \sum_{j=0}^\infty\ee^{-Bj} Z_j})$ \
 is strictly less than $1$ for all $\theta\in\RR\setminus\{0\}$, yielding that $\sum_{j=0}^\infty\ee^{-Bj} Z_j$ cannot be the identically 0 random variable.
Indeed, using that $Z_j$, $j\in\ZZ_+$, are $\PP$-independent and identically distributed, we have
 \begin{align*}
   \left\vert \EE(\ee^{\ii\theta \sum_{j=0}^\infty\ee^{-Bj} Z_j}) \right\vert
    &=\left\vert \prod_{j=0}^\infty \exp\left\{ \int_0^1 \int_0^\infty \Big( \ee^{\ii \theta \ee^{-Bj} r \ee^{Bu}} - 1 - \ii \theta \ee^{-Bj}r \ee^{Bu} \Big)\dd u\,\nu(\dd r)\right\} \right\vert
 \end{align*}
 \begin{align*}
  &= \prod_{j=0}^\infty \exp\left\{ \int_0^1 \int_0^\infty \Big( \cos(\theta \ee^{-Bj} r \ee^{Bu}) - 1\Big)\dd u\,\nu(\dd r)\right\}\\
    &=\exp\left\{ \int_0^1 \int_0^\infty \sum_{j=0}^\infty \Big( \cos(\theta \ee^{-Bj} r \ee^{Bu}) - 1\Big)\,\dd u\,\nu(\dd r) \right\},
    \qquad \theta\in\RR.
 \end{align*}
\proofend
\end{Rem}

\begin{Rem}
Under the assumptions of Theorem \ref{CBI_exists} together with \ $C=0$, \ the SDE \eqref{SDE_atirasa_dim1_mod} takes the form
  \begin{align*}
   X_t =X_0 + \int_0^t (A + b X_s) \, \dd s
         + \int_0^t \int_0^\infty r \, \tM(\dd s, \dd r), \qquad t\in\RR_+,
 \end{align*}
and in this case we have
 \[
  \int_0^\infty \ee^{uy} P_t(x, \dd y)
  = \exp\biggl\{x u\ee^{tb} + \int_0^t F(u\ee^{sb}) \, \dd s\biggr\} ,
     \qquad x \in \RR_+, \quad u \in \CC_- , \quad t \in \RR_+ ,
 \]
 since the unique locally bounded solution to the differential equation \ $\partial_t \psi(t,u) = b \psi(t,u)$ \ with initial
 condition \ $\psi(0,u) = u$, \ takes the form \ $\psi(t,u) = u \ee^{tb}$, \ $t\in\RR_+$.
\proofend
\end{Rem}

\begin{Rem}\label{REMARK2}
Clearly, \eqref{htBhtb} and \eqref{htBhtb1} imply that the estimator $\hB_n$ is weakly consistent as $n\to\infty$.
In fact, we also have that the estimator $\hB_n$ is strongly consistent as $n\to\infty$, see our forthcoming
 Theorem \ref{Thm_strong_consistency}.
Moreover, \eqref{htBhtb} and \eqref{htBhtb2} yield that the estimator $\hA_n$ is not (weakly) consistent as $n\to\infty$
 in case of $C \ne 0$, and consistent if $C = 0$.
Indeed, in case of $C\ne 0$, if $\hA_n$ were consistent as $n\to\infty$, then
 $n \ee^{-Bn/2} (\hA_n - A)$ would converge in probability under $\PP$ to $0$ as $n\to\infty$,
 which cannot happen in view of \eqref{htBhtb}, since $C\ne0$ and $\PP(w_{X_0}=0)=0$ (following from $A>0$ by Theorem \ref{w=0}).
Note that the same phenomenon occurs for the weighted LSE of $B$ and $a$ proved by Huang et al.\ \cite{HuaMaZhu}
 (for more details, see Section \ref{section_HuangMaZhu}).
In case of $C=0$, the question of strong consistency of $\hA_n$ as $n\to\infty$ remains open.
\proofend
\end{Rem}

\begin{Thm}\label{Thm_strong_consistency}
Let \ $(X_t)_{t\in\RR_+}$ \ be a supercritical CBI process with parameters
 \ $(c, a, b, \nu, \mu)$ \ such that \ $\EE(X_0) < \infty$, \ the
 moment condition \eqref{moment_condition_2} and $A>0$ \ (i.e., $(X_t)_{t\in\RR_+}$ is not a pure branching process) hold.
Then the CLS estimator $\hB_n$ of $B$, and the CLS estimator \ $\hvarrho_n$ \ of \ $\varrho$ \ are strongly consistent,
 i.e., \ $\hB_n\asP B$ \ as $n\to\infty$, \ and \ $\hvarrho_n \asP \varrho$ \ as \ $n\to\infty$, \ respectively.
\end{Thm}

The proof of Theorem \ref{Thm_strong_consistency} can be found in Section \ref{section_proof_Thm_strong_consistency}.
Theorem \ref{main} will follow from the following statement.

\begin{Thm}\label{main_rb}
Under the assumptions of Theorem \ref{main}, the probability of the existence
 of a unique CLS estimator \ $(\hvarrho_n, \hcA_n)$ \ of \ $(\varrho,\cA)$ \ converges to 1
 as \ $n \to \infty$ \ and
 \begin{equation}\label{hvarrhohobeta}
  \begin{bmatrix}
   \ee^{Bn/2} (\hvarrho_n - \varrho) \\[1mm]
   n \ee^{-Bn/2} (\hcA_n - \cA)
  \end{bmatrix}
  \to \tbS^{1/2} \tbN
  \qquad \text{$\cF_\infty^X$-stably as \ $n \to \infty$,}
 \end{equation}
 where \ $\tbN$ \ is a 2-dimensional random vector $\PP$-independent of
  \ $\cF_\infty^X$ \ such that \ $\tbN \distre \cN_2(\bzero, V \bI_2)$ \ with \ $V = C \int_0^1 \ee^{B(1+u)} \, \dd u$ \
 (appearing in Proposition \ref{moment_formula_2}), and the random matrix \ $\tbS$, \ defined by
 \[
   \tbS
   := \begin{bmatrix}
       \frac{(\ee^{2B}-1)^2}{\ee^{3B}-1} w_{X_0}^{-1}
        & - \frac{\ee^B(\ee^B-1)}{\ee^{3B}-1} \\[1mm]
       - \frac{\ee^B(\ee^B-1)}{\ee^{3B}-1}
        & \frac{\ee^B}{\ee^{3B}-1} w_{X_0}
      \end{bmatrix},
 \]
 is \ $\PP$-independent of \ $\tbN$.

If, in addition, \ $C = 0$ \ (i.e., $c=0$ and $\mu=0$), \ then
 \begin{gather}\label{hvarrhohobeta1}
  \ee^{Bn} (\hvarrho_n - \varrho)
  \to \frac{\ee^{2B}-1}{\ee^B} w_{X_0}^{-1}
      \sum_{j=0}^\infty \ee^{-Bj} Z_j
  \qquad \text{$\cF_\infty^X$-stably as \ $n \to \infty$,} \\
  n^{1/2} (\hcA_n - \cA)
  \to \tN_1
  \qquad \text{$\cF_\infty^X$-mixing as \ $n \to \infty$,} \label{hvarrhohobeta2}
 \end{gather}
 where \ $\tN_1$ \ is a random variable $\PP$-independent of $\cF_\infty^X$
 such that \ $\tN_1 \distre \cN\bigl(0, \frac{\ee^{2B}-1}{2B} \int_0^\infty r^2 \, \nu(\dd r)\bigr)$,
 and \ $(Z_j)_{j\in\ZZ_+}$ \  are \ $\PP$-independent and identically distributed random variables being \ $\PP$-independent from
 \ $\cF_\infty^X$ \ such that \ $Z_1$ \ has a characteristic function
 \[
         \EE(\ee^{\ii\theta Z_1})
          = \exp\biggl\{\int_0^1 \int_0^\infty
                    \Big(\ee^{\ii\theta r \ee^{Bu} }- 1  - \ii \theta r \ee^{Bu} \Big)
                    \, \dd u \, \nu(\dd r)\biggr\},
                    \qquad \theta\in\RR.
 \]
In particular, we have \ $Z_1 \distre M_1$.
\ Further, the series in \eqref{hvarrhohobeta1} is absolutely convergent $\PP$-almost surely.
\end{Thm}

The proof of Theorem \ref{main_rb} can be found in Section \ref{section_proof_main_rb}.

\begin{Rem}
In Theorem \ref{main_rb} the random variable \ $\tbN$ \ and the random matrix \ $\tbS$ \ are indeed \ $\PP$-independent,
 since \ $\tbN$ \ is \ $\PP$-independent of \ $\cF_\infty^X$ \ and \ $\tbS$ \ is \ $\cF_\infty^X$-measurable
 (following from the \ $\cF_\infty^X$-measurability of \ $w_{X_0}$, \ which is explained in Step 1 of the proof of Theorem \ref{main_Ad}).
The \ $\cF_\infty^X$-measurability of \ $w_{X_0}$ also implies that in \eqref{hvarrhohobeta1} the sequence of random variables \ $(Z_j)_{j\in\ZZ_+}$ \ and the random variable
 \ $w_{X_0}$ \ are \ $\PP$-independent.
\proofend
\end{Rem}

\begin{Rem}\label{REMARK3}
Clearly, \eqref{hvarrhohobeta} and \eqref{hvarrhohobeta1} imply that the estimator \ $\hvarrho_n$ \ is weakly consistent as $n\to\infty$.
In fact, we also have that the estimator \ $\hvarrho_n$ \ is strongly consistent as \ $n\to\infty$, \ see Theorem \ref{Thm_strong_consistency}.
Moreover, \eqref{hvarrhohobeta} and \eqref{hvarrhohobeta2} yield that the estimator \ $\hcA_n$ \ is not (weakly) consistent
 in case of \ $C \ne 0$, \ and consistent if \ $C = 0$.
\proofend
\end{Rem}

Theorem \ref{main_rb} will follow from Theorem \ref{main_Ad0} and the forthcoming Theorems \ref{main_Ad} and \ref{main1_Ad}.

\begin{Thm}\label{main_Ad}
Under the assumptions of Theorem \ref{main}, we have
 \begin{equation}\label{U}
  \begin{split}
   \begin{bmatrix}
      \ee^{-Bn/2} & 0 \\
      0 & \ee^{-3Bn/2} \\
    \end{bmatrix}
   \begin{bmatrix}
   \sum_{k=1}^n M_k \\
   \sum_{k=1}^n M_k X_{k-1}
   \end{bmatrix}
  \to \bR^{1/2} \tbN
  \qquad \text{$\cF_\infty^X$-stably as \ $n \to \infty$,}
  \end{split}
 \end{equation}
 where \ $\tbN$ \ is a 2-dimensional random vector $\PP$-independent of \ $\cF_\infty^X$ \ such that \ $\tbN \distre \cN_2(\bzero, V \bI_2)$ \
  with \ $V = C \int_0^1 \ee^{B(1+u)} \, \dd u$ \ (appearing in Proposition \ref{moment_formula_2}), and the random matrix \ $\bR$, \ defined by
 \[
   \bR
   := \begin{bmatrix}
       \frac{w_{X_0}}{\ee^B-1}
        & \frac{w_{X_0}^2}{\ee^{2B}-1} \\[1mm]
       \frac{w_{X_0}^2}{\ee^{2B}-1}
        & \frac{w_{X_0}^3}{\ee^{3B}-1}
      \end{bmatrix},
 \]
 is \ $\PP$-independent of \ $\tbN$.
\end{Thm}

The proof of Theorem \ref{main_Ad} can be found in Section \ref{section_proof_main_Ad}.

\begin{Rem}\label{REMARK4}
(i): Note that the random matrix \ $\bR$ \ given in Theorem \ref{main_Ad} is symmetric and positive definite \ $\PP$-almost surely,
 and hence its unique symmetric and positive definite square root \ $\bR^{1/2}$ \ is well-defined \ $\PP$-almost surely.
Indeed, we have \ $w_{X_0}/(\ee^B-1)$ \ is positive \ $\PP$-almost surely, and
 \begin{align*}
  \det(\bR)
   = w_{X_0}^4\left(\frac{1}{(\ee^{B}-1)(\ee^{3B}-1)} - \frac{1}{(\ee^{2B}-1)^2}\right)
   = w_{X_0}^4\frac{\ee^B(\ee^B-1)}{(\ee^{2B}-1)^2(\ee^{3B}-1)}
   > 0
 \end{align*}
 $\PP$-almost surely, where we used that \ $\PP(w_{X_0}>0)=1$ \ (following from Theorem \ref{w=0}).
Further, we have \ $\bR^{1/2} \tbN \distre \bkappa \tbN$, \ where the random matrix \ $\bkappa$ \ is defined by
 \begin{align}\label{kappa_2}
   \bkappa
   := \begin{bmatrix}
       \sqrt{\frac{w_{X_0}}{\ee^B-1}}
        & 0 \\[2mm]
        \sqrt{\frac{(\ee^B-1) w_{X_0}^3 }
                   {(\ee^{2B}-1)^2}}
         & \sqrt{\frac{\ee^B(\ee^B-1)^2w_{X_0}^3}
                   { (\ee^{2B}-1)^2(\ee^{3B}-1)}}
      \end{bmatrix} .
 \end{align}
Indeed, we have
 \begin{align*}
  \bkappa \bkappa^\top
     = \begin{bmatrix}
       \frac{w_{X_0}}{\ee^B-1}
        & \frac{w_{X_0}^2}{\ee^{2B}-1} \\[1mm]
       \frac{w_{X_0}^2}{\ee^{2B}-1}
        & \left(\frac{\ee^B-1}{(\ee^{2B}-1)^2} +  \frac{\ee^B(\ee^B-1)^2 }{(\ee^{2B}-1)^2(\ee^{3B}-1)} \right)w_{X_0}^3
      \end{bmatrix}
     = \bR.
 \end{align*}

(ii): In Theorem \ref{main_Ad} the random variable \ $\tbN$ \ and the random matrix \ $\bR$ \ are indeed \ $\PP$-independent,
 since \ $\tbN$ \ is \ $\PP$-independent of \ $\cF_\infty^X$ \ and \ $\bR$ \ is \ $\cF_\infty^X$-measurable
  (following from the \ $\cF_\infty^X$-measurability of \ $w_{X_0}$, \ which is explained in Step 1 of the proof of Theorem \ref{main_Ad}).
\proofend
\end{Rem}

In case of \ $C = 0$ \ we have \ $V = 0$, \ yielding that the limit distribution in Theorem \ref{main_Ad} is the two-dimensional zero vector,
 and hence other scaling factors should be found in Theorem \ref{main_Ad} in order to get a non-trivial limit distribution.
This motivates the following theorem.

\begin{Thm}\label{main1_Ad}
Suppose that the assumptions of Theorem \ref{main} and \ $C = 0$ \ hold.
Then
 \begin{itemize}
   \item[(i)] we have
      \[
        n^{-1/2} \sum_{k=1}^n M_k \to \tN_1
        \qquad \text{$\cF_\infty^X$-mixing as \ $n \to \infty$,}
       \]
      where \ $\tN_1$ \ is a random variable \ $\PP$-independent of \ $\cF_\infty^X$ \ such that
      \[
         \tN_1 \distre \cN\left(0, \frac{\ee^{2B}-1}{2B} \int_0^\infty r^2 \, \nu(\dd r)\right),
      \]
   \item[(ii)] we have
       \[
        \ee^{-Bn} \sum_{k=1}^n M_k X_{k-1}
         \to \frac{w_{X_0}}{\ee^B} \sum_{j=0}^\infty \ee^{-Bj} Z_j
         \qquad \text{$\cF_\infty^X$-stably as \ $n \to \infty$,}
       \]
       where \ $(Z_j)_{j\in\ZZ_+}$ \ are \ $\PP$-independent and identically distributed random variables being \ $\PP$-independent from
       \ $\cF_\infty^X$ \ such that \ $Z_1$ \ has a characteristic function
       \[
         \EE(\ee^{\ii\theta Z_1})
          = \exp\biggl\{\int_0^1 \int_0^\infty
                    \Big(\ee^{\ii\theta r \ee^{Bu} }- 1  - \ii \theta r \ee^{Bu} \Big)
                    \, \dd u \, \nu(\dd r)\biggr\},
                    \qquad \theta\in\RR.
        \]
        In particular, we have $Z_1 \distre M_1$.
        Further, the series $\sum_{j=0}^\infty \ee^{-Bj} Z_j$ is absolutely convergent $\PP$-almost surely.
  \end{itemize}
\end{Thm}

The proof of Theorem \ref{main1_Ad} can be found in Section \ref{section_proof_main1_Ad}.

\begin{Rem}
In part (ii) of Theorem \ref{main1_Ad} the sequence of random variables \ $(Z_j)_{j\in\ZZ_+}$ \ and the random variable \ $w_{X_0}$ \ are \ $\PP$-independent,
 since \ $(Z_j)_{j\in\ZZ_+}$ \ is \ $\PP$-independent of \ $\cF_\infty^X$ \ and \ $w_{X_0}$ \ is \ $\cF_\infty^X$-measurable
 (as it is explained in Step 1 of the proof of Theorem \ref{main1_Ad}).
\proofend
\end{Rem}

\section{Comparison of our results on the CLS estimator of \ $(B,A)$ \ and
  those of Huang, Ma and Zhu (2011) on the weighted CLS estimator of \ $(B,a)$}
\label{section_HuangMaZhu}

First, we  will recall the results of Huang et al.\ \cite{HuaMaZhu} on the asymptotic behaviour of
 the weighted CLS estimator of \ $(B,a)$ \ for a supercritical CBI process \ $(X_t)_{t\in\RR_+}$ \
 with parameters \ $(c, a, b, \nu, \mu)$.
We emphasize that our proof technique and that of Huang et al.\ \cite{HuaMaZhu} are completely different.
They use a martingale technique following Wei and Winnicki \cite[Theorem 3.5]{WW},
 while we apply a multidimensional stable limit theorem (see Theorem \ref{MSLTES}).

Recall that the supercritical case means that \ $B>0$.
In this section, further, we suppose that $X_0=x_0\in\RR_+$, \ $A>0$, \ $C>0$ \ and the moment condition \eqref{moment_condition_2} hold.
Then the weighted CLS estimator \ $(\widetilde B_n,\widetilde a_n)$ \ of \ $(B,a)$ \ based on the observations
 $X_1,\ldots,X_n$ takes the following form
 \[
   \widetilde B_n = \log\left(\frac{\sum_{k=1}^n X_k \sum_{k=1}^n \frac{1}{X_{k-1}+1} - n \sum_{k=1}^n \frac{X_k}{X_{k-1}+1}}
                               {\sum_{k=1}^n (X_{k-1}+1)\sum_{k=1}^n \frac{1}{X_{k-1}+1}-n^2}\right),
 \]
 and
 \[
   \widetilde a_n = \frac{\frac{1}{n}\left(\sum_{k=1}^n X_k - \ee^{\widetilde B_n}\sum_{k=1}^n X_{k-1}\right)}
                         {\ee^{\widetilde B_n} - 1}\widetilde B_n - \int_0^\infty r \, \nu(\dd r),
 \]
 provided that the denominators are not zero and the expression after the logarithm is positive,
 see formulae (1.8) and (1.9) with $\Delta t=1$ in Huang et al.\ \cite{HuaMaZhu}.
In the supercritical case, i.e., when $B>0$, Huang et al.\ \cite[Theorem 2.2 and its proof]{HuaMaZhu} proved
 that \ $\widetilde B_n$ \ is a strongly consistent estimator of \ $B$,
 \ $\widetilde a_n$ is not a (weakly) consistent estimator of $a$,
 and the following two convergences in distribution hold
 \begin{align}\label{wCLSE_B_asymp}
   \left(\sum_{k=1}^n (X_{k-1}+1)\right)^{1/2} (\widetilde B_n - B)
      \distrP \ee^{-B}\cN(0,\tau)\qquad \text{as \ $n\to\infty$,}
 \end{align}
 and
 \begin{align*}
   n \left(\sum_{k=1}^n (X_{k-1}+1)\right)^{-1/2}(\widetilde a_n -a)
     \distrP \frac{2B}{\ee^{B}-1}\cN(0,\tau) \qquad \text{as \ $n\to\infty$,}
 \end{align*}
 where
 \[
   \tau:= C\frac{\ee^B(\ee^B -1)}{B}.
 \]
Since \ $\widetilde a_n - a = \widetilde A_n -A$, \ where
 \[
   \widetilde A_n := \frac{\frac{1}{n}\left(\sum_{k=1}^n X_k - \ee^{\widetilde B_n}\sum_{k=1}^n X_{k-1}\right)}
                         {\ee^{\widetilde B_n} - 1}\widetilde B_n ,
 \]
 provided that the denominator is not zero, we also have
 \begin{align}\label{wCLSE_a_asymp}
   n \left(\sum_{k=1}^n (X_{k-1}+1)\right)^{-1/2}(\widetilde A_n -A)
     \distrP \frac{2B}{\ee^{B}-1}\cN(0,\tau) \qquad \text{as \ $n\to\infty$.}
 \end{align}

Next, using \eqref{htBhtb}, we derive convergence in distribution for
 \ $\left(\sum_{k=1}^n (X_{k-1}+1)\right)^{1/2} (\widehat B_n - B)$ \ and
 for \ $n \left(\sum_{k=1}^n (X_{k-1}+1)\right)^{-1/2}(\widehat A_n -A)$ \ as \ $n\to\infty$.
\ Recall that, by \eqref{help_main_Ad0_1}, we have
 \[
   \ee^{-Bn}\sum_{k=1}^n (X_{k-1}+1)\asP \frac{w_{X_0}}{\ee^B -1}
    \qquad \text{as \ $n\to\infty$,}
 \]
 where the random variable $w_{X_0}$ is given in Theorem \ref{convX}, and, by Theorem \ref{w=0}, we have \ $\PP(w_{X_0} > 0) = 1$.
\ Using \eqref{htBhtb} and that stable convergence yields convergence in distribution, we get that
 \[
    \ee^{Bn/2} (\hB_n - B) \distrP  \left(\frac{(\ee^B-1)(\ee^{2B}-1)^2}{B\ee^B(\ee^{3B}-1)} w_{X_0}^{-1}C \right)^{1/2}\zeta_B
    \qquad \text{as \ $n\to\infty$,}
 \]
 and
 \[
  n \ee^{-Bn/2} (\hA_n - A)\distrP \left(\frac{B\ee^{2B}}{(\ee^B-1)(\ee^{3B}-1)} w_{X_0}C\right)^{1/2}\zeta_A
   \qquad \text{as \ $n\to\infty$,}
 \]
 where \ $\zeta_B$ \ and \ $\zeta_A$ \ are standard normally distributed random variables being \ $\PP$-independent of \ $w_{X_0}$.
\ Consequently, using Slutsky's lemma, we have
 \begin{align}\label{wCLSE_B_asymp_2}
  \begin{split}
  \left(\sum_{k=1}^n (X_{k-1}+1)\right)^{1/2} (\widehat B_n - B)
     &= \left(\ee^{-Bn}\sum_{k=1}^n (X_{k-1}+1)\right)^{1/2} \cdot \ee^{Bn/2}(\widehat B_n - B)\\
     & \distrP \left(\frac{w_{X_0}}{\ee^B -1}\right)^{1/2} \left(\frac{(\ee^B-1)(\ee^{2B}-1)^2}{B\ee^B(\ee^{3B}-1)} w_{X_0}^{-1}C \right)^{1/2}\zeta_B\\
     & \phantom{ \distrP \;} = \left(\frac{C(\ee^{2B}-1)^2}{B\ee^B(\ee^{3B}-1)}\right)^{1/2}\,\zeta_B
       \qquad \text{as \ $n\to\infty$,}
  \end{split}
 \end{align}
 and
 \begin{align}\label{wCLSE_a_asymp_2}
  \begin{split}
  n \left(\sum_{k=1}^n (X_{k-1}+1)\right)^{-1/2}(\widehat A_n -A)
   &= \left(\ee^{-Bn}\sum_{k=1}^n (X_{k-1}+1)\right)^{-1/2}\cdot n \ee^{-Bn/2}(\widehat A_n -A) \\
   &\distrP \left(\frac{w_{X_0}}{\ee^B -1}\right)^{-1/2} \left(\frac{B\ee^{2B}}{(\ee^B-1)(\ee^{3B}-1)} w_{X_0}C\right)^{1/2}\zeta_A\\
   &\phantom{ \distrP \;} =  \left(\frac{CB\ee^{2B}}{\ee^{3B}-1}\right)^{1/2}\zeta_A
   \qquad \text{as \ $n\to\infty$.}
  \end{split}
 \end{align}

Note that the random scaling factors $\left(\sum_{k=1}^n (X_{k-1}+1)\right)^{1/2}$ and
 $n \left(\sum_{k=1}^n (X_{k-1}+1)\right)^{-1/2}$ in \eqref{wCLSE_B_asymp_2} and \eqref{wCLSE_a_asymp_2}
 for the CLS estimators \ $\hB_n$ \ and \ $\hA_n$ \ are the same as the scaling factors in \eqref{wCLSE_B_asymp}
 and \eqref{wCLSE_a_asymp} for the weighted CLS estimators \ $\widetilde B_n$ \ and \ $\widetilde A_n$, \ respectively.
Remark also that Huang et al.\ \cite{HuaMaZhu} did not describe the joint asymptotic behavior of $(\widetilde B_n, \widetilde a_n)$,
 they proved convergence of distribution of $\widetilde B_n$ and $\widetilde a_n$, separately.
However, for the estimator \ $(\hB_n,\hA_n)$, \ we prove (joint) stable convergence using deterministic scaling,
 see \eqref{htBhtb}.

Next, we compare the variances of the limit distributions in \eqref{wCLSE_B_asymp_2} and \eqref{wCLSE_a_asymp_2}
 with the corresponding ones \eqref{wCLSE_B_asymp} and \eqref{wCLSE_a_asymp}.
By \eqref{wCLSE_B_asymp_2}, the variance of the limit distribution of $\left(\sum_{k=1}^n (X_{k-1}+1)\right)^{1/2}(\hB_n - B)$ is
 \[
    \frac{C(\ee^{2B}-1)^2}{B\ee^B(\ee^{3B}-1)},
 \]
 and, by \eqref{wCLSE_B_asymp}, the variance of the limit distribution of \ $\left(\sum_{k=1}^n (X_{k-1}+1)\right)^{1/2}(\widetilde B_n - B)$ \ is
 \[
   \ee^{-2B} \tau = \ee^{-2B}\frac{C\ee^B(\ee^B-1)}{B} = \frac{C(\ee^B-1)}{B\ee^B}.
 \]
One can easily check that
 \[
     \frac{C(\ee^{2B}-1)^2}{B\ee^B(\ee^{3B}-1)} >  \frac{C(\ee^B-1)}{B\ee^B}
 \]
 holds if and only if \ $\ee^B(\ee^B-1)^2>0$, \ which is satisfied due to \ $B>0$.
Hence the variance of the limit distribution of  \ $\left(\sum_{k=1}^n (X_{k-1}+1)\right)^{1/2}(\hB_n - B)$ \ is strictly greater then
 that of \ $\left(\sum_{k=1}^n (X_{k-1}+1)\right)^{1/2}(\tB_n - B)$.
\ It is not so surprising, since weighted estimators usually outperform non-weighted ones.
By \eqref{wCLSE_a_asymp_2}, the variance of the limit distribution of \ $n \left(\sum_{k=1}^n (X_{k-1}+1)\right)^{-1/2} (\hA_n - A)$ \ is
 \[
     \frac{CB\ee^{2B}}{\ee^{3B}-1}
 \]
 and, by \eqref{wCLSE_a_asymp}, the variance of the limit distribution of \ $n \left(\sum_{k=1}^n (X_{k-1}+1)\right)^{-1/2} (\widetilde A_n - A)$ \ is
 \[
     \frac{4B^2}{(\ee^B-1)^2} \cdot \tau = \frac{4B^2}{(\ee^B-1)^2}\cdot\frac{C\ee^B(\ee^B-1)}{B} = \frac{4CB\ee^B}{\ee^B-1}.
 \]
One can check that
 \[
  \frac{CB\ee^{2B}}{\ee^{3B}-1} < \frac{4CB\ee^B}{\ee^B-1}
 \]
  holds if and only if
  \[
   4\ee^{3B} - \ee^{2B} + \ee^B -4 = (\ee^B-1)(4\ee^{2B} + 3\ee^{B}+1)>0,
  \]
  which is satisfied due to \ $B>0$.
\ Hence the variance of the limit distribution of \ $n \left(\sum_{k=1}^n (X_{k-1}+1)\right)^{-1/2} (\hA_n - A)$ \ is strictly less
 then that of  \ $n \left(\sum_{k=1}^n (X_{k-1}+1)\right)^{-1/2} (\widetilde A_n - A)$.
This is somewhat surprising, however, remember that \ $\hA_n$ \ and \ $\ta_n$ \ (and hence $\tA_n$) \
 are not (weakly) consistent estimators as $n\to\infty$.
In case of \ $\hA_n$, \ it follows from the assumption \ $C>0$ \ and Remark \ref{REMARK2}, and in case of  \ $\ta_n$, \
 from Huang et al.\ \cite{HuaMaZhu}.

\section{Proof of Theorem \ref{main_Ad0}}
\label{section_proof_Theorem_main_Ad0}

First, we prove \eqref{help_main_Ad0_1}.
Let $\ell\in\NN$.
\ By Theorem \ref{convX}, we have $\ee^{-\ell Bn} X_n^\ell \asP w_{X_0}^\ell$ \ as \ $n \to \infty$.
Recall that, by Kronecker's lemma, if $(\kappa_n)_{n\in\NN}$ is a sequence of positive real numbers such that
 $\kappa_n\uparrow \infty$ as $n\to\infty$ and $(s_n)_{n\in\NN}$ is a sequence of real numbers such that $s_n\to s\in\RR$
 as $n\to\infty$, then
 \[
   \frac{1}{\kappa_n} \sum_{k=1}^n (\kappa_k - \kappa_{k-1})s_k\to s
   \qquad \text{as \ $n\to\infty$,}
 \]
 where \ $\kappa_0:=0$.
\ Let us choose $\kappa_n:=\sum_{k=1}^n \ee^{\ell B k}$, $n\in\NN$, \ with $\kappa_0=0$,
 and \ $s_n:=\ee^{-\ell B n}X_n^\ell$, $n\in\NN$.
\ Since $B>0$, we have
 \[
  \kappa_n=\ee^{\ell B} \frac{\ee^{\ell B n} - 1}{\ee^{\ell B} -1}
            \to \infty \qquad \text{as \ $n\to\infty$,}
 \]
 and it holds that \ $s_n\asP w_{X_0}^\ell$ as $n\to\infty$.
\ Hence, by Kronecker's lemma, we get
 \[
   \frac{\ee^{\ell B} - 1}{\ee^{\ell B}}
     \frac{1}{\ee^{\ell Bn} - 1}
     \sum_{k=1}^n \ee^{\ell B k} (\ee^{-\ell B k} X_k^\ell) \asP w_{X_0}^\ell
       \qquad \text{as \ $n\to\infty$,}
 \]
 yielding that
 \begin{align}\label{help_main_Ad0_3}
     \ee^{-\ell B n} \sum_{k=1}^n X_k^\ell \asP \frac{\ee^{\ell B}}{\ee^{\ell B} -1} w_{X_0}^\ell
         \qquad \text{as \ $n\to\infty$.}
 \end{align}
Hence, using that $\ee^{-\ell B(n-1)} X_0\asP 0$ as $n\to\infty$, we have
 \begin{align*}
     \ee^{-\ell B n} \sum_{k=1}^n X_{k-1}^\ell
        = \ee^{-\ell B n} \sum_{k=0}^{n-1} X_k^\ell
        = \ee^{-\ell B} \ee^{-\ell B (n-1)}
          \sum_{k=0}^{n-1} X_k^\ell
        \asP \ee^{-\ell B} \frac{\ee^{\ell B}}{\ee^{\ell B} -1} w_{X_0}^\ell
            = \frac{w_{X_0}^\ell}{\ee^{\ell B} -1}
 \end{align*}
 as \ $n\to\infty$, \ as desired.

Next, we turn to prove \eqref{help_main_Ad0_2}.
By Theorem \ref{convX}, we have $(\ee^{-Bn} X_n) (\ee^{-B(n-1)}X_{n-1})  \asP w_{X_0}^2$ \ as \ $n \to \infty$, \ and hence
 $\ee^{-2Bn} X_n X_{n-1}  \asP \ee^{-B}w_{X_0}^2$ \ as \ $n \to \infty$.
\ Let us apply the above recalled Kronecker's lemma by choosing $\kappa_n:=\sum_{k=1}^n \ee^{2B k}$, $n\in\NN$, \ with $\kappa_0=0$,
 and \ $s_n:=\ee^{-2Bn} X_n X_{n-1} $, $n\in\NN$.
\ Since $B>0$, we have
 \[
  \kappa_n=\ee^{2B} \frac{\ee^{2B n} - 1}{\ee^{2B} -1}
            \to \infty \qquad \text{as \ $n\to\infty$,}
 \]
 and it holds that $s_n\asP \ee^{-B}w_{X_0}^2$  \ as \ $n \to \infty$.
\ Hence, by Kronecker's lemma, we get
 \[
   \frac{\ee^{2B} - 1}{\ee^{2B}}
     \frac{1}{\ee^{2Bn} - 1}
     \sum_{k=1}^n \ee^{2B k} (\ee^{-2Bk} X_{k-1}X_k) \asP \ee^{-B}w_{X_0}^2
       \qquad \text{as \ $n\to\infty$,}
 \]
 yielding that
 \[
     \ee^{-2Bn} \sum_{k=1}^n X_{k-1}X_k \asP \frac{\ee^{2B}}{\ee^{2B} -1} \ee^{-B}w_{X_0}^2
                                             = \frac{\ee^B}{\ee^{2B} -1} w_{X_0}^2
         \qquad \text{as \ $n\to\infty$,}
 \]
 as desired.

\section{Proof of Lemma \ref{LEMMA_CLSE_exist_discrete}}
\label{section_proof_LEMMA_CLSE_exist_discrete}

By Theorem 3.2 in Barczy et al.\ \cite{BarPalPap2} (see also Theorem \ref{w=0}),
 we have \ $\PP(w_{X_0}=0)=0$, \ yielding that \ $\PP(w_{X_0}>0)=1$ \ (since \ $w_{X_0}$ \ is non-negative).
For each \ $n \in \NN$, \ we have
 \begin{align*}
  \Omega \setminus H_n
  &= \left\{\omega \in \Omega
             : \sum_{k=1}^n X_{k-1}^2(\omega)
               - \frac{1}{n}
                 \left(\sum_{i=1}^n X_{i-1}(\omega)\right)^2 \leq 0\right\} \\
  &= \left\{\omega \in \Omega
             : \sum_{k=1}^n \left( X_{k-1}(\omega)
               - \frac{1}{n} \sum_{i=1}^n X_{i-1}(\omega)\right)^2 \leq  0\right\} \\
  &= \left\{\omega \in \Omega
             : \sum_{k=1}^n \left( X_{k-1}(\omega)
               - \frac{1}{n} \sum_{i=1}^n X_{i-1}(\omega)\right)^2 = 0\right\}
 \end{align*}
 \begin{align*}
  &= \left\{\omega \in \Omega
             : X_{k-1}(\omega) = \frac{1}{n} \sum_{i=1}^n X_{i-1}(\omega) ,
               \,\, k \in \{1, \ldots, n\} \right\} \\
  &= \left\{\omega \in \Omega
             : X_0(\omega) = X_1(\omega) = \cdots = X_{n-1}(\omega) \right\} .
 \end{align*}
Hence
 \[
   \bigcap_{n\in\NN} (\Omega\setminus H_n )
      = \big\{ \omega\in\Omega : X_n(\omega) = X_0(\omega),\; n\in\NN \big\}
      \subset \Big\{  \omega\in\Omega : \lim_{n\to\infty}\ee^{-Bn} \sum_{k=1}^n X_{k-1}(\omega)=0 \Big\},
 \]
 and then, by the continuity and monotonicity of probability, we have
 \begin{align}\label{help_LEMMA_CLSE_exist_discrete_1}
  \begin{split}
   \lim_{n\to\infty}\PP(\Omega\setminus H_n)
      = \PP\left( \bigcap_{n\in\NN} (\Omega\setminus H_n) \right)
      \leq \PP\left( \lim_{n\to\infty} \ee^{-Bn}  \sum_{k=1}^n X_{k-1} = 0\right).
 \end{split}
 \end{align}
By Theorem \ref{main_Ad0}, we get
 \[
  \PP\left( \lim_{n\to\infty} \ee^{-Bn}  \sum_{k=1}^n X_{k-1} = \frac{w_{X_0}}{\ee^B-1} \right) = 1,
 \]
 where \ $\PP(w_{X_0}>0)=1$.
\ Consequently, we have \ $\PP\left( \lim_{n\to\infty} \ee^{-Bn}  \sum_{k=1}^n X_{k-1} = 0\right)=0$, \ and then
 \eqref{help_LEMMA_CLSE_exist_discrete_1} yields that \ $\lim_{n\to\infty}\PP(\Omega\setminus H_n) =0 $, \ i.e,
 \ $\lim_{n\to\infty}\PP(H_n) =1$, \ as desired.
Taking into account \eqref{CLSErb}, this also yields the last statement of the assertion.

\section{Proof of Theorem \ref{Thm_strong_consistency}}
\label{section_proof_Thm_strong_consistency}

First, we prove an auxiliary lemma.

\begin{Lem}\label{Lem_strong_consistency}
Let \ $\xi_n$, \ $\eta_n$, \ $n\in\NN$, \ and \ $\eta$ \ be random variables on a probability space \ $(\Omega,\cF,\PP)$ \
 such that \ $\eta_n\asP \eta$ \ as \ $n\to\infty$.
\ Let \ $H_n$, \ $n\in\NN$, \ be a sequence of events in \ $\cF$ \ such that \ $H_n\subset H_m$ \ for each \ $n\leq m$, \ $n,m\in\NN$, \
 $H_n\subset \{\xi_n=\eta_n\}$ \ for each \ $n\in\NN$, \ and $\lim_{n\to\infty}\PP(H_n)=1$.
\ Then \ $\xi_n\asP \eta$ \ as \ $n\to\infty$.
\end{Lem}

\noindent{\bf Proof.}
For each $N\in\NN$, we have
 \begin{align*}
   \PP\left(\lim_{n\to\infty}\xi_n = \eta \right)
    & = \PP\left( \Big\{\lim_{n\to\infty}\xi_n = \eta \Big\}\cap H_N \right)
        + \PP\left( \Big\{\lim_{n\to\infty}\xi_n = \eta \Big\}\cap \overline{H_N} \right) \\
    & = \PP\left( \Big\{\lim_{n\to\infty}\xi_n = \eta \Big\}\cap \bigcap_{\ell=N}^\infty H_\ell \right)
         + \PP\left( \Big\{\lim_{n\to\infty}\xi_n = \eta \Big\}\cap \overline{H_N} \right)  \\
    & = \PP\left( \Big\{\lim_{n\to\infty}\eta_n = \eta \Big\}\cap \bigcap_{\ell=N}^\infty H_\ell \right)
         + \PP\left( \Big\{\lim_{n\to\infty}\xi_n = \eta \Big\}\cap \overline{H_N}   \right)  \\
    & = \PP\left( \Big\{\lim_{n\to\infty}\eta_n = \eta \Big\}\cap H_N \right)
         + \PP\left( \Big\{\lim_{n\to\infty}\xi_n = \eta \Big\}\cap \overline{H_N} \right),
 \end{align*}
 where \ $\overline{H_N}:=\Omega\setminus H_N$.
\ Here
 \[
   \PP\left(\Big\{\lim_{n\to\infty}\xi_n = \eta \Big\}\cap \overline{H_N} \right)
     \leq \PP(\overline{H_N})\to 0   \qquad \text{as \ $N\to\infty$.}
 \]
Further, we have \ $H_N \uparrow \bigcup_{\ell=1}^\infty H_\ell$ \ as \ $N\to\infty$,
 and hence, by the assumption \ $\lim_{n\to\infty}\PP(H_n)=1$ \ and the continuity of probability,
 we get
 $1=\lim_{N\uparrow \infty} \PP(H_N) = \PP\left(  \bigcup_{\ell=1}^\infty H_\ell \right)$.
Consequently, we have
 \begin{align*}
   \PP\left(\lim_{n\to\infty}\xi_n = \eta \right)
     & = \lim_{N\to\infty } \PP\left( \Big\{\lim_{n\to\infty}\eta_n = \eta \Big\}\cap H_N \right)\\
     & = \PP\left( \Big\{\lim_{n\to\infty}\eta_n = \eta \Big\}\cap \bigcup_{\ell=1}^\infty H_\ell \right)
       = \PP\left(\lim_{n\to\infty}\eta_n = \eta \right)=1,
 \end{align*}
 where at the second equality we used again the continuity of probability and at the last equality the assumption
 \ $\PP(\lim_{n\to\infty}\eta_n = \eta)=1$.
\proofend

\bigskip

Now we turn to prove Theorem \ref{Thm_strong_consistency}.
Recall that \ $\hB_n = \log(\hvarrho_n)$ \ on the set \ $\{\hvarrho_n\in\RR_{++}\}$ \ (see \eqref{help_B_A_est_definition})
 and $B=\log(\varrho)$, \ where \ $B>0$ \ yielding that $\varrho>1$.
Since the function $\log$ is continuous, it is enough to prove that the CLS estimator \ $\hvarrho_n$ \ of \ $\varrho$ \ is strongly consistent,
 i.e., \ $\hvarrho_n\asP \varrho$ \ as $n\to\infty$.
Recall that, by \eqref{CLSErb}, for each \ $n \in \NN$, \ we have
 \begin{align}\label{help_strong_consistency2}
 \begin{split}
    \hvarrho_n
     & = \frac{1}{n \sum_{k=1}^n X_{k-1}^2
                 - \left( \sum_{k=1}^n X_{k-1} \right)^2}
         \left( n \sum_{k=1}^n X_k X_{k-1}
               - \sum_{k=1}^n X_k \sum_{k=1}^n X_{k-1}
         \right) \\
     &= \frac{1}{\ee^{-2Bn}\sum_{k=1}^n X_{k-1}^2
                 - n^{-1}\left( \ee^{-Bn} \sum_{k=1}^n X_{k-1} \right)^2} \\
     &\phantom{=\;}\times\left( \ee^{-2Bn} \sum_{k=1}^n X_k X_{k-1}
                - n^{-1}\left(\ee^{-Bn} \sum_{k=1}^n X_k \right) \left(\ee^{-Bn} \sum_{k=1}^n X_{k-1} \right)
         \right)
 \end{split}
 \end{align}
 on the set \ $H_n$.
\ Further, by Lemma \ref{LEMMA_CLSE_exist_discrete} and its proof, we get \ $\PP(w_{X_0}>0)=1$, \ $\PP(H_n)\to 1$ \ as \ $n\to\infty$,
 \ and \ $\Omega\setminus H_n = \{\omega\in \Omega : X_0(\omega) = X_1(\omega) = \cdots = X_{n-1}(\omega)\}$, $n\in\NN$,
 \ and hence \ $\Omega\setminus H_m \subset \Omega\setminus H_n$ \ for each \ $n\leq m$, $n,m\in\NN$,
 \ yielding that \ $H_n\subset H_m$ \ for each \ $n\leq m$, $n,m\in\NN$.
\ Moreover, \eqref{help_main_Ad0_3}, \eqref{help_main_Ad0_1} and \eqref{help_main_Ad0_2} yield that
 \begin{align*}
   &\ee^{-Bn} \sum_{k=1}^n X_k   \asP \frac{\ee^B}{\ee^B-1}w_{X_0}\qquad \text{as \ $n\to\infty$,}\\
   &\ee^{-Bn} \sum_{k=1}^n X_{k-1} \asP \frac{1}{\ee^B-1}w_{X_0}\qquad \text{as \ $n\to\infty$,}\\
   &\ee^{-2Bn} \sum_{k=1}^n X_{k-1}^2 \asP \frac{1}{\ee^{2B}-1}w_{X_0}^2\qquad \text{as \ $n\to\infty$,}\\
   &\ee^{-2Bn} \sum_{k=1}^n X_kX_{k-1} \asP \frac{\ee^B}{\ee^{2B}-1}w_{X_0}^2\qquad \text{as \ $n\to\infty$.}
 \end{align*}
Using Lemma \ref{Lem_strong_consistency}, \eqref{help_strong_consistency2},  $\PP(H_n)\to 1$ \ as \ $n\to\infty$, \
 $H_n\subset H_m$ \ for each \ $n\leq m$, $n,m\in\NN$, \ and that \ $\PP(w_{X_0}>0)=1$, \ we get
 \begin{align*}
   \hvarrho_n
     \asP \frac{1}{\frac{1}{\ee^{2B}-1}w_{X_0}^2 - 0}\left( \frac{\ee^B}{\ee^{2B}-1}w_{X_0}^2 - 0\right)
         = \ee^B = \varrho \qquad \text{as \ $n\to\infty$,}
 \end{align*}
 as desired.

\section{Proof of Theorem \ref{main}}
\label{section_proof_main}

First, note that the random matrix \ $\bS$ \ is symmetric and positive definite \ $\PP$-almost surely, since
 \ $\frac{(\ee^B-1)(\ee^{2B}-1)^2}{B\ee^B(\ee^{3B}-1)}w_{X_0}^{-1}>0$ \ $\PP$-almost surely (due to Theorem \ref{w=0}) and
 \begin{align*}
  \det(\bS)
    = \frac{\ee^B(\ee^{2B}-1)^2}{(\ee^{3B}-1)^2}
      - \frac{\ee^{2B}(\ee^B-1)^2}{(\ee^{3B}-1)^2}
    = \frac{\ee^B(\ee^B-1)}{\ee^{3B}-1} > 0.
 \end{align*}
Consequently, there exists a unique symmetric and positive definite square root \ $\bS^{1/2}$ \ of \ $\bS$ \
 \ $\PP$-almost surely.

Using Theorem \ref{Thm_strong_consistency} we have already checked that
 for all \ $\delta\in(0,1)$ \ there exists an \ $n_0\in\NN$ \ such that
 \ $\PP(\hvarrho_n \in\RR_{++}) > 1-\delta$ for each \ $n\geq n_0$,
 see \eqref{help_est_existence_1}.
Consequently, due to the definition of \ $(\hB_n, \hA_n)$ \ given in \eqref{help_B_A_est_definition},
 we have that for all \ $\delta\in(0,1)$ \ there exists an \ $n_0\in\NN$ \ such that
 the probability of the existence of the estimator \ $(\hB_n, \hA_n)$ \ is greater than
 $1-\delta$ for each \ $n\geq n_0$.
Hence, the probability of the existence of the estimator \ $(\hB_n, \hA_n)$ \ converges to 1 as $n\to\infty$.

Next, we apply Lemma \ref{Lem_Kallenberg_stable} with
 \[
   \xi_n := \begin{bmatrix}
           \ee^{Bn/2} (\hvarrho_n - \varrho) \\[1mm]
           n \ee^{-Bn/2} (\hcA_n - \cA)
           \end{bmatrix} , \qquad
   \xi := \tbS^{1/2} \tbN ,
 \]
  and with functions \ $f : \RR^2 \to \RR^2$ \ and \ $f_n : \RR^2 \to \RR^2$,
 \ $n \in \NN$, \ given by
 \[
   f\Biggl( \begin{bmatrix} x \\ y \end{bmatrix} \Biggr)
   := \begin{bmatrix}
       \frac{x}{\varrho} \\
       \frac{y}{\int_0^1 \varrho^s \, \dd s} \end{bmatrix} , \quad (x, y) \in \RR^2 , \qquad
   f_n\Biggl( \begin{bmatrix} x \\ y \end{bmatrix} \Biggr)
   := \begin{bmatrix}
       \ee^{Bn/2} \log\bigl( 1 + \frac{x}{\ee^{Bn/2}\varrho} \bigr) \\[1mm]
       {\DS\frac{y+n\ee^{-Bn/2}\cA}
            {\int_0^1 (\varrho + \frac{x}{\ee^{Bn/2}})^s \, \dd s}}
       - n \ee^{-Bn/2} A
      \end{bmatrix}
 \]
 for \ $(x, y) \in \RR^2$ \ with \ $x > -\ee^{Bn/2}\varrho$, \ and \ $f_n(x, y) := 0$
 \ otherwise.
We have
 \[
   f_n\Biggl( \begin{bmatrix}
               \ee^{Bn/2} (\hvarrho_n - \varrho) \\
               n \ee^{-Bn/2} (\hcA_n - \cA)
              \end{bmatrix} \Biggr)
   = \begin{bmatrix}
      \ee^{Bn/2} (\hB_n - B) \\
      n \ee^{-Bn/2} (\hA_n - A)
     \end{bmatrix}
 \]
 on the set \ $\{\omega \in \Omega : \hvarrho_n(\omega) \in \RR_{++}\}$.
Indeed, using \eqref{help_B_A_est_definition}, one can easily see that
 \begin{align*}
    \ee^{Bn/2} \log\left( 1 + \frac{\ee^{Bn/2}(\hvarrho_n - \varrho)}{\ee^{Bn/2} \varrho}  \right)
    = \ee^{Bn/2} (\log(\hvarrho_n) - \log(\varrho) )
    = \ee^{Bn/2} (\hB_n - B),
 \end{align*}
 and
 \begin{align*}
   \frac{ n\ee^{-Bn/2}(\hcA_n - \cA) +n\ee^{-Bn/2}\cA}
             {\int_0^1 (\varrho + \frac{\ee^{Bn/2}(\hvarrho_n - \varrho)}{\ee^{Bn/2}})^s \, \dd s}
         - n \ee^{-Bn/2} A
   =  \frac{n\ee^{-Bn/2}\hcA_n}{\int_0^1 (\hvarrho_n)^s\,\dd s}  - n \ee^{-Bn/2} A
   = n\ee^{-Bn/2}(\hA_n - A)
 \end{align*}
 on the set \ $\{\omega \in \Omega : \hvarrho_n(\omega) \in \RR_{++}\}$.
Further, we get \ $f_n(x_n, y_n) \to f(x, y)$ \ as \ $n \to \infty$ \ if
 \ $(x_n, y_n) \to (x, y)\in\RR^2$ \ as \ $n \to \infty$.
Indeed, if \ $x_n \to x\in\RR$ \ as \ $n \to \infty$, then \ $x_n > -\ee^{Bn/2}\varrho$ \ for sufficiently large \ $n\in\NN$,
 \ $x_n\ee^{-Bn/2}\varrho^{-1} \to 0$ \ as \ $n\to\infty$,
 and, using that $\frac{z}{z+1}\leq \log(1+z)\leq z$ for $z>-1$, we have
 \[
   \frac{\frac{x_n}{\varrho}}{\frac{x_n}{\ee^{Bn/2}\varrho} +1}
       \leq \ee^{Bn/2} \log\left( 1 + \frac{x_n}{\ee^{Bn/2}\varrho} \right)
       \leq \frac{x_n}{\varrho}
 \]
 for sufficiently large $n\in\NN$, which, by the sandwich theorem, implies that
 \[
   \lim_{n\to\infty} \ee^{Bn/2} \log\Bigl( 1 + \frac{x_n}{\ee^{Bn/2}\varrho} \Bigr) = \frac{x}{\varrho} .
 \]
We also have if \ $x_n \to x\in\RR$ \ as \ $n \to \infty$, then
 \ $\lim_{n\to\infty} \int_0^1 (\varrho + \frac{x_n}{\ee^{Bn/2}})^s \, \dd s = \int_0^1 \varrho^s \, \dd s$, \
 \ since the function \ $\RR_{++} \ni u \mapsto \int_0^1 u^s \, \dd s \in \RR$ \ is continuous
 or one can argue directly as follows
 \begin{align*}
  \left\vert \int_0^1 \Big(\varrho + \frac{x_n}{\ee^{Bn/2}}\Big)^s \, \dd s
             - \int_0^1 \varrho^s \, \dd s  \right\vert
   &= \left\vert \frac{x_n}{\ee^{Bn/2}} \int_0^1 s\Big(\varrho + \frac{\theta_{n,s}}{\ee^{Bn/2}}\Big)^{s-1} \, \dd s \right\vert\\
   &\leq \sup_{s\in[0,1]} \Big(\varrho + \frac{\theta_{n,s}}{\ee^{Bn/2}}\Big)^{s-1} \frac{\vert x_n\vert}{\ee^{Bn/2}}
   \to 0 \qquad \text{as \ $n\to\infty$,}
 \end{align*}
 where $\theta_{n,s}$ (depending on $x_n$ and $s$) lies between $x_n$ and $0$, since
 for sufficiently large $n\in\NN$ we have \ $\vert \theta_{n,s}/\ee^{Bn/2}\vert <1$ \ (following from \ $x_n\ee^{-Bn/2}\to 0$ \ as \ $n\to\infty$)
 \ and then
 \[
   0\leq \limsup_{n\to\infty}\sup_{s\in[0,1]}\Big(\varrho + \frac{\theta_{n,s}}{\ee^{Bn/2}}\Big)^{s-1}
     \leq \sup_{s\in[0,1]} (\varrho - 1)^{s-1} < \infty,
 \]
 since the function \ $\RR_{++}\ni u \mapsto u^{s-1}$ \ is monotone decreasing for all $s\in[0,1]$.
\ Hence \eqref{hvarrhohobeta}, Lemma \ref{Lem_Kallenberg_stable}, \ $\lim_{n\to\infty}\PP(\hvarrho_n \in \RR_{++})=1$ \
 and what is explained before Theorem \ref{main} imply that
 \begin{align*}
    \begin{bmatrix}
      \ee^{Bn/2} (\hB_n - B) \\
      n \ee^{-Bn/2} (\hA_n - A)
     \end{bmatrix}
     \to
     f(\tbS^{1/2} \tbN)
   \qquad \text{$\cF_\infty^X$-stably as \ $n \to \infty$,}
 \end{align*}
 where \ $\tbS$ \ and \ $\tbN$ \ are given in Theorem \ref{main_rb}.
Note that \ $\int_0^1 \varrho^s \, \dd s = \frac{\varrho-1}{\log(\varrho)}$, \ since \ $\varrho = \ee^B\in(1,\infty)$,
 \ and hence
 \[
 f\Biggl( \begin{bmatrix} x \\ y \end{bmatrix} \Biggr)
   = \begin{bmatrix}
       \frac{1}{\varrho} & 0 \\
       0 & \frac{\log(\varrho)}{\varrho-1} \\
     \end{bmatrix}
      \begin{bmatrix}
       x \\
       y
       \end{bmatrix} , \quad (x, y) \in \RR^2 .
 \]
Consequently, in order to prove \eqref{htBhtb}, it remains to check that
 \begin{align}\label{help3}
    \begin{bmatrix}
       \frac{1}{\varrho} & 0 \\
       0 & \frac{\log(\varrho)}{\varrho-1} \\
     \end{bmatrix}
      \tbS^{1/2} \tbN
       \distre
    \bS^{1/2} \bN.
 \end{align}
Since $\bS$ and $\bN$ are \ $\PP$-independent, $\tbS$ and $\tbN$ are \ $\PP$-independent,
 the conditional distribution of $\bS^{1/2}\bN$ given $\bS$ (equivalently, given \ $w_{X_0}$)
 and that of
 \[
 \begin{bmatrix}
       \frac{1}{\varrho} & 0 \\
       0 & \frac{\log(\varrho)}{\varrho-1} \\
     \end{bmatrix}
      \tbS^{1/2} \tbN
 \]
 given \ $\tbS$ \ (equivalently, given \ $w_{X_0}$) \
 are two-dimensional normal distributions with mean vector $\bzero\in\RR^2$,
 in order to show \eqref{help3}, it is enough to verify that the covariance matrices of the conditional distributions in question coincide.
Since \ $V = C\int_0^1 \ee^{B(1+u)}\,\dd u = C\frac{\ee^B(\ee^{B} - 1)}{B}$ \ and \ $\varrho = \ee^B$,
 \ the covariance matrix of
  \[
 \text{the conditional distribution of} \quad
        \begin{bmatrix}
       \frac{1}{\varrho} & 0 \\
       0 & \frac{\log(\varrho)}{\varrho-1} \\
       \end{bmatrix}
       \tbS^{1/2} \tbN \quad
   \text{given \ $\tbS$}
 \]
 takes the form
 \begin{align*}
    &C\frac{\ee^B(\ee^{B} - 1)}{B}
      \begin{bmatrix}
        \frac{1}{\varrho} & 0 \\
        0 & \frac{\log(\varrho)}{\varrho-1} \\
      \end{bmatrix}
      \tbS
      \begin{bmatrix}
        \frac{1}{\varrho} & 0 \\
        0 & \frac{\log(\varrho)}{\varrho-1} \\
      \end{bmatrix}\\
    &\qquad = C\frac{\ee^B(\ee^{B} - 1)}{B}
       \begin{bmatrix}
         \frac{(\ee^{2B} - 1)^2}{\ee^B (\ee^{3B} - 1)} w_{X_0}^{-1} & -\frac{\ee^{B} - 1}{\ee^{3B} - 1} \\
         -\frac{B\ee^{B}}{\ee^{3B} - 1} & \frac{B\ee^B}{(\ee^{B} - 1)(\ee^{3B} - 1)}w_{X_0} \\
       \end{bmatrix}
      \begin{bmatrix}
        \frac{1}{\varrho} & 0 \\
        0 & \frac{\log(\varrho)}{\varrho-1} \\
      \end{bmatrix}\\
    &\qquad = C\bS,
 \end{align*}
 which is nothing else but the covariance matrix of the conditional distribution of $\bS^{1/2} \bN$ given $\bS$, as desired.

In case of \ $C = 0$, \ to prove \eqref{htBhtb1}, we apply Lemma \ref{Lem_Kallenberg_stable} with
 \[
   \xi_n := \ee^{Bn} (\hvarrho_n - \varrho) , \qquad
   \xi := \frac{\ee^{2B}-1}{\ee^B}w_{X_0}^{-1}
         \sum_{j=0}^\infty \ee^{-Bj} Z_j ,
 \]
 and with functions \ $f : \RR \to \RR$ \ and \ $f_n : \RR \to \RR$,
 \ $n \in \NN$, \ given by
 \[
   f(x) := \frac{x}{\varrho} , \quad x \in \RR , \qquad
   f_n(x)
   := \ee^{Bn} \log\biggl( 1 + \frac{x}{\ee^{Bn}\varrho} \biggr)
 \]
 for \ $x \in \RR$ \ with \ $x > -\ee^{Bn} \varrho$, \ and \ $f_n(x) := 0$
 \ otherwise.
We have again \ $f_n(\ee^{Bn} (\hvarrho_n - \varrho) )= \ee^{Bn}(\hB_n - B)$ \ on the set \ $\{\omega\in\Omega: \hvarrho_n(\omega)\in\RR_{++}\}$,
 \ and \ $f_n(x_n) \to f(x)$ \ as \ $n \to \infty$ \ if \ $x_n \to x\in\RR$ \ as \ $n \to \infty$.
\ Consequently, \eqref{hvarrhohobeta1}, Lemma \ref{Lem_Kallenberg_stable}, \ $\lim_{n\to\infty}\PP(\hvarrho_n \in \RR_{++})=1$ \
 and what is explained before Theorem \ref{main} imply \eqref{htBhtb1}.

Finally, we turn to prove \eqref{htBhtb2}.
Using the definitions of $\hA_n$, $n\in\NN$ \ (given in \eqref{help_B_A_est_definition}) and of \ $\cA$ \ (given in \eqref{ttBbeta}),
 for each \ $n \in \NN$, \ we have
 \begin{align}\label{help12}
   n^{1/2} (\hA_n - A)
   = \frac{n^{1/2}(\hcA_n-\cA)}
          {\int_0^1(\hvarrho_n)^s\,\dd s}
     + \frac{A n^{1/2}
             \bigl(\int_0^1\varrho^s\,\dd s
                   -\int_0^1(\hvarrho_n)^s\,\dd s\bigr)}
            {\int_0^1(\hvarrho_n)^s\,\dd s}
 \end{align}
 on the set \ $\{\omega \in \Omega : \hvarrho_n(\omega) \in \RR_{++}\}$.
\ Using Theorem \ref{Thm_strong_consistency}, we get \ $\hvarrho_n\asP\varrho$ \ as \ $n \to \infty$, \
\ and hence \ $\int_0^1(\hvarrho_n)^s\,\dd s\asP \int_0^1\varrho^s\,\dd s = \frac{\varrho -1}{\log(\varrho)}$ \
 as \ $n \to \infty$, \ since the function \ $\RR_{++}\ni u \mapsto \int_0^1 u^s \, \dd s \in \RR$ \ is continuous.
Next, we apply Lemma \ref{Lem_Kallenberg} with
 \[
   \xi_n := \ee^{Bn}(\hvarrho_n - \varrho)  , \qquad
   \xi := \frac{\ee^{2B}-1}{\ee^B}w_{X_0}^{-1} \sum_{j=0}^\infty \ee^{-Bj} Z_j,
 \]
 and with functions \ $f : \RR \to \RR$ \ and \ $f_n : \RR \to \RR$,
 \ $n \in \NN$, \ given by
 \[
   f(x) := 0 , \quad x \in \RR , \qquad
   f_n(x)
   := n^{1/2}
      \biggl(\int_0^1\Bigl(\varrho
             + \frac{x}{\ee^{Bn}}\Bigr)^s\,\dd s
             - \int_0^1 \varrho^s\,\dd s \biggr)
 \]
 for \ $x \in \RR$ \ with \ $x > -\ee^{Bn} \varrho$, \ and \ $f_n(x) := 0$
 \ otherwise.
We have
 \[
   f_n( \ee^{Bn}(\hvarrho_n - \varrho) )
      =  n^{1/2} \left( \int_0^1 (\hvarrho_n)^s\,\dd s - \int_0^1 \varrho^s\,\dd s \right)
 \]
 on the set \ $\{\omega \in \Omega : \hvarrho_n(\omega) \in \RR_{++}\}$, \ and
 \ $f_n(x_n) \to f(x)$ \ as \ $n \to \infty$ \ if \ $x_n \to x\in\RR$ \ as \ $n \to \infty$.
\ Indeed, if \ $x_n\to x\in\RR$ \ as \ $n\to\infty$, \ then \ $x_n> -\ee^{Bn}\varrho$ \ for sufficiently large
 \ $n\in\NN$, and
 \begin{align*}
  |f_n(x_n) - f(x)|
  &= n^{1/2}
     \biggl|\int_0^1
             \Bigl(\varrho
                   + \frac{x_n}
                          {\ee^{Bn}}
             \Bigr)^s \dd s
            -\int_0^1 \varrho^s
                \dd s\biggr| \\
  &= n^{1/2}
     \biggl|\frac{x_n}{\ee^{Bn}}
            \int_0^1s\Bigl(\varrho + \frac{\theta_{n,s}}{\ee^{Bn}}\Bigr)^{s-1}\dd s\biggr|
  \leq \sup_{s\in[0,1]} \Bigl(\varrho + \frac{\theta_{n,s}}{\ee^{Bn}}\Bigr)^{s-1}
        n^{1/2} \frac{|x_n|}{\ee^{Bn}}
  \to 0
 \end{align*}
 as \ $n \to \infty$, \ where \ $\theta_{n,s}$ \ (depending on \ $x_n$ \ and \ $s$) \ lies between \ $x_n$ \ and \ $0$, \ since
 for sufficiently large \ $n\in\NN$ \ we have \ $\vert \theta_{n,s}/\ee^{Bn}\vert<1$ \ (following from \ $x_n\ee^{-Bn}\to 0$ \ as \ $n\to\infty$) \
 and then
 \[
   0\leq \limsup_{n\to\infty} \sup_{s\in[0,1]} \Bigl(\varrho + \frac{\theta_{n,s}}{\ee^{Bn}}\Bigr)^{s-1}
    \leq \sup_{s\in[0,1]}(\varrho - 1)^{s-1}<\infty,
 \]
 since the function \ $\RR_{++}\ni u \mapsto u^{s-1}$ \ is monotone decreasing for all $s\in[0,1]$.
Consequently, using \eqref{hvarrhohobeta1}, the fact that stable convergence yields convergence in distribution,
  \ $\lim_{n\to\infty}\PP(\hvarrho_n \in \RR_{++})=1$ \ and what is explained before Theorem \ref{main}, we get
 \[
  n^{1/2}\biggl( \int_0^1(\hvarrho_n)^s\,\dd s - \int_0^1\varrho^s\,\dd s \biggr)
    \distrP
    f\left(\frac{\ee^{2B}-1}{\ee^B} w_{X_0}^{-1}
      \sum_{j=0}^\infty \ee^{-Bj} Z_j\right)
      =0
    \qquad \text{as \ $n\to\infty$.}
 \]
Hence we have
 \[
   n^{1/2}\biggl( \int_0^1(\hvarrho_n)^s\,\dd s - \int_0^1\varrho^s\,\dd s \biggr)
   \stoch 0 \qquad \text{as \ $n \to \infty$.}
 \]
This together with \ $\int_0^1 (\hvarrho_n)^s\,\dd s\asP \int_0^1 \varrho^s\,\dd s = \frac{\varrho-1}{\log(\varrho)}=\frac{\ee^B-1}{B}>0$
 as $n\to\infty$ \ (which was already checked) yield that
 \begin{align}\label{help1_new}
   \frac{n^{1/2}\biggl( \int_0^1(\hvarrho_n)^s\,\dd s - \int_0^1\varrho^s\,\dd s \biggr)}{\int_0^1 (\hvarrho_n)^s\,\dd s}
    \stoch 0 \qquad \text{as \ $n \to \infty$.}
 \end{align}
Using Theorem 3.18 in H\"ausler and Luschgy \cite{HauLus} (see also Theorem \ref{Thm_HL_Thm3_18}),
  \eqref{hvarrhohobeta2} and \ $\big(\int_0^1(\hvarrho_n)^s\,\dd s\big)^{-1}\asP \frac{B}{\ee^B-1}$
 \ as \ $n\to\infty$, we conclude that
 \[
  \frac{n^{1/2}(\hcA_n-\cA)}{\int_0^1(\hvarrho_n)^s\,\dd s}\to \frac{B}{\ee^B-1}\tN_1   \qquad \text{$\cF_\infty^X$-stably as \ $n \to \infty$,}
 \]
  where \ $\tN_1$ \ is a random variable \ $\PP$-independent of \ $\cF_\infty^X$ \ such that
 \ $\tN_1 \distre \cN\bigl(0, \frac{\ee^{2B}-1}{2B} \int_0^\infty r^2 \, \nu(\dd r)\bigr)$.
Hence, using \eqref{help12}, \eqref{help1_new},  Theorem 3.18 in H\"ausler and Luschgy \cite{HauLus} (see also Theorem \ref{Thm_HL_Thm3_18}),
 \ $\lim_{n\to\infty}\PP(\hvarrho_n\in\RR_{++})=1$ \ and what is explained before Theorem \ref{main}, we obtain that
 \[
     n^{1/2} (\hA_n - A) \to \frac{B}{\ee^B -1}\tN_1 \qquad \text{$\cF_\infty^X$-stably as \ $n \to \infty$.}
 \]
Since \ $\frac{B}{\ee^B-1}\tN_1$ \ is \ $\PP$-independent of \ $\cF_\infty^X$ \ as well, we get the \ $\cF_\infty^X$-mixing
 convergence \eqref{htBhtb2}.

\section{Proof of Theorem \ref{main_rb}}
\label{section_proof_main_rb}

First, note that the random matrix \ $\tbS$ \ is symmetric and positive definite \ $\PP$-almost surely, since
 \ $\frac{(\ee^{2B}-1)^2}{\ee^{3B}-1}w_{X_0}^{-1}>0$ \ $\PP$-almost surely (due to Theorem \ref{w=0}) and
 \begin{align*}
  \det(\tbS)
    = \frac{\ee^B(\ee^{2B}-1)^2}{(\ee^{3B}-1)^2}
      - \frac{\ee^{2B}(\ee^B-1)^2}{(\ee^{3B}-1)^2}
    = \frac{\ee^B(\ee^B-1)}{\ee^{3B}-1} > 0.
 \end{align*}
Consequently, there exists a unique symmetric and positive definite square root \ $\tbS^{1/2}$ \ of \ $\tbS$ \
 \ $\PP$-almost surely.

The statement about the existence of a unique CLS estimator
 \ $(\hvarrho_n, \hcA_n)$ \ of \ $(\varrho,\cA)$ \ under the given conditions follows from Lemma
 \ref{LEMMA_CLSE_exist_discrete}.

Next, we derive \eqref{hvarrhohobeta} using Theorems \ref{main_Ad0} and \ref{main_Ad}.
Using \eqref{Mk} and \eqref{CLSErb}, for each \ $n \in \NN$, \ we have
 \begin{align*}
  & \begin{bmatrix}
    \hvarrho_n - \varrho \\[1mm]
    \hcA_n - \cA
   \end{bmatrix} = \\
  & = \frac{1}{n \sum\limits_{k=1}^n X_{k-1}^2
              - \left( \sum\limits_{k=1}^n X_{k-1} \right)^2}
      \begin{bmatrix}
         n \sum\limits_{k=1}^n X_{k-1}(X_k - \varrho X_{k-1})
             -  \sum\limits_{k=1}^n X_{k-1} \sum\limits_{k=1}^n (X_k - \varrho X_{k-1})  \\[3mm]
         \sum\limits_{k=1}^n X_{k-1}^2 \sum\limits_{k=1}^n (X_k - \cA) - \sum\limits_{k=1}^n X_{k-1}\sum\limits_{k=1}^n X_{k-1}(X_k-\cA) \\
     \end{bmatrix}\\
  & = \frac{1}{n \sum\limits_{k=1}^n X_{k-1}^2
              - \left( \sum\limits_{k=1}^n X_{k-1} \right)^2}
      \begin{bmatrix}
         n \sum\limits_{k=1}^n X_{k-1}(M_k+\cA)
             -  \sum\limits_{k=1}^n X_{k-1} \sum\limits_{k=1}^n (M_k + \cA)  \\[3mm]
         \sum\limits_{k=1}^n X_{k-1}^2 \sum\limits_{k=1}^n (M_k + \varrho X_{k-1}) - \sum\limits_{k=1}^n X_{k-1}\sum\limits_{k=1}^n X_{k-1}(M_k+\varrho X_{k-1}) \\
     \end{bmatrix}
 \end{align*}
 \begin{align*}
  & = \frac{1}{n \sum\limits_{k=1}^n X_{k-1}^2
              - \left( \sum\limits_{k=1}^n X_{k-1} \right)^2}
     \begin{bmatrix}
      n \sum\limits_{k=1}^n M_k X_{k-1}
      - \sum\limits_{k=1}^n M_k \sum\limits_{k=1}^n X_{k-1} \\[3mm]
      \sum\limits_{k=1}^n M_k \sum\limits_{k=1}^n X_{k-1}^2
      - \sum\limits_{k=1}^n M_k X_{k-1} \sum\limits_{k=1}^n X_{k-1}
     \end{bmatrix}
 \end{align*}
 on the set \ $H_n$, \ where \ $\PP(H_n)\to 1$ \ as \ $n\to\infty$ \ (see Lemma \ref{LEMMA_CLSE_exist_discrete}).
Consequently, for each \ $n\in\NN$, \ we have
 \begin{align*}
  \begin{bmatrix}
   \ee^{Bn/2} (\hvarrho_n - \varrho) \\[1mm]
   n \ee^{-Bn/2} (\hcA_n - \cA)
  \end{bmatrix}
  &= \frac{1}{\ee^{-2Bn} \sum\limits_{k=1}^n X_{k-1}^2
              - n^{-1} \Bigl(\ee^{-Bn} \sum\limits_{k=1}^n X_{k-1} \Bigr)^2} \\
  &\quad\times
   \begin{bmatrix}
    - n^{-1} \ee^{-Bn} \sum\limits_{k=1}^n X_{k-1} & 1 \\[1mm]
    \ee^{-2Bn} \sum\limits_{k=1}^n X_{k-1}^2
     & - \ee^{-Bn}\sum\limits_{k=1}^n X_{k-1}
   \end{bmatrix}
   \begin{bmatrix}
    \ee^{-Bn/2} \sum\limits_{k=1}^n M_k \\[1mm]
    \ee^{-3Bn/2} \sum\limits_{k=1}^n M_k X_{k-1}
   \end{bmatrix}
 \end{align*}
 on the set \ $H_n$.
By Theorem \ref{w=0}, we have $\PP(w_{X_0}>0)=1$, and hence, by Theorem \ref{main_Ad0}, we get
 \begin{align}\label{help2}
   \frac{1}{\ee^{-2Bn} \sum\limits_{k=1}^n X_{k-1}^2
              - n^{-1} \Bigl(\ee^{-Bn} \sum\limits_{k=1}^n X_{k-1} \Bigr)^2}
    \asP \frac{\ee^{2B}-1}{w_{X_0}^2}
    \qquad \text{as \ $n \to \infty$,}
 \end{align}
 and
 \[
   \begin{bmatrix}
    - n^{-1} \ee^{-Bn} \sum\limits_{k=1}^n X_{k-1} & 1 \\[1mm]
    \ee^{-2Bn} \sum\limits_{k=1}^n X_{k-1}^2
     & - \ee^{-Bn}\sum\limits_{k=1}^n X_{k-1}
   \end{bmatrix}
   \asP \begin{bmatrix}
        0 & 1 \\[1mm]
        \frac{w_{X_0^2}}{\ee^{2B}-1}
        & - \frac{w_{X_0}}{\ee^B-1}
       \end{bmatrix} \qquad
   \text{as \ $n \to \infty$.}
 \]
Thus, using Theorem 3.18 in H\"ausler and Luschgy \cite{HauLus} (see also Theorem \ref{Thm_HL_Thm3_18}),
 the fact that \ $w_{X_0}$ \ is \ $\cF_\infty^X$-measurable (see Step 1 of the proof of Theorem \ref{main_Ad}),
 \ $\PP(H_n)\to 1$ \ as \ $n\to\infty$, \ and what is explained before Theorem \ref{main},
 the \ $\cF_\infty^X$-stable convergence \eqref{U} yields that
 \[
   \begin{bmatrix}
    \ee^{Bn/2} (\hvarrho_n - \varrho) \\[1mm]
    n \ee^{-Bn/2} (\hcA_n - \cA)
   \end{bmatrix}
   \to \frac{\ee^{2B}-1}{w_{X_0}^2}
       \begin{bmatrix}
        0 & 1 \\[1mm]
        \frac{w_{X_0^2}}{\ee^{2B}-1}
        & - \frac{w_{X_0}}{\ee^B-1}
       \end{bmatrix} \bR^{1/2}\tbN
  \qquad \text{$\cF_\infty^X$-stably as \ $n \to \infty$,}
 \]
 where \ $\tbN$ \ is a 2-dimensional random vector $\PP$-independent of \ $\cF_\infty^X$ \ such that \ $\tbN \distre \cN_2(0, V \bI_2)$
 \ and the random matrix \ $\bR$ \ is given in Theorem \ref{main_Ad}.
This implies \eqref{hvarrhohobeta}, since
 \[
  \frac{\ee^{2B}-1}{w_{X_0}^2}
       \begin{bmatrix}
        0 & 1 \\[1mm]
        \frac{w_{X_0^2}}{\ee^{2B}-1}
        & - \frac{w_{X_0}}{\ee^B-1}
       \end{bmatrix} \bR^{1/2}\tbN
    \distre \tbS^{1/2} \tbN,
 \]
 following from
  \begin{align*}
     \frac{(\ee^{2B}-1)^2}{w_{X_0}^4}
      \begin{bmatrix}
        0 & 1 \\[1mm]
        \frac{w_{X_0^2}}{\ee^{2B}-1}
        & - \frac{w_{X_0}}{\ee^B-1}
       \end{bmatrix}
      & \bR
       \begin{bmatrix}
        0 & \frac{w_{X_0^2}}{\ee^{2B}-1} \\[1mm]
        1 & - \frac{w_{X_0}}{\ee^B-1}
       \end{bmatrix} =
   \end{align*}
  \begin{align*}
     &= \frac{(\ee^{2B}-1)^2}{w_{X_0}^4}
        \begin{bmatrix}
         \frac{w_{X_0}^2}{\ee^{2B}-1}  & \frac{w_{X_0}^3}{\ee^{3B}-1} \\[1mm]
          0 & -\frac{\ee^B (\ee^B-1)w_{X_0}^4}{(\ee^{2B}-1)^2(\ee^{3B}-1)} \\
        \end{bmatrix}
        \begin{bmatrix}
        0 & \frac{w_{X_0^2}}{\ee^{2B}-1} \\[1mm]
        1 & - \frac{w_{X_0}}{\ee^B-1}
       \end{bmatrix}\\
  & = \frac{(\ee^{2B}-1)^2}{w_{X_0}^4}
       \begin{bmatrix}
         \frac{w_{X_0}^3}{\ee^{3B}-1}  & -\frac{\ee^B(\ee^B-1)}{(\ee^{2B}-1)^2(\ee^{3B}-1)}w_{X_0}^4 \\[1mm]
          -\frac{\ee^B(\ee^B-1)}{(\ee^{2B}-1)^2(\ee^{3B}-1)}w_{X_0}^4  & \frac{\ee^B}{(\ee^{2B}-1)^2(\ee^{3B}-1)}w_{X_0}^5 \\
        \end{bmatrix}
    = \tbS.
  \end{align*}

Under the additional assumption \ $C = 0$, \ we derive \eqref{hvarrhohobeta1} using Theorems \ref{main_Ad0} and \ref{main1_Ad}.
Using \eqref{Mk} and \eqref{CLSErb}, for each \ $n \in \NN$, \ we have
 \[
   \ee^{Bn} (\hvarrho_n - \varrho)
   = \frac{\ee^{-Bn}\sum\limits_{k=1}^n M_k X_{k-1}-n^{-1/2}\Bigl(n^{-1/2}\sum\limits_{k=1}^n M_k\Bigr)
            \Bigl(\ee^{-Bn}\sum\limits_{k=1}^nX_{k-1}\Bigr)}
          {\ee^{-2Bn} \sum\limits_{k=1}^n X_{k-1}^2
              - n^{-1} \Bigl(\ee^{-Bn} \sum\limits_{k=1}^n X_{k-1} \Bigr)^2}
 \]
 on the set \ $H_n$.
\ By part (i) of Theorem  \ref{main1_Ad}, we have \ $n^{-1/2}\sum_{k=1}^n M_k\to \widetilde N_1$ \ $\cF_\infty^X$-mixing as \ $n\to\infty$,
 \ where \ $\widetilde N_1$ \ is given in Theorem \ref{main1_Ad}; and, by Theorem \ref{main_Ad0},
 we get \ $\ee^{-Bn}\sum_{k=1}^nX_{k-1}\asP w_{X_0}/(\ee^B-1)$ \ as \ $n\to\infty$.
\ Consequently, using Slutsky's lemma and that mixing convergence yields convergence in distribution,
 we get
 \[
    n^{-1/2}\Bigl(n^{-1/2}\sum\limits_{k=1}^n M_k\Bigr)
   \Bigl(\ee^{-Bn}\sum\limits_{k=1}^nX_{k-1}\Bigr)
    \stoch 0 \qquad \text{as \ $n\to\infty$.}
 \]
Using part (ii) of Theorem \ref{main1_Ad}, \eqref{help2}, Theorem 3.18 in H\"ausler and Luschgy \cite{HauLus}
 (see also Theorem \ref{Thm_HL_Thm3_18}), the fact that \ $w_{X_0}$ \ is \ $\cF_\infty^X$-measurable (see Step 1 of the proof of Theorem \ref{main_Ad}),
 \ $\PP(H_n)\to 1$ \ as \ $n\to\infty$, \ and what is explained before Theorem \ref{main}, we obtain
 \[
    \ee^{Bn} (\hvarrho_n - \varrho)
      \to \frac{\ee^{2B}-1}{w_{X_0}^2}\cdot \frac{w_{X_0}}{\ee^B} \sum_{j=0}^\infty \ee^{-Bj}Z_j
      \qquad \text{$\cF_\infty^X$-stably as \ $n \to \infty$,}
 \]
 yielding \eqref{hvarrhohobeta1}.
Finally, using again \eqref{Mk} and \eqref{CLSErb}, for each \ $n \in \NN$, \ we have
 \[
   n^{1/2} (\hcA_n - \cA)
   = \frac{\Bigl(n^{-1/2}\sum\limits_{k=1}^n M_k\Bigr)
   \Bigl(\ee^{-2Bn}\sum\limits_{k=1}^nX_{k-1}^2\Bigr)
   -n^{-1/2} \Bigl(\ee^{-Bn}\sum\limits_{k=1}^n M_k X_{k-1}\Bigr) \Bigl(\ee^{-Bn}\sum\limits_{k=1}^nX_{k-1}\Bigr)}
          {\ee^{-2Bn} \sum\limits_{k=1}^n X_{k-1}^2
              - n^{-1} \Bigl(\ee^{-Bn} \sum\limits_{k=1}^n X_{k-1} \Bigr)^2}
 \]
 on the set \ $H_n$.
\ Similarly as above, by Theorems \ref{main_Ad0} and \ref{main1_Ad}, Theorem 3.18 in H\"ausler and Luschgy \cite{HauLus}
 (see also Theorem \ref{Thm_HL_Thm3_18}), the fact that \ $w_{X_0}$ \ is \ $\cF_\infty^X$-measurable (see Step 1 of the proof of Theorem \ref{main_Ad}),
 \ $\PP(H_n)\to 1$ \ as \ $n\to\infty$, \ and what is explained before Theorem \ref{main}, we obtain
 \begin{align*}
 n^{1/2} (\hcA_n - \cA)
    \to \frac{\ee^{2B}-1}{w_{X_0}^2}\left(\widetilde N_1 \cdot \frac{w_{X_0}^2}{\ee^{2B}-1} -0\right) = \widetilde N_1
         \qquad \text{$\cF_\infty^X$-stably as \ $n \to \infty$,}
 \end{align*}
 where \ $\tN_1$ \ is a random variable \ $\PP$-independent of \ $\cF_\infty^X$ \ such that
 \ $\tN_1 \distre \cN\bigl(0, \frac{\ee^{2B}-1}{2B} \int_0^\infty r^2 \, \nu(\dd r)\bigr)$.
\ Since \ $\tN_1$ \ is \ $\PP$-independent of \ $\cF_\infty^X$, \ we get the \ $\cF_\infty^X$-mixing convergence \eqref{hvarrhohobeta2}.

\section{Proof of Theorem \ref{main_Ad}}
\label{section_proof_main_Ad}

Let $\bB_0:=\bzero$, \ $\bU_0:=\bzero$, and for each \ $n \in \NN$, \ let
 \[
   \bU_n
   := \begin{bmatrix}
       U_n^{(1)} \\ U_n^{(2)}
      \end{bmatrix}
   := \begin{bmatrix}
       \sum_{k=1}^n M_k \\
        \sum_{k=1}^n M_k X_{k-1}
      \end{bmatrix}
 \]
 and
 \[
   \bB_n
   := \begin{bmatrix}
      \big(\sum_{k=1}^n (X_{k-1}+1) \big)^{-1/2} & 0 \\
      0 & \big(\sum_{k=1}^n (X_{k-1}^3+1)\big)^{-1/2} \\
    \end{bmatrix}.
 \]
Further, let
 \[
  \bQ_n
   := \begin{bmatrix}
       \ee^{-Bn/2} & 0 \\
       0 & \ee^{-3Bn/2}
      \end{bmatrix},\qquad n\in\ZZ_+.
 \]

We are going to apply Theorem \ref{MSLTES} with $d:=2$, $(\bU_n)_{n\in\ZZ_+}$, $(\bB_n)_{n\in\ZZ_+}$ and $(\bQ_n)_{n\in\ZZ_+}$ are given above,
  \ $(\cF_n)_{n\in\ZZ_+}:=(\sigma(X_0,\ldots,X_n))_{n\in\ZZ_+}$ and $G:=\Omega$
 \ (the matrices \ $\Beta$ \ and \ $\bP$, \ and the probability measure \ $\mu$ \ appearing
 in the assumptions (i), (iii) and (iv) of Theorem \ref{MSLTES} will be chosen later on).
Note that $G\in\cF_\infty(=\cF_\infty^X)$, $\PP(G)=1>0$, \ $\bB_n$ \ and \ $\bQ_n$ \ are invertible for each $n\in\NN$,
 \ and \ $(\bU_n)_{n\in\ZZ_+}$ \ and \ $(\bB_n)_{n\in\ZZ_+}$ \ are adapted to the filtration \ $(\cF_n)_{n\in\ZZ_+}$.

The remaining part of the proof is divided into five steps.

{\sl Step 1 (checking condition (i) of Theorem \ref{MSLTES}):}
For each $n\in\NN$, we have
 \[
   \bQ_n\bB_n^{-1}
     =  \begin{bmatrix}
          \ee^{-Bn/2}\big(\sum_{k=1}^n (X_{k-1}+1) \big)^{1/2} & 0 \\
          0 & \ee^{-3Bn/2}\big(\sum_{k=1}^n (X_{k-1}^3+1)\big)^{1/2} \\
        \end{bmatrix}.
 \]
By Theorems \ref{main_Ad0} and \ref{w=0}, we have
 \begin{align}\label{help_M_1}
   \ee^{-Bn} \sum_{k=1}^n X_{k-1} \asP \frac{w_{X_0}}{\ee^B-1}  \qquad
   \text{as \ $n \to \infty$,}
 \end{align}
 and
 \begin{align}\label{help_M_2}
   \ee^{-3Bn} \sum_{k=1}^n X_{k-1}^3 \asP \frac{w_{X_0}^3}{\ee^{3B}-1}  \qquad
   \text{as \ $n \to \infty$,}
 \end{align}
 where \ $\PP(w_{X_0}>0)=1$.
\ Hence, using also that $n\ee^{-Bn}\to 0$ and $n\ee^{-3Bn}\to 0$ as $n\to\infty$, we get
 \[
   \bQ_n\bB_n^{-1}
      \asP \begin{bmatrix}
          \big( \frac{w_{X_0}}{\ee^B-1} \big)^{1/2} & 0 \\
          0 & \big( \frac{w_{X_0}^3}{\ee^{3B}-1} \big)^{1/2} \\
        \end{bmatrix}
   \qquad \text{as \ $n \to \infty$.}
 \]
Consequently, since almost sure convergence yields convergence in probability, we obtain that condition
 (i) of Theorem \ref{MSLTES} holds with
 \[
   \Beta:= \begin{bmatrix}
          \big( \frac{w_{X_0}}{\ee^B-1} \big)^{1/2} & 0 \\
          0 & \big( \frac{w_{X_0}^3}{\ee^{3B}-1} \big)^{1/2} \\
        \end{bmatrix}.
 \]
Indeed, the random  matrix \ $\Beta$ \ is invertible if and only if \ $w_{X_0}\ne0$, \ so
 \ $\PP(G\cap\{\exists\; \Beta^{-1}\}) = \PP(w_{X_0}\ne 0)=1>0$ \ holds.
\ Further, \ $\Beta$ \ is \ $\cF_\infty$-measurable, since \ $w_{X_0}$ \ is \ $\cF_\infty$-mesaurable
 that can be checked as follows.
By \eqref{help_M_1}, we have $\ee^{-Bn} \sum_{k=1}^n X_{k-1}\asP \frac{w_{X_0}}{\ee^B-1}$ \
 as \ $n\to\infty$, \
 \ $\ee^{-Bn} \sum_{k=1}^n X_{k-1}$ \ is \ $\cF_\infty$-measurable for all $n\in\NN$,
 and the probability space \ $(\Omega,\cF_\infty^X ,\PP)$ \ is supposed to be complete (see, e.g., Cohn \cite[Corollary 2.2.3]{Coh}).

{\sl Step 2 (checking condition (ii) of Theorem \ref{MSLTES}):}
Condition (ii) of Theorem \ref{MSLTES} holds if and only if \ $\bigl(\ee^{-Bn/2} U_n^{(1)}\bigr)_{n\in\NN}$
 \ and \ $\bigl(\ee^{-3Bn/2} U_n^{(2)}\bigr)_{n\in\NN}$ \ are stochastically bounded in \ $\PP$-probability.
Indeed, \ $\PP(G\cap\{\exists \,\Beta^{-1}\})=1$ \ yields that \ $\PP_{G\cap\{\exists\,\Beta^{-1}\}}=\PP$ \ and for all \ $K>0$, \ we have
 \begin{align*}
    \{\bx=(x_1,x_2)\in\RR^2 : \vert x_j\vert >K\}
       \subset \{ \bx\in\RR^2 : \Vert \bx\Vert > K \},\qquad j=1,2,
 \end{align*}
 and
 \begin{align*}
   \{ \bx\in\RR^2 : \Vert \bx\Vert > K \}
      \subset \Big\{\bx=(x_1,x_2)\in\RR^2 : \vert x_1\vert > \frac{K}{\sqrt{2}}\Big\}
          \cup \Big\{\bx=(x_1,x_2)\in\RR^2 : \vert x_2\vert > \frac{K}{\sqrt{2}} \Big\}.
 \end{align*}
The process \ $(U_n^{(1)})_{n\in\ZZ_+}$ \ is a square integrable martingale
 with respect to the filtration \ $(\cF_n)_{n\in\ZZ_+}$ \ and it has quadratic characteristic
 \begin{align}\label{help13}
   \langle U^{(1)} \rangle_n
   = \sum_{k=1}^n \EE(M_k^2 \mid \cF_{k-1})
   = V \sum_{k=1}^n X_{k-1} + n V_0 , \quad n \in \NN ,
   \quad \text{and} \quad \langle U^{(1)} \rangle_0:=0,
 \end{align}
 where the second equality follows by Proposition \ref{moment_formula_2}.
For each \ $n \in \NN$ \ and \ $K \in (0, \infty)$, \ by Corollary \ref{Lenglart_Y} (Lenglart's inequality) with $a:=K^2 \ee^{Bn}$
 and $b:=K\ee^{Bn}$, we get
 \[
   \PP(\ee^{-Bn/2} |U_n^{(1)}| \geq K)
   = \PP((U_n^{(1)})^2 \geq K^2 \ee^{Bn})
   \leq \frac{1}{K}
        + \PP(\langle U^{(1)} \rangle_n \geq K \ee^{Bn}) ,
 \]
 and hence for all \ $K \in (0, \infty)$, \ we have
 \[
   \sup_{n\in\NN} \PP(\ee^{-Bn/2} |U_n^{(1)}| \geq K)
   \leq \frac{1}{K}
        + \sup_{n\in\NN} \PP(\ee^{-Bn} \langle U^{(1)} \rangle_n \geq K) .
 \]
By \eqref{help13} and Theorem \ref{main_Ad0}, we get
 \[
   \ee^{-Bn} \langle U^{(1)} \rangle_n
   = V \ee^{-Bn} \sum_{k=1}^n X_{k-1} + n \ee^{-Bn} V_0
      \asP \frac{V}{\ee^B-1} w_{X_0}
   \qquad \text{as \ $n \to \infty$,}
 \]
 and hence \ $\bigl(\ee^{-Bn} \langle U^{(1)} \rangle_n\bigr)_{n\in\NN}$ \ is stochastically bounded in $\PP$-probability,
 i.e.,
 \[
  \lim_{K\to\infty} \sup_{n\in\NN} \PP(\ee^{-Bn} \langle U^{(1)} \rangle_n \geq K) = 0.
 \]
Consequently, \ $\lim_{K\to\infty} \sup_{n\in\NN} \PP(\ee^{-Bn/2} |U_n^{(1)}| \geq K) = 0$, \ i.e.,
 \ $\bigl(\ee^{-Bn/2} U_n^{(1)}\bigr)_{n\in\NN}$ \ is stochastically bounded in $\PP$-probability, as desired.
In a similar way, the process \ $(U_n^{(2)})_{n\in\ZZ_+}$ \ is a square integrable martingale
 with respect to the filtration \ $(\cF_n)_{n\in\ZZ_+}$ \ and it has quadratic characteristic
 \begin{align}\label{help14}
  \begin{split}
   \langle U^{(2)} \rangle_n
   & = \sum_{k=1}^n \EE(M_k^2 X_{k-1}^2 \mid \cF_{k-1})
    = \sum_{k=1}^n X_{k-1}^2 \EE(M_k^2 \mid \cF_{k-1})\\
   & = V \sum_{k=1}^n X_{k-1}^3 + V_0 \sum_{k=1}^n X_{k-1}^2 , \qquad n \in \NN ,
  \end{split}
 \end{align}
 and \ $\langle U^{(2)} \rangle_0:=0$, \ where we used again Proposition \ref{moment_formula_2}.
For each \ $n \in \NN$ \ and \ $K \in (0, \infty)$, \ by Corollary \ref{Lenglart_Y} (Lenglart's inequality)
 with $a:=K^2 \ee^{3Bn}$ and $b:=K\ee^{3Bn}$, we get
 \[
   \PP(\ee^{-3Bn/2} |U_n^{(2)}| \geq K)
   = \PP((U_n^{(2)})^2 \geq K^2 \ee^{3Bn})
   \leq \frac{1}{K}
        + \PP(\langle U^{(2)} \rangle_n \geq K \ee^{3Bn}) ,
 \]
 and hence for all \ $K \in (0, \infty)$, \ we have
 \[
   \sup_{n\in\NN} \PP(\ee^{-3Bn/2} |U_n^{(2)}| \geq K)
   \leq \frac{1}{K}
        + \sup_{n\in\NN} \PP(\ee^{-3Bn} \langle U^{(2)} \rangle_n \geq K) .
 \]
By \eqref{help14} and Theorem \ref{main_Ad0}, we get
 \begin{align*}
   \ee^{-3Bn} \langle U^{(2)} \rangle_n
   & = V \ee^{-3Bn} \sum_{k=1}^n X_{k-1}^3 + V_0 \ee^{-3Bn} \sum_{k=1}^n X_{k-1}^2\\
   & = V \ee^{-3Bn} \sum_{k=1}^n X_{k-1}^3 + V_0 \ee^{-Bn} \cdot \ee^{-2Bn} \sum_{k=1}^n X_{k-1}^2\\
   & \asP \frac{V}{\ee^{3B}-1} w_{X_0}^3 + V_0 \cdot 0 \cdot \frac{w_{X_0}^2}{\ee^{2B}-1}
    = \frac{V}{\ee^{3B}-1} w_{X_0}^3
   \qquad \text{as \ $n \to \infty$,}
 \end{align*}
 and hence \ $\bigl(\ee^{-3Bn} \langle U^{(2)} \rangle_n\bigr)_{n\in\NN}$ \ is stochastically bounded in $\PP$-probability,
 i.e.,
 \[
    \lim_{K\to\infty} \sup_{n\in\NN} \PP(\ee^{-3Bn} \langle U^{(2)} \rangle_n \geq K) = 0.
 \]
Consequently, $\lim_{K\to\infty} \sup_{n\in\NN} \PP(\ee^{-3Bn/2} |U_n^{(2)}| \geq K) = 0$, i.e.,
 $\bigl(\ee^{-3Bn/2} U_n^{(2)}\bigr)_{n\in\NN}$ is stochastically bounded in $\PP$-probability,
 as desired.
All in all, we conclude that condition (ii) of Theorem \ref{MSLTES} holds.

{\sl Step 3 (checking condition (iii) of Theorem \ref{MSLTES}):}
For each \ $n,r \in \NN$ \ with \ $r<n$, \ we have
 \begin{align*}
   &\bB_n \bB_{n-r}^{-1} = \\
   & = \diag\left( \Big(\sum_{k=1}^n (X_{k-1}+1) \Big)^{-1/2},
       \Big(\sum_{k=1}^n (X_{k-1}^3+1)\Big)^{-1/2} \right)\\
   &\phantom{=\;}
     \times \left(\diag\left(
      \Big(\sum_{k=1}^{n-r} (X_{k-1}+1) \Big)^{-1/2} ,
      \Big(\sum_{k=1}^{n-r} (X_{k-1}^3+1)\Big)^{-1/2} \right)\right)^{-1} \\
   & = \diag\left(\!
      \Big(\sum_{k=1}^n (X_{k-1}+1) \Big)^{-1/2} \Big(\sum_{k=1}^{n-r} (X_{k-1}+1) \Big)^{1/2} ,
       \Big(\sum_{k=1}^n (X_{k-1}^3+1)\Big)^{-1/2} \Big(\sum_{k=1}^{n-r} (X_{k-1}^3+1)\Big)^{1/2}
       \right) \\
   & = \diag\Bigg(
        \ee^{-Br/2} \Big(\ee^{-Bn}\sum_{k=1}^n (X_{k-1}+1) \Big)^{-1/2} \Big(\ee^{-B(n-r)}\sum_{k=1}^{n-r} (X_{k-1}+1) \Big)^{1/2},\\
   &\phantom{=\diag\Bigg( \;}
         \ee^{-3Br/2}\Big(\ee^{-3Bn}\sum_{k=1}^n (X_{k-1}^3+1)\Big)^{-1/2} \Big(\ee^{-3B(n-r)}\sum_{k=1}^{n-r} (X_{k-1}^3+1)\Big)^{1/2}
          \Bigg).
 \end{align*}
Hence, using \eqref{help_M_1}, \eqref{help_M_2}, $\lim_{n\to\infty} n\ee^{-Bn}=0$, $\lim_{n\to\infty} n\ee^{-3Bn}=0$, \
 and \ $\PP(w_{X_0}>0)=1$, \ for each \ $r\in\NN$, \ we get
 \[
   \bB_n \bB_{n-r}^{-1}
       \asP \diag\left( \ee^{-Br/2} ,\ee^{-3Br/2} \right)
            = \left( \diag\left( \ee^{-B/2} ,\ee^{-3B/2} \right) \right)^r
   \qquad \text{as \ $n\to\infty$.}
 \]
This implies that condition (iii) of Theorem \ref{MSLTES} holds with
 \[
   \bP
   := \begin{bmatrix}
      \ee^{-B/2} & 0 \\
      0 & \ee^{-3B/2}
     \end{bmatrix} ,
 \]
 where \ $\varrho(\bP)=\max(\ee^{-B/2},\ee^{-3B/2}) = \ee^{-B/2}<1$.

{\sl Step 4 (checking condition (iv) of Theorem \ref{MSLTES}):}
For all \ $\btheta = (\theta_1,\theta_2)^\top\in \RR^2$ \ and \ $n \in \NN$, \ we have
 \begin{align}\label{help9}
  \begin{split}
   &\EE(\ee^{\ii\langle\btheta,\bB_n\Delta\bU_n\rangle} \mid \cF_{n-1})\\
   &= \EE\Bigl(\exp\Bigl\{\ii\Big(\theta_1  \Big(\sum_{k=1}^n (X_{k-1}+1) \Big)^{-1/2}
             + \theta_2 \Big(\sum_{k=1}^n (X_{k-1}^3+1)\Big)^{-1/2} X_{n-1}\Big) M_n \Bigr\} \,\big|\, \cF_{n-1}\Bigr) .
  \end{split}
 \end{align}
For all \ $\theta \in \RR$, \ by \eqref{regr}, we have
 \[
   \EE(\ee^{\ii\theta M_n} \mid \cF_{n-1})
    = \EE(\ee^{\ii\theta(X_n - \varrho X_{n-1} - \cA)} \mid \cF_{n-1})
    = \ee^{-\ii\theta(\varrho X_{n-1} + \cA)} \EE(\ee^{\ii\theta X_n} \mid \cF_{n-1}) ,
 \]
 where \ $\EE(\ee^{\ii\theta X_n} \mid \cF_{n-1}) = \EE(\ee^{\ii\theta X_n} \mid X_{n-1})$ \ and
 \ $\EE(\ee^{\ii\theta X_n} \mid X_{n-1} = x) = \EE(\ee^{\ii\theta X_1} \mid X_0 = x)$ \ for all \ $n \in \NN$, \ $\theta \in \RR$ \ and \
 \ $x \in \RR_+$, \ since \ $(X_t)_{t\in\RR_+}$ \ is a time-homogeneous
 Markov process.
Hence, by Theorem \ref{CBI_exists}, we get
 \begin{align*}
  \EE(\ee^{\ii\theta X_n} \mid \cF_{n-1})
  = \exp\biggl\{X_{n-1} \psi(1, \ii\theta) + \int_0^1 F(\psi(s, \ii\theta)) \, \dd s\biggr\} ,
  \qquad \theta \in \RR ,
 \end{align*}
 and
 \begin{align}\label{help15}
  \EE(\ee^{\ii\theta X_1} \mid X_0 = 0)
  = \exp\biggl\{\int_0^1 F(\psi(s, \ii\theta)) \, \dd s\biggr\} ,
  \qquad \theta \in \RR ,
 \end{align}
 where $\psi$ and $F$ are given in Theorem \ref{CBI_exists}.
Let \ $(Y_t)_{t\in\RR_+}$ \ be a CBI process with parameters \ $(c, 0, b, 0, \mu)$ \ (i.e., it is a pure branching process, a CB process).
Note that the same function $\psi$ (appearing in Theorem \ref{CBI_exists}) corresponds to the processes $(X_t)_{t\in\RR_+}$ and $(Y_t)_{t\in\RR_+}$,
 since $\psi$ does not depend on the immigration mechanism $F$.
Hence, by Theorem \ref{CBI_exists}, we have
 \begin{align*}
   \EE(\ee^{\ii\theta Y_1} \mid Y_0 = 1)
   = \ee^{\psi(1, \ii\theta)} , \qquad \theta \in \RR .
 \end{align*}
By \eqref{EXcond} and \eqref{ttBbeta}, we also have \ $\EE(Y_1 \mid Y_0 = 1) = \ee^B = \varrho$ \
 and \ $\EE(X_1 \mid X_0 = 0) = A \int_0^1 \ee^{Bs} \, \dd s = \cA$, \ thus
 \begin{align}\label{M_n}
 \begin{split}
   \EE(\ee^{\ii\theta M_n} \mid \cF_{n-1})
   & = \exp\Big\{(\psi(1,\ii\theta) - \ii\theta \varrho)X_{n-1}\Big\}
       \exp\biggl\{ \int_0^1 F(\psi(s, \ii\theta)) \, \dd s - \ii\theta \cA \biggr\}\\
   & = (\kappa_Y(\theta))^{X_{n-1}} \kappa_X(\theta) , \qquad
       \theta \in \RR , \quad n \in \NN ,
 \end{split}
 \end{align}
 where
 \begin{align}\label{help_kappaX_kappaY}
   \kappa_Y(\theta) := \EE(\ee^{\ii\theta(Y_1-\varrho)} \mid Y_0 = 1) = \ee^{\psi(1, \ii\theta) - \ii \theta \varrho} , \quad
   \kappa_X(\theta) := \EE(\ee^{\ii\theta(X_1-\cA)} \mid X_0 = 0) , \quad \theta \in \RR .
 \end{align}
Here we call the attention that \ $(\kappa_Y(\theta))^{X_{n-1}}$, \ the nonnegative $X_{n-1}$-th power of the conditional characteristic function
 of $Y_1-\varrho$ given $Y_0=1$ at the point \ $\theta$, \ is meant in the sense of so-called distinguished $X_{n-1}$-th power of $\kappa_Y(\theta)$,
 see, e.g., Sato \cite[Lemma 7.6 and the paragrapgh after it]{Sat}.
Indeed, \ $\kappa_Y(\theta)\ne 0$, \ $\theta\in\RR$ \ (following from $\psi(1,\ii\theta)\in\CC_{--}$) and since the function \ $\RR\ni\theta \mapsto \psi(1,\ii\theta) -\ii\theta \varrho$ \
 is continuous (following from Theorem \ref{CBI_exists}) and \ $\psi(1,\ii\cdot 0) -\ii\cdot 0 \varrho =0$, \
 we get that the function \ $\RR\ni\theta \mapsto \psi(1,\ii\theta) -\ii\theta \varrho$ \ coincides with the unique continuous function
 \ $f:\RR\to\CC$ \ satisfying \ $f(0)=0$ \ and \ $\ee^{f(\theta)} = \kappa_Y(\theta)$, \ $\theta\in\RR$.
\ That is, the function \ $\RR\ni\theta \mapsto \psi(1,\ii\theta) -\ii\theta \varrho$ \ is nothing else but the distinguished
 logarithm of \ $\kappa_Y$.

Consequently, using \eqref{help9} and that $\theta_1  \big(\sum_{k=1}^n (X_{k-1}+1) \big)^{-1/2}
 + \theta_2 \big(\sum_{k=1}^n (X_{k-1}^3+1)\big)^{-1/2} X_{n-1}$ is $\cF_{n-1}$-measurable,
 by the properties of conditional expectation, we get
 \begin{align}\label{help6}
  \begin{split}
  & \EE(\ee^{\ii\langle\btheta,\bB_n\Delta\bU_n\rangle} \mid \cF_{n-1})\\
  &\qquad = \Big(\kappa_Y\Big( \theta_1  \Big(\sum_{k=1}^n (X_{k-1}+1) \Big)^{-1/2}
                           + \theta_2 \Big(\sum_{k=1}^n (X_{k-1}^3+1)\Big)^{-1/2} X_{n-1}
                     \Big)\Big)^{X_{n-1}} \\
  &\phantom{\qquad=\;} \times \kappa_X\Big( \theta_1  \Big(\sum_{k=1}^n (X_{k-1}+1) \Big)^{-1/2}
                           + \theta_2 \Big(\sum_{k=1}^n (X_{k-1}^3+1)\Big)^{-1/2} X_{n-1} \Big)
 \end{split}
 \end{align}
 for all \ $\btheta \in \RR^2$ \ and \ $n \in \NN$.
\ Next, we verify that
 \begin{align}\label{help4}
   \kappa_X\Big( \theta_1  \Big(\sum_{k=1}^n (X_{k-1}+1) \Big)^{-1/2}
                           + \theta_2 \Big(\sum_{k=1}^n (X_{k-1}^3+1)\Big)^{-1/2} X_{n-1} \Big)
     \asP 1 \qquad \text{as \ $n \to \infty$,}
 \end{align}
 and
 \begin{align}\label{help17}
  \begin{split}
   &\Big(\kappa_Y\Big( \theta_1  \Big(\sum_{k=1}^n (X_{k-1}+1) \Big)^{-1/2}
                           + \theta_2 \Big(\sum_{k=1}^n (X_{k-1}^3+1)\Big)^{-1/2} X_{n-1}
                     \Big)\Big)^{X_{n-1}}\\
   &\qquad  \asP \exp\left\{-\frac{V}{2} \left( \theta_1 \left(\frac{\ee^B -1}{\ee^B}\right)^{1/2}
                                                + \theta_2 \left(\frac{\ee^{3B} -1}{\ee^{3B}}\right)^{1/2}  \right)^2\right\}
    \qquad \text{as \ $n \to \infty$.}
  \end{split}
 \end{align}
Note that for all \ $\btheta \in \RR^2$, \ using Theorem \ref{convX}, \eqref{help_M_1}, \eqref{help_M_2},
 \ $\lim_{n\to\infty} n\ee^{-Bn} = 0$ \ and that \ $\PP(w_{X_0}>0)=1$, \ we have
 \begin{align*}
  &\theta_1 \Big(\sum_{k=1}^n (X_{k-1}+1) \Big)^{-1/2}
        + \theta_2 \Big(\sum_{k=1}^n (X_{k-1}^3+1)\Big)^{-1/2} X_{n-1} \\
  &\qquad  =   \theta_1 \ee^{-Bn/2}\Big( \ee^{-Bn} \sum_{k=1}^n (X_{k-1}+1) \Big)^{-1/2}\\
  &\phantom{\qquad =\; } + \theta_2 \ee^{-B(n+2)/2} \Big(\ee^{-3Bn} \sum_{k=1}^n (X_{k-1}^3+1)\Big)^{-1/2} (\ee^{-B(n-1)} X_{n-1})
 \end{align*}
 \begin{align*}
  &\phantom{\qquad } \asP \theta_1 \cdot 0\cdot (w_{X_0}/(\ee^{B}-1))^{-1/2}
                                + \theta_2 \cdot 0\cdot (w_{X_0}^3/(\ee^{3B}-1))^{-1/2} w_{X_0} =0
 \end{align*}
 as \ $n \to \infty$.
\ Thus, using that the conditional characteristic function \ $\kappa_X$ \ of \ $X_1-\cA$ \ given \ $X_0=0$ \ is continuous,
 we get \eqref{help4}.

Now, we turn to verify \eqref{help17}.
Recall that \ $\kappa_Y$ \ is the conditional characteristic function of \ $Y_1-\varrho$ \ given \ $Y_0=1$ \ (see \eqref{help_kappaX_kappaY}).
Since $\EE(Y_1\mid Y_0=1)=\varrho$, \ we have \ $\EE(Y_1 - \varrho \mid Y_0 = 1) = 0$, and,
 by Remark \ref{REMARK_par}, we get $\var(Y_1 - \varrho \mid Y_0 = 1) = \var(Y_1 \mid Y_0 = 1) = V$.
Hence the classical central limit theorem for sums of independent and identically distributed random variables
 and the continuity theorem yield that
 \begin{equation*}
   \left(\kappa_Y\left(\frac{\theta}{\sqrt{n}}\right)\right)^n
          = \left(\ee^{\psi(1,\ii\frac{\theta}{\sqrt{n}}) - \ii \frac{\theta}{\sqrt{n}} \varrho}\right)^n
          = \ee^{\big(\psi(1,\ii\frac{\theta}{\sqrt{n}}) - \ii \frac{\theta}{\sqrt{n}} \varrho\big)n}
   \to \ee^{-V\theta^2/2} \qquad \text{as \ $n \to \infty$}
 \end{equation*}
 uniformly in \ $\theta \in \RR$ \ on compact intervals,
 where, for each \ $n\in\NN$,
 \vspace*{-2mm}
 \begin{itemize}
   \item[(a)] the function
         \[
         \RR\ni\theta \mapsto  \ee^{\big(\psi(1,\ii\frac{\theta}{\sqrt{n}}) - \ii \frac{\theta}{\sqrt{n}} \varrho\big)n} =: \varphi_n(\theta)
         \]
         is the characteristic function of a random variable \ $\frac{1}{\sqrt{n}} \sum_{j=1}^n (\tY_j-\varrho)$,
         \ where \ $\tY_j-\varrho$, $j=1,\ldots,n$, \ are independent and identically distributed random variables
          with a common characteristic function \ $\RR\ni\theta\mapsto \ee^{\psi(1,\ii\theta) - \ii\theta \varrho}$,
   \item[(b)] the function \ $\RR\ni\theta\mapsto \ee^{-V\theta^2/2}=: \varphi(\theta)$ \ is the characteristic function of a normally distributed random
         variable with mean zero and variance \ $V$.
 \end{itemize}
Then we can apply Lemma 7.7 in Sato \cite{Sat} with the choices \ $\varphi_n$, $n\in\NN$, \ and \  $\varphi$ \ (given above in (a) and (b)),
 since \ $\varphi_n(0) = 1$, \ $n\in\NN$, \ $\varphi(0) =1$, \ $\varphi_n(\theta)\ne 0$, \ $\theta\in\RR$, \ $n\in\NN$ \
 (due to $\psi(1,\ii\frac{\theta}{\sqrt{n}})\in\CC_{--}$, $\theta\in\RR$, $n\in\NN$),
 \ $\varphi(\theta)\ne 0$, \ $\theta\in\RR$, \ and \ $\varphi_n$, $n\in\NN$, \ and \ $\varphi$ \ are continuous (in case of \ $\varphi_n$, \ it follows by Theorem \ref{CBI_exists}).
So Lemma 7.7 in Sato \cite{Sat} implies that the distinguished logarithm of \ $\varphi_n$ \ converges to the distinguished logarithm of \ $\varphi$ \ as \ $n\to\infty$ \
 uniformly on compact intervals.
Since, for each $n\in\NN$, the functions \ $\RR\ni\theta\mapsto \big(\psi(1,\ii\frac{\theta}{\sqrt{n}}) - \ii \frac{\theta}{\sqrt{n}} \varrho\big)n$
 \ and \ $\RR\ni\theta \mapsto -V\frac{\theta^2}{2}$ \ are continuous and take the value $0$ at $0$, using Lemma 7.6 in Sato \cite{Sat},
 we get that these two functions are the distinguished logarithm of $\varphi_n$ and $\varphi$, respectively.
Consequently, we have
 \[
   \left(\psi\left(1,\ii\frac{\theta}{\sqrt{n}}\right) - \ii \frac{\theta}{\sqrt{n}} \varrho\right)n \to -V\frac{\theta^2}{2}
    \qquad \text{as \ $n\to\infty$}
 \]
 uniformly in \ $\theta \in \RR$ \ on compact intervals.
Note that \ $X_n\asP \infty$ \ as $n\to\infty$.
Indeed, by Theorems \ref{convX} and \ref{w=0},
 we have \ $\ee^{-Bn}X_n\asP w_{X_0}$ \ as $n\to\infty$, where $\PP(w_{X_0}>0)=1$, and hence
 $X_n=\ee^{Bn}(\ee^{-Bn}X_n)\asP\infty\cdot w_{X_0}\aseP\infty$ as $n\to\infty$.
\ In particular, we also get \ $\lfloor X_n\rfloor\asP \infty$ \ and \ $\frac{\lfloor X_n\rfloor}{X_n}\asP 1$ \ as \ $n\to\infty$.
Consequently, we obtain that
 \begin{align}\label{help18}
 \left(\psi\left(1,\ii\frac{\theta}{\sqrt{\lfloor X_{n-1}\rfloor}}\right) - \ii \frac{\theta}{\sqrt{\lfloor X_{n-1}\rfloor}} \varrho\right)\lfloor X_{n-1}\rfloor
    \asP -V\frac{\theta^2}{2}
    \qquad \text{as \ $n\to\infty$}
 \end{align}
 uniformly in \ $\theta \in \RR$ \ on compact intervals.
Moreover, using Theorem \ref{convX}, \eqref{help_M_1}, \eqref{help_M_2} and that $\PP(w_{X_0}>0)=1$,
 \ for all \ $\btheta \in \RR^2$, \ we have that
 \begin{align*}
   &\xi_n \lfloor X_{n-1}\rfloor^{1/2}:= \Big(\theta_1 \Big(\sum_{k=1}^n (X_{k-1}+1) \Big)^{-1/2}
        + \theta_2 \Big(\sum_{k=1}^n (X_{k-1}^3+1)\Big)^{-1/2} X_{n-1} \Big) \lfloor X_{n-1}\rfloor^{1/2} \\
    & = \Bigg( \theta_1 \ee^{-B/2} \left(\frac{\ee^{-B(n-1)}X_{n-1}}{\ee^{-Bn}\sum_{k=1}^n (X_{k-1}+1)} \right)^{1/2}\\
    &\phantom{=\Bigg(\;} + \theta_2 \ee^{-3B/2} \left(\frac{\ee^{-3B(n-1)} X_{n-1}^3}{\ee^{-3Bn} \sum_{k=1}^n (X_{k-1}^3+1)} \right)^{1/2} \Bigg)
            \left(\frac{\lfloor X_{n-1}\rfloor}{X_{n-1}}\right)^{1/2}  \\
    &\asP \theta_1 \ee^{-B/2} \left(\frac{w_{X_0}}{w_{X_0}/(\ee^B -1)}\right)^{1/2}
          + \theta_2 \ee^{-3B/2} \left(\frac{w_{X_0}^3}{w_{X_0}^3/(\ee^{3B} -1)}\right)^{1/2}\\
    &\phantom{\asP}
        = \theta_1 \left(\frac{\ee^B -1}{\ee^B}\right)^{1/2}
          + \theta_2 \left(\frac{\ee^{3B} -1}{\ee^{3B}}\right)^{1/2} =: L
          \qquad \text{as \ $n\to\infty$.}
 \end{align*}
Consequently, using that the convergence in \eqref{help18} is uniform in \ $\theta \in \RR$ \ on compact intervals, we get that
 \[
   \left(\psi\left(1,\ii \xi_n \right) - \ii \xi_n \varrho\right)\lfloor X_{n-1}\rfloor \asP -V\frac{L^2}{2}
    \qquad \text{as \ $n\to\infty$.}
 \]
Hence, using that \ $\frac{\lfloor X_n\rfloor}{X_n}\asP 1$ \ as \ $n\to\infty$, \ we have that
 \begin{align*}
  \left(\psi\left(1,\ii \xi_n \right) - \ii \xi_n \varrho\right) X_{n-1}
    =  \left(\psi\left(1,\ii \xi_n \right) - \ii \xi_n \varrho\right)\lfloor X_{n-1}\rfloor \frac{X_{n-1}}{\lfloor X_{n-1}\rfloor}
    \asP -V\frac{L^2}{2}\qquad \text{as \ $n\to\infty$.}
 \end{align*}
Since the function \ $\exp$ \ is continuous on \ $\CC$, \ we get \eqref{help17}, as desired.

Using \eqref{help6}, \eqref{help4} and \eqref{help17}, for all $\btheta\in\RR^2$, we have
 \begin{align*}
  \EE(\ee^{\ii\langle\btheta,\bB_n\Delta\bU_n\rangle} \mid \cF_{n-1})
     &\asP \exp\left\{-\frac{V}{2} \left( \theta_1 \left(\frac{\ee^B -1}{\ee^B}\right)^{1/2}
                                           + \theta_2 \left(\frac{\ee^{3B} -1}{\ee^{3B}}\right)^{1/2}  \right)^2\right\}\\
     &\phantom{\asP\;} = \int_{\RR^2} \ee^{\ii \langle\btheta, \bx\rangle} \, \PP^{\bzeta}(\dd\bx)
     \qquad \text{as \ $n \to \infty$,}
 \end{align*}
  where \ $\bzeta$ \ is an \ $\RR^2$-valued random variable having a 2-dimensional normal distribution
 \[
  \cN_2\left(\bzero, V
                    \begin{bmatrix}
                     \frac{\ee^B -1}{\ee^B} & \frac{(\ee^B -1)^{1/2}(\ee^{3B} -1)^{1/2}}{\ee^{2B}}  \\[2mm]
                     \frac{(\ee^B -1)^{1/2}(\ee^{3B} -1)^{1/2}}{\ee^{2B}} & \frac{\ee^{3B} -1}{\ee^{3B}}
                    \end{bmatrix}\right).
 \]
Hence we obtain that condition (iv) of Theorem \ref{MSLTES} holds with \ $\mu := \PP^{\bzeta}$, \ since
 \begin{align*}
  \int_{\RR^2} \log^+(\Vert \bx\Vert) \,\mu(\dd\bx)
    & = \int_{\{ \bx\in \RR^2 : \Vert \bx\Vert\geq 1\}} \log(\Vert \bx\Vert) \,\mu(\dd\bx)
      \leq \int_{\{\bx\in \RR^2 : \Vert \bx\Vert\geq 1\}} \Vert \bx\Vert \,\mu(\dd\bx) \\
    & \leq \int_{\{\bx=(x_1,x_2)\in \RR^2 : \Vert \bx\Vert\geq 1\}} (\vert x_1\vert + \vert x_2\vert) \,\mu(\dd\bx)
      <\infty
 \end{align*}
 due to the fact that all the mixed moments of \ $\mu$ \ (being a \ $2$-dimensional normal distribution) are finite.

{\sl Step 5 (application of Theorem \ref{MSLTES}):}
Using Steps 1--4, we can apply Theorem \ref{MSLTES} with our choices given before  Step 1 and in Steps 1, 3 and 4
 (for $\Beta$, $\bP$ and $\mu$), and we conclude that
 \begin{align}\label{help7}
   \bQ_n \bU_n \to \Beta\sum_{j=0}^\infty \bP^j \bZ_j \qquad \text{$\cF_\infty$-stably as \ $n \to \infty$,}
 \end{align}
 where $(\bZ_j)_{j\in\ZZ_+}$ is a sequence of $\PP$-independent and identically distributed $\RR^2$-valued
 random vectors being $\PP$-independent of $\cF_\infty$ (in particular of $\Beta$, since $\Beta$ is $\cF_\infty$-measurable, see Step 1)
 with $\PP(\bZ_0 \in B) = \mu(B)$ for all $B\in\cB(\RR^2)$,
 and the series in \eqref{help7} is absolutely convergent $\PP$-almost surely.
The distribution of the limit random variable in \eqref{help7} can be written in the form $\Beta\bZ$, where
 $\bZ$ is a 2-dimensional random vector \ $\PP$-independent of $\cF_\infty$ (in particular of $\Beta$) such that
 $\bZ$ has a 2-dimensional normal distribution with mean vector $\bzero\in\RR^2$ and with covariance matrix
 \begin{align*}
  &V\sum_{j=0}^\infty
      \bP^j \begin{bmatrix}
                     \frac{\ee^B -1}{\ee^B} & \frac{(\ee^B -1)^{1/2}(\ee^{3B} -1)^{1/2}}{\ee^{2B}}  \\[2mm]
                     \frac{(\ee^B -1)^{1/2}(\ee^{3B} -1)^{1/2}}{\ee^{2B}} & \frac{\ee^{3B} -1}{\ee^{3B}}
                    \end{bmatrix}
            (\bP^j)^\top\\[1mm]
  &\qquad = V\sum_{j=0}^\infty
      \begin{bmatrix}
                     \frac{\ee^B -1}{\ee^B} \,\ee^{-Bj} & \frac{(\ee^B -1)^{1/2}(\ee^{3B} -1)^{1/2}}{\ee^{2B}}\, \ee^{-2Bj}  \\[2mm]
                     \frac{(\ee^B -1)^{1/2}(\ee^{3B} -1)^{1/2}}{\ee^{2B}}\, \ee^{-2Bj} & \frac{\ee^{3B} -1}{\ee^{3B}}\, \ee^{-3Bj}
                    \end{bmatrix} \\[1mm]
 &\qquad =V \begin{bmatrix}
       1 & \left( \frac{(\ee^B -1)(\ee^{3B} -1)}{(\ee^{2B}-1)^2} \right)^{1/2} \\[2mm]
       \left( \frac{(\ee^B -1)(\ee^{3B} -1)}{(\ee^{2B}-1)^2} \right)^{1/2} & 1
      \end{bmatrix}.
 \end{align*}
Since
 \begin{align*}
    &\begin{bmatrix}
       1 & \left( \frac{(\ee^B -1)(\ee^{3B} -1)}{(\ee^{2B}-1)^2} \right)^{1/2} \\[2mm]
       \left( \frac{(\ee^B -1)(\ee^{3B} -1)}{(\ee^{2B}-1)^2} \right)^{1/2} & 1
      \end{bmatrix} \\
   &\qquad = \begin{bmatrix}
         1 & 0 \\[2mm]
        \left( \frac{(\ee^B -1)(\ee^{3B} -1)}{(\ee^{2B}-1)^2} \right)^{1/2} & \frac{\ee^{B/2}(\ee^B -1)}{\ee^{2B}-1}
     \end{bmatrix}
     \begin{bmatrix}
         1 & 0 \\[2mm]
        \left( \frac{(\ee^B -1)(\ee^{3B} -1)}{(\ee^{2B}-1)^2} \right)^{1/2} & \frac{\ee^{B/2}(\ee^B -1)}{\ee^{2B}-1}
      \end{bmatrix}^\top,
 \end{align*}
 using also the formula for \ $\Beta$ \ (see Step 1),
 we get that the distribution of $\Beta\bZ$ coincides with the distribution of
 \begin{align*}
   &\begin{bmatrix}
          \big( \frac{w_{X_0}}{\ee^B-1} \big)^{1/2} & 0 \\
          0 & \big( \frac{w_{X_0}^3}{\ee^{3B}-1} \big)^{1/2} \\
   \end{bmatrix}
   \begin{bmatrix}
         1 & 0 \\[2mm]
        \left( \frac{(\ee^B -1)(\ee^{3B} -1)}{(\ee^{2B}-1)^2} \right)^{1/2} & \frac{\ee^{B/2}(\ee^B -1)}{\ee^{2B}-1}
     \end{bmatrix}
   \tbN  =
 \end{align*}
 \begin{align*}
 &= \begin{bmatrix}
      \frac{1}{(\ee^B-1)^{1/2}}\,w_{X_0}^{1/2} & 0 \\[2mm]
       \frac{(\ee^B-1)^{1/2}}{\ee^{2B}-1}\,w_{X_0}^{3/2}  &  \frac{\ee^{B/2}(\ee^B-1)}{(\ee^{2B}-1)(\ee^{3B}-1)^{1/2}} \,w_{X_0}^{3/2}
    \end{bmatrix}
   \tbN
 = \bkappa\tbN  ,
 \end{align*}
 where \ $\tbN$ \ is given in Theorem \ref{main_Ad} and $\bkappa$ is given in \eqref{kappa_2}.
By Remark \ref{REMARK4}, we have \ $\bR^{1/2}\tbN \distre \bkappa\tbN$, \ and hence we conclude the statement of Theorem \ref{main_Ad}.

\section{Proof of Theorem \ref{main1_Ad}}
\label{section_proof_main1_Ad}

The \ $\cF_\infty^X$-mixing convergence in part (i) follows from a mixing limit theorem for an ergodic sequence of martingale differences,
 see, e.g., Corollary 6.26 in H\"ausler and Luschgy \cite{HauLus} (see also Theorem \ref{Thm_HL_Cor6_26}).
Indeed, \ $(M_k)_{k\in\NN}$ \ is a sequence of martingale differences with respect to the filtration
 \ $(\cF_k^X)_{k\in\ZZ_+}$.
Further, by Lemma \ref{SDE_transform_sol_1t}, the sequence \ $(M_k)_{k\in\NN}$ \ consists of \ $\PP$-independent
 and identically distributed random variables yielding that \ $(M_k)_{k\in\NN}$ \ is stationary and any event belonging
 to the tail $\sigma$-algebra corresponding to $(M_k)_{k\in\NN}$ has probability 0 or 1 (due to Kolmogorov's 0 or 1 law).
Since for any $D\in\cI_S$ (where $\cI_S$ is defined in Appendix \ref{App_1} before Theorem \ref{Thm_HL_Cor6_26}),
 the set \ $\{\omega\in\Omega : (M_k(\omega))_{k\in\NN}\in D \}$ \ belongs to the tail-$\sigma$-field in question
 (see, e.g., H\"ausler and Luschgy \cite[page 57]{HauLus}), we get that \ $\PP^{(M_k)_{k\in\NN}}(D)\in\{0,1\}$ \
 for all \ $D\in \cI_S$, \ i.e., \ $(M_k)_{k\in\NN}$ \ is ergodic as well.
Moreover, we have \ $\EE(M_1) = 0$, \ and, by the assumption $C=0$ and Proposition \ref{moment_formula_2},
 we also have \ $V=0$ \ and
 \[
    \EE(M_k^2\mid \cF_{k-1}^X)
     = \var(M_k\mid \cF_{k-1}^X)
     = V_0
     = \int_0^\infty r^2 \, \nu(\dd r) \int_0^1 \ee^{2Bu} \, \dd u
     = \frac{\ee^{2B}-1}{2B} \int_0^\infty r^2 \, \nu(\dd r)
 \]
 for \ $k\in\NN$, \ yielding that \ $\EE(M_1^2) = \var(M_1) =  \frac{\ee^{2B}-1}{2B} \int_0^\infty r^2 \, \nu(\dd r)<\infty$.
\ All in all, one can apply Corollary 6.26 in H\"ausler and Luschgy \cite{HauLus} (see also Theorem \ref{Thm_HL_Cor6_26}).

The \ $\cF_\infty^X$-stable convergence in part (ii) can be shown similarly as the \ $\cF_\infty^X$-stable convergence in Theorem \ref{main_Ad}.
For each \ $n \in \NN$, \ we have \ $\ee^{-Bn} \sum_{k=1}^n M_k X_{k-1} = Q_n U_n$, \ where
 \ $Q_n := \ee^{-Bn}$, \ $n \in \NN$, \ and \ $U_n := \sum_{k=1}^n M_k X_{k-1}$, \ $n \in \NN$ \ with $U_0:=0$.
Further, let $B_n:=\left(\sum_{k=1}^n (X_{k-1}+1)\right)^{-1}$, $n\in\NN$, with $B_0:=0$.
\ We are going to apply Theorem \ref{MSLTES} with $d:=1$, $(U_n)_{n\in\ZZ_+}$, $(B_n)_{n\in\ZZ_+}$ and $(Q_n)_{n\in\ZZ_+}$ are given above,
  \ $(\cF_n)_{n\in\ZZ_+}:=(\sigma(X_0,\ldots,X_n))_{n\in\ZZ_+}$ and $G:=\Omega$
 \ (and \ $\eta$, \ $P$ \ and the probability measure \ $\mu$ \ appearing in the assumptions (i), (iii) and (iv) of Theorem \ref{MSLTES} will be chosen later on).
Since we apply Theorem \ref{main_Ad} in case of dimension $1$, we do not use boldface notation for $Q_n, B_n, U_n,\eta$ and $P$ for simplicity.
Note that $G\in\cF_\infty(=\cF_\infty^X)$, $\PP(G)=1>0$, \ $B_n$ and $Q_n$ are non-zero (invertible) for each $n\in\NN$,
 \ and \ $(U_n)_{n\in\ZZ_+}$ \ and \ $(B_n)_{n\in\ZZ_+}$ \ are adapted to the filtration \ $(\cF_n)_{n\in\ZZ_+}$.

The remaining part of the proof is divided into five steps.

{\sl Step 1 (checking condition (i) of Theorem \ref{MSLTES}):}
For each $n\in\NN$, using Theorem \ref{main_Ad0} and that \ $n\ee^{-Bn}\to 0$ \ as \ $n\to\infty$,
 \ we have
 \begin{align*}
   Q_nB_n^{-1}
     = \ee^{-Bn} \sum_{k=1}^n (X_{k-1}+1) \asP \frac{w_{X_0}}{\ee^B -1}
       \qquad \text{as \ $n\to\infty$.}
 \end{align*}
Hence, since almost sure convergence yields convergence in probability, we obtain that
 condition (i) of Theorem \ref{MSLTES} holds with
 \[
   \eta:= \frac{w_{X_0}}{\ee^B -1}.
 \]
Indeed, \ $\eta$ \ is non-zero (invertible) if and only if \ $w_{X_0}\ne0$, \ so
 \ $\PP(G\cap\{\exists\; \eta^{-1}\}) = \PP(w_{X_0}\ne 0)=1>0$ \ holds.
\ Further, \ $\eta$ \ is \ $\cF_\infty$-measurable, since \ $w_{X_0}$ \ is \ $\cF_\infty$-measurable
 that can be checked as follows.
By Theorem \ref{main_Ad0}, we have $\ee^{-Bn} \sum_{k=1}^n X_{k-1}\asP\frac{w_{X_0}}{\ee^B-1}$ as $n\to\infty$,
 \ $\ee^{-Bn} \sum_{k=1}^n X_{k-1}$ is $\cF_\infty$-measurable for all $n\in\NN$,
 and the probability space $(\Omega,\cF_\infty^X ,\PP)$ is supposed to be complete (see, e.g., Cohn \cite[Corollary 2.2.3]{Coh}).

{\sl Step 2 (checking condition (ii) of Theorem \ref{MSLTES}):}
The process \ $(U_n)_{n\in\ZZ_+}$ \ is a square integrable martingale with respect to the filtration \ $(\cF_n)_{n\in\ZZ_+}$ \
 and it has quadratic characteristic
 \begin{align}\label{help16}
   \langle U\rangle_n
   = \sum_{k=1}^n \EE(M_k^2 X_{k-1}^2 \mid \cF_{k-1})
   = \sum_{k=1}^n X_{k-1}^2 \EE(M_k^2 \mid \cF_{k-1})
   = V_0 \sum_{k=1}^n X_{k-1}^2 , \qquad n \in \NN,
 \end{align}
 and $\langle U \rangle_0:=0$, where the last equality follows by Proposition \ref{moment_formula_2} with $C=0$.
For each $n \in \NN$ and $K \in (0, \infty)$, by Corollary \ref{Lenglart_Y} (Lenglart's inequality) with
 $a:=K^2\ee^{2Bn}$ and $b:=K\ee^{2Bn}$, we get
 \[
   \PP(\ee^{-Bn} |U_n| \geq K)
   = \PP(U_n^2 \geq K^2 \ee^{2Bn})
   \leq \frac{1}{K}
        + \PP(\langle U\rangle_n \geq K \ee^{2Bn}) ,
 \]
 and hence for all \ $K \in (0, \infty)$, \ we have
 \[
   \sup_{n\in\NN} \PP(\ee^{-Bn} |U_n| \geq K)
   \leq \frac{1}{K}
        + \sup_{n\in\NN} \PP(\ee^{-2Bn} \langle U\rangle_n \geq K) .
 \]
By \eqref{help16} and Theorem \ref{main_Ad0}, we get
 \[
   \ee^{-2Bn} \langle U\rangle_n
   = V_0 \ee^{-2Bn} \sum_{k=1}^n X_{k-1}^2
      \asP V_0\frac{w_{X_0}^2}{\ee^{2B}-1}
   \qquad \text{as \ $n \to \infty$,}
 \]
 and hence \ $\bigl(\ee^{-2Bn} \langle U\rangle_n\bigr)_{n\in\NN}$ \ is stochastically bounded in $\PP$-probability,
 i.e., it holds that \ $\lim_{K\to\infty} \sup_{n\in\NN} \PP(\ee^{-2Bn} \langle U\rangle_n \geq K) = 0$.
\ Consequently, \ $\lim_{K\to\infty} \sup_{n\in\NN} \PP(\ee^{-Bn} |U_n| \geq K) = 0$, \ i.e.,
 \ $\bigl(\ee^{-Bn} U_n\bigr)_{n\in\NN}$ \ is stochastically bounded in $\PP$-probability,
 and we conclude that condition (ii) of Theorem \ref{MSLTES} holds.

{\sl Step 3 (checking condition (iii) of Theorem \ref{MSLTES}):}
For each \ $n,r \in \NN$ \ with \ $n>r$, \ we have
 \begin{align*}
  B_n B_{n-r}^{-1}
  & = \left(\sum_{k=1}^n (X_{k-1}+1)\right)^{-1}
      \sum_{k=1}^{n-r} (X_{k-1}+1) \\
  & = \ee^{-Br}\left(\ee^{-Bn}\sum_{k=1}^n (X_{k-1}+1)\right)^{-1}
      \ee^{-B(n-r)} \sum_{k=1}^{n-r} (X_{k-1}+1).
 \end{align*}
Hence, using Theorem \ref{main_Ad0}, \ $\lim_{n\to\infty} n\ee^{-Bn}=0$, \ and
 \ $\PP(w_{X_0}>0)=1$, \ for each \ $r\in\NN$, \ we get
 \[
   B_n B_{n-r}^{-1} \asP \ee^{-Br} \left(\frac{w_{X_0}}{\ee^B-1}\right)^{-1}
                                   \frac{w_{X_0}}{\ee^B-1}
                     = (\ee^{-B})^r \qquad \text{as \ $n\to\infty$.}
 \]
This implies that condition (iii) of Theorem \ref{MSLTES} holds with \ $P := \ee^{-B}$, where $\varrho(P)=\ee^{-B}<1$.

{\sl Step 4 (checking condition (iv) of Theorem \ref{MSLTES}):}
For all \ $\theta \in \RR$ \ and \ $n \in \NN$, \ we have
 \begin{align}\label{help11}
   \EE(\ee^{\ii\theta B_n\Delta U_n} \mid \cF_{n-1})
   = \EE\left(\exp\left\{\ii\theta \left(\sum_{k=1}^n (X_{k-1}+1)\right)^{-1} M_n X_{n-1}\right\} \Bigg|\, \cF_{n-1}\right) .
 \end{align}
Recall that, by \eqref{M_n}, we get
 \[
   \EE(\ee^{\ii\theta M_n} \mid \cF_{n-1} )
        = \exp\Big\{(\psi(1,\ii\theta) - \ii\theta \varrho)X_{n-1}\Big\}
        \exp\biggl\{ \int_0^1 F(\psi(s, \ii\theta)) \, \dd s - \ii\theta \cA \biggr\}
 \]
 for all \ $\theta\in\RR$ \ and \ $n\in\NN$.
\ Using that $C=0$, i.e., $c=0$ and $\mu=0$, by Theorem \ref{CBI_exists}, we get that
 for all $u\in\CC_-$, the function $\RR_+\ni t \mapsto \psi(t,u)$ \ is the unique solution to the differential equation
 \[
 \partial_t \psi(t,u) = b \psi(t,u), \quad t\in\RR_+, \qquad \psi(0,u)=u.
 \]
Hence $\psi(t,u)=u\ee^{bt}=u\ee^{Bt}$, $t\in\RR_+$, where we used that $b=B$ following from \eqref{tBbeta} and $C=0$.
Consequently, \ $\psi(1,\ii \theta) - \ii \theta \varrho = \ii\theta \ee^B - \ii \theta \varrho = 0$, $\theta\in\RR$, \
 and then, using also \eqref{help15}, we have
 \begin{align}\label{help2_new}
   \EE(\ee^{\ii\theta M_n} \mid \cF_{n-1} )
       = \exp\biggl\{ \int_0^1 F(\ii \theta \ee^{Bs}) \, \dd s - \ii\theta \cA \biggr\}
       = \kappa_X(\theta),
       \qquad \theta\in\RR,\quad n\in\NN,
 \end{align}
 where \ $\kappa_X$ \ is given in \eqref{help_kappaX_kappaY}.

Therefore, using \eqref{help11} and that $\left(\sum_{k=1}^n (X_{k-1}+1)\right)^{-1}$ is $\cF_{n-1}$-measurable,
 by the properties of conditional expectation, we get that
 \[
  \EE(\ee^{\ii\theta B_n \Delta U_n}\mid \cF_{n-1})
      = \kappa_X\left( \theta \left( \sum_{k=1}^n (X_{k-1} +1) \right)^{-1} X_{n-1} \right),
      \qquad \theta\in\RR,\; n\in\NN.
 \]
Here, by Theorems \ref{convX}, \ref{main_Ad0}, and \ $\PP(w_{X_0}>0)=1$, \ for all \ $\theta \in \RR$, \ we have that
 \begin{align*}
   \theta \left(\sum_{k=1}^n (X_{k-1}+1)\right)^{-1} X_{n-1}
     &= \theta \ee^{-B} \left(\ee^{-Bn}\sum_{k=1}^n (X_{k-1}+1)\right)^{-1} (\ee^{-B(n-1)} X_{n-1})\\
     & \asP \theta \ee^{-B} \left(\frac{w_{X_0}}{\ee^B -1}\right)^{-1} w_{X_0}
             = \theta \frac{\ee^B - 1}{\ee^B}
     \qquad \text{as \ $n \to \infty$.}
 \end{align*}
Since $\kappa_X$ is the conditional characteristic function of \ $X_1-\cA$ \ given \ $X_0=0$ \ (see \eqref{help_kappaX_kappaY}),
 we have \ $\kappa_X$ \ is continuous, and hence
 \[
   \kappa_X\left(\theta \left(\sum_{k=1}^n (X_{k-1}+1)\right)^{-1} X_{n-1}\right)
   \asP \kappa_X\left(\theta \frac{\ee^B - 1}{\ee^B}\right) \qquad \text{as \ $n \to \infty$.}
 \]
By \eqref{help_kappaX_kappaY}, \eqref{regr} and the independence of \ $M_1$ \ and \ $X_0$ \ in case of \ $C=0$ \
 (which follows from the form of \ $M_1$ \ given in Lemma \ref{SDE_transform_sol_1t}
 and the independence of the Poisson random measure \ $M$ \ and \ $X_0$), \ for all \ $\theta\in\RR$, \ we have
 \begin{align*}
   \kappa_X(\theta)
    = \EE(\ee^{\ii\theta(X_1-\cA)}\mid X_0=0)
    = \EE(\ee^{\ii\theta(M_1 + \varrho X_0)}\mid X_0=0)
    = \EE(\ee^{\ii\theta M_1}\mid X_0=0)
    = \EE(\ee^{\ii\theta M_1}).
 \end{align*}
Using again Lemma \ref{SDE_transform_sol_1t}, this yields that
 \begin{align*}
   \kappa_X(\theta)
     = \exp\biggl\{\int_0^1 \int_0^\infty
                    \big(\ee^{\ii\theta r \ee^{Bu} }- 1  - \ii \theta r \ee^{Bu} \big)
                    \, \dd u \, \nu(\dd r)\biggr\},
                    \qquad \theta\in\RR.
 \end{align*}
In fact, the above formula for \ $\kappa_X(\theta)$, $\theta\in\RR$, \ follows directly from
 \ $\psi(t,u) = u\ee^{Bt}$, \ $t\in\RR_+$, $u\in\CC_{-}$, \ since, using \eqref{help2_new} and \eqref{help1_BarKorPap}, we get
 \begin{align*}
    \kappa_X(\theta) &= \exp\biggl\{ \int_0^1 F(\ii \theta \ee^{Bs}) \, \dd s - \ii\theta \cA \biggr\}\\
                     &= \exp\biggl\{ \int_0^1 \Big( a\ii\theta \ee^{Bs} + \int_0^\infty (\ee^{\ii\theta r \ee^{Bs}} - 1)\, \nu(\dd r)\Big) \, \dd s
                                    -  \ii\theta A \int_0^1 \ee^{Bs}\,\dd s \biggr\}\\
                     &=  \exp\biggl\{ \int_0^1 \int_0^\infty (\ee^{\ii\theta r \ee^{Bu}} - 1) \, \nu(\dd r)\,\dd u
                         - \ii\theta \int_0^\infty r\, \nu(\dd r) \int_0^1 \ee^{Bu}\,\dd u \biggr\},
 \end{align*}
 as desired.

Consequently,
 \begin{align*}
  \kappa_X\left(\theta \frac{\ee^B - 1}{\ee^B}\right)
   &= \EE\left( \ee^{\ii\theta \frac{\ee^B-1}{\ee^B} M_1}\right) \\
   & = \exp\biggl\{\int_0^1 \int_0^\infty
                    \big(\ee^{\ii\theta (\ee^B - 1) r \ee^{-B(1-u)} }- 1  - \ii \theta (\ee^B - 1) r \ee^{-B(1-u)} \big)
                    \, \dd u \, \nu(\dd r)\biggr\} \\
   & = \exp\biggl\{\int_0^1 \int_0^\infty
                    \big(\ee^{\ii\theta (\ee^B - 1) r \ee^{-Bu} }- 1  - \ii \theta (\ee^B - 1) r \ee^{-Bu} \big)
                    \, \dd u \, \nu(\dd r)\biggr\}.
 \end{align*}
This implies that condition (iv) of Theorem \ref{MSLTES} holds with $\mu:=\PP^{\frac{\ee^B-1}{\ee^B}M_1}$, where
 \begin{align*}
  \int_{\RR} \log^+(\vert x\vert) \,\mu(\dd x)
    & = \int_{\{ x\in \RR : \vert x\vert\geq 1\}} \log(\vert x\vert) \,\mu(\dd x)
      \leq \int_{\{x\in \RR : \vert x\vert\geq 1\}} \vert x\vert \,\mu(\dd x) \\
    &\leq \EE\left(\left\vert \frac{\ee^B-1}{\ee^B} M_1\right\vert\right)<\infty,
 \end{align*}
 since $\EE(M_1^2) =  \frac{\ee^{2B}-1}{2B} \int_0^\infty r^2 \, \nu(\dd r)<\infty$
 \ (as we have seen at the beginning of the proof).

{\sl Step 5 (application of Theorem \ref{MSLTES}):}
Using Steps 1--4, we can apply Theorem \ref{MSLTES} with our choices given before  Step 1 and in Steps 1, 3 and 4
 (for \ $\eta$, \ $P$ \ and \ $\mu$), \ and we conclude that
 \begin{align}\label{help8}
   Q_n U_n \to \eta\sum_{j=0}^\infty P^j \tZ_j
               = \frac{w_{X_0}}{\ee^B -1} \sum_{j=0}^\infty \ee^{-Bj} \tZ_j   \qquad \text{$\cF_\infty$-stably as \ $n \to \infty$,}
 \end{align}
 where \ $(\tZ_j)_{j\in\ZZ_+}$ \ is a sequence of \ $\PP$-independent and identically distributed
 random variables being \ $\PP$-independent of \ $\cF_\infty$ \ (in particular of $\eta$ or equivalently of $w_{X_0}$,
 since $w_{X_0}$ is \ $\cF_\infty$-measurable, see Step 1)
 with \ $\PP(\tZ_0 \in B) = \mu(B)$ \ for all \ $B\in\cB(\RR)$,
 \ and the series in \eqref{help8} is absolutely convergent \ $\PP$-almost surely.
This readily yields the statement, since \ $\widetilde Z_1\distre \frac{\ee^B-1}{\ee^B} M_1\distre \frac{\ee^B-1}{\ee^B} Z_1$.


\appendix

\vspace*{5mm}

\noindent{\bf\Large Appendices}

\section{Stable convergence and a  multidimensional stable limit theorem}
\label{App_1}

First, we recall the notions of stable and mixing convergences.

\begin{Def}\label{Def_HL_stable_conv}
Let \ $(\Omega,\cF,\PP)$ \ be a probability space and \ $\cG\subset \cF$ \ be a sub-$\sigma$-field.
Let \ $(\bX_n)_{n\in\NN}$ \ and \ $\bX$ \ be \ $\RR^d$-valued random variables defined on $(\Omega,\cF,\PP)$, where \ $d\in\NN$.

\noindent (i) We say that \ $\bX_n$ \ converges \ $\cG$-stably to \ $\bX$ \ as \ $n\to\infty$, \ if the conditional distribution
 \ $\PP^{\bX_n\mid \cG}$ \ of \ $\bX_n$ \ given \ $\cG$ \ converges weakly to the conditional distribution
 \ $\PP^{\bX\mid \cG}$ \ of \ $\bX$ \ given \ $\cG$ \ as \ $n\to\infty$ \ in the sense of weak convergence of Markov kernels.
It equivalently means that
 \[
   \lim_{n\to\infty} \EE_\PP(\xi \EE_\PP(h(\bX_n) \mid \cG ) )
      = \EE_\PP( \xi \EE_\PP(h(\bX) \mid \cG ) )
 \]
 for all random variables \ $\xi:\Omega\to\RR$ \ with \ $\EE_\PP(\vert \xi\vert)<\infty$ \ and for all bounded and continuous functions
 \ $h:\RR^d\to\RR$.

\noindent (ii) We say that \ $\bX_n$ \ converges \ $\cG$-mixing to \ $\bX$ \ as \ $n\to\infty$, \
 if \ $\bX_n$ \ converges \ $\cG$-stably to \ $\bX$ \ as \ $n\to\infty$, \ and \ $\PP^{\bX\mid \cG} = \PP^\bX$ \ $\PP$-almost surely, where
 \ $\PP^\bX$ \ denotes the distribution of \ $\bX$ \ on \ $(\RR^d,\cB(\RR^d))$ \ under \ $\PP$.
\ Equivalently, we can say that \ $\bX_n$ \ converges \ $\cG$-mixing to \ $\bX$ \ as \ $n\to\infty$,
 \ if \ $\bX_n$ \ converges \ $\cG$-stably to \ $\bX$ \ as \ $n\to\infty$, \ and \ $\sigma(\bX)$ \ and \ $\cG$ \ are independent,
 which equivalently means that
 \[
   \lim_{n\to\infty} \EE_\PP(\xi \EE_\PP(h(\bX_n) \mid \cG ) )
      = \EE_\PP(\xi) \EE_\PP(h(\bX))
 \]
 for all random variables \ $\xi:\Omega\to\RR$ \ with \ $\EE_\PP(\vert \xi\vert)<\infty$ \ and for all bounded and continuous functions
 \ $h:\RR^d\to\RR$.
\end{Def}

In Definition \ref{Def_HL_stable_conv}, \ $\PP^{\bX_n\mid \cG}$, \ $n\in\NN$, \ and \ $\PP^{\bX\mid \cG}$ \ are
 the \ $\PP$-almost surely unique \ $\cG$-measurable Markov kernels from \ $(\Omega,\cF)$ \ to \ $(\RR^d,\cB(\RR^d))$ \ such that for each \ $n\in\NN$,
 we have
 \[
   \int_G \PP^{\bX_n\mid \cG}(\omega,B)\,\PP(\dd \omega)
     = \PP(\bX_n^{-1}(B) \cap G)
     \qquad \text{for all \ $G\in \cG$, \ $B\in\cB(\RR^d)$,}
 \]
 and
 \[
   \int_G \PP^{\bX\mid \cG}(\omega,B)\,\PP(\dd \omega)
     = \PP(\bX^{-1}(B) \cap G)
     \qquad \text{for all \ $G\in \cG$, \ $B\in\cB(\RR^d)$,}
 \]
 respectively.
For the notion of weak convergence of Markov kernels towards a Markov kernel, see H\"ausler and Luschgy \cite[Definition 2.2]{HauLus}.
For more details on stable convergence, see H\"ausler and Luschgy \cite[Chapter 3 and Appendix A]{HauLus}.

Next, we recall a result about stable convergence of random variables, which plays an important role
 in the proof of Theorem \ref{main}.

\begin{Thm}[H\"ausler and Luschgy {\cite[Theorem 3.18]{HauLus}}]\label{Thm_HL_Thm3_18}
Let \ $\bX_n$, \ $n\in\NN$, \ $\bX$, \ $\bY_n$, \ $n\in\NN$, \ and \ $\bY$ \ be \ $\RR^d$-valued random variables on a probability space \ $(\Omega,\cF,\PP)$, \
  and \ $\cG\subset \cF$ \ be a sub-$\sigma$-field.
Assume that \ $\bX_n\to \bX$ \ $\cG$-stably as \ $n\to\infty$.
 \begin{itemize}
   \item[(a)] If \ $\Vert \bX_n - \bY_n\Vert \stoch 0$ \ as \ $n\to\infty$, \ then \ $\bY_n\to \bX$ \ $\cG$-stably as \ $n\to\infty$.
   \item[(b)] If \ $\bY_n \stoch \bY$ \ as \ $n\to\infty$, \ and \ $\bY$ \ is \ $\cG$-measurable, then
               \ $(\bX_n,\bY_n)\to (\bX,\bY)$ \ $\cG$-stably as \ $n\to\infty$.
   \item[(c)] If \ $g:\RR^d\to\RR^d$ \ is a Borel measurable function such that \ $\PP^\bX(\{ \bx\in\RR^d : \text{$g$ \ is not continuous at \ $\bx$}\})=0$, \
              then \ $g(\bX_n)\to g(\bX)$ \  $\cG$-stably as \ $n\to\infty$.
              \ Here recall that \ $\PP^\bX$ \ denotes the distribution of \ $\bX$ \ on \ $(\RR^d,\cB(\RR^d))$ \ under \ $\PP$.
 \end{itemize}
\end{Thm}

Let $(Z_n)_{n\in\NN}$ be a real-valued stochastic process on a probability space $(\Omega,\cF,\PP)$.
Let $\RR^\NN:=\{ (z_n)_{n\in\NN} : z_n\in\RR, n\in\NN \}$, and denote $(\cB(\RR))^\NN$ the product $\sigma$-field
 (also called the cylinder $\sigma$-field) on $\RR^\NN$.
We say that $(Z_n)_{n\in\NN}$ is stationary if
 $\PP^{S((Z_n)_{n\in\NN})} = \PP^{(Z_n)_{n\in\NN}}$ on the $\sigma$-field
 $(\cB(\RR))^{\NN}$, where the mapping $S:\RR^\NN\to\RR^\NN$,
 $S((z_n)_{n\in\NN}):=(z_{n+1})_{n\in\NN}$ for $(z_n)_{n\in\NN}\in\RR^\NN$, denotes the shift operator.
Clearly, $S$ is measurable, i.e., $S^{-1}(D)\in(\cB(\RR))^{\NN}$ for all $D\in (\cB(\RR))^{\NN}$.
Let $\cI_S:=\{ D\in (\cB(\RR))^{\NN} : D=S^{-1}(D) \}$.
We say that $(Z_n)_{n\in\NN}$ is ergodic if it is stationary and
 $\PP^{(Z_n)_{n\in\NN}}(D) \in\{0,1\}$ for all $D\in \cI_S$.
Next, we recall Corollary 6.26 in H\"ausler and Luschgy \cite{HauLus} about a mixing limit theorem for
 an ergodic sequence of martingale differences, which plays a role in the proof of Theorem \ref{main1_Ad}.

\begin{Thm}[H\"ausler and Luschgy {\cite[Corollary 6.26]{HauLus}}]\label{Thm_HL_Cor6_26}
Let \ $(M_k)_{k\in\NN}$ \ be an ergodic sequence of martingale differences on a probability space \ $(\Omega,\cF,\PP)$
 \ with respect to a filtration \ $(\cF_k)_{k\in\ZZ_+}$ \ (i.e., $M_k$ is $\cF_k$-measurable
 and \ $\EE_\PP(M_k\mid \cF_{k-1})=0$, \ $k\in\NN$) \ such that $\EE_\PP(M_1^2)<\infty$.
\ Then we have
 \[
    \frac{1}{\sqrt{n}} \sum_{k=1}^n M_k \to \sqrt{\EE_\PP(M_1^2)}\, N
       \qquad \text{$\cF_\infty$-mixing as \ $n \to \infty$,}
 \]
 where \ $\cF_\infty:=\sigma(\bigcup_{k=0}^\infty \cF_k)$ \ and \ $N$ \ is a standard normally
 distributed random variable independent of \ $\cF_\infty$.
\end{Thm}

For an \ $\RR^d$-valued stochastic process \ $(\bU_n)_{n\in\ZZ_+}$, \ the increments
 \ $\Delta \bU_n$, \ $n \in \ZZ_+$, \ are defined by \ $\Delta \bU_0 := \bzero$ \ and
 \ $\Delta \bU_n := \bU_n - \bU_{n-1}$ \ for \ $n \in \NN$.

Next, we recall a multidimensional stable limit theorem (a multidimensional analogue
 of a one-dimensional result due to H\"ausler and Luschgy \cite[Theorem 8.2]{HauLus}),
 which plays a crucial role in the proofs of Theorems \ref{main_Ad} and \ref{main1_Ad}.

\begin{Thm}[Barczy and Pap {\cite[Theorem 1.4]{BarPap1}}]\label{MSLTES}
Let \ $(\bU_n)_{n\in\ZZ_+}$ \ and \ $(\bB_n)_{n\in\ZZ_+}$ \ be \ $\RR^d$-valued and \ $\RR^{d\times d}$-valued stochastic processes, respectively,
 defined on a probability space \ $(\Omega,\cF,\PP)$ \ and adapted to a filtration \ $(\cF_n)_{n\in\ZZ_+}$.
 \ Suppose that \ $\bB_n$ \ is invertible for sufficiently large \ $n \in \NN$.
\ Let \ $(\bQ_n)_{n\in\NN}$ \ be a sequence in \ $\RR^{d\times d}$ \ such that
 \ $\bQ_n \to \bzero$ \ as \ $n \to \infty$ \ and \ $\bQ_n$ \ is invertible for
 sufficiently large \ $n \in \NN$.
\ Let \ $G \in \cF_\infty := \sigma(\bigcup_{k=0}^\infty \cF_k)$ \ with \ $\PP(G) > 0$.
\ Assume that the following conditions are satisfied:
 \begin{enumerate}
  \item[\textup{(i)}]
   there exists an \ $\RR^{d\times d}$-valued, \ $\cF_\infty$-measurable random matrix \ $\Beta:\Omega\to\RR^{d\times d}$ \ such that
    \ $\PP(G \cap \{\exists\,\Beta^{-1}\}) > 0$ \ and
    \[
      \bQ_n \bB_n^{-1} \stochG \Beta \qquad \text{as \ $n \to \infty$,}
    \]
  \item[\textup{(ii)}]
   $(\bQ_n \bU_n)_{n\in\NN}$ \ is stochastically bounded in \ $\PP_{G\cap\{\exists\,\Beta^{-1}\}}$-probability, i.e.,
   \[
     \lim_{K\to\infty} \sup_{n\in\NN} \PP_{G\cap\{\exists\,\Beta^{-1}\}}(\|\bQ_n \bU_n\| > K) = 0,
   \]
  \item[\textup{(iii)}]
   there exists an invertible matrix \ $\bP \in \RR^{d\times d}$ \ with \ $\varrho(\bP) < 1$
    \ such that
    \[
      \bB_n \bB_{n-r}^{-1} \stochG \bP^r \qquad
      \text{as \ $n \to \infty$ \ for every \ $r \in \NN$,}
    \]
  \item[\textup{(iv)}]
   there exists a probability measure \ $\mu$ \ on \ $(\RR^d, \cB(\RR^d))$ \ with
    \ $\int_{\RR^d} \log^+(\|\bx\|) \, \mu(\dd\bx) < \infty$
    \ such that
    \[
      \EE_\PP\bigl(\ee^{\ii\langle\btheta,\bB_n\Delta\bU_n\rangle}
                     \mid \cF_{n-1}\bigr)
      \stochGeta \int_{\RR^d} \ee^{\ii\langle\btheta,\bx\rangle} \, \mu(\dd\bx)
      \qquad \text{as \ $n \to \infty$}
    \]
    for all \ $\btheta \in \RR^d$.
 \end{enumerate}
Then
 \begin{equation}\label{conv_BU}
  \bB_n \bU_n \to \sum_{j=0}^\infty \bP^j \bZ_j \qquad
  \text{$\cF_\infty$-mixing under \ $\PP_{G \cap \{\exists\,\Beta^{-1}\}}$ \ as \ $n \to \infty$,}
 \end{equation}
 and
 \begin{equation}\label{conv_QU}
   \bQ_n \bU_n \to \Beta \sum_{j=0}^\infty \bP^j \bZ_j \qquad
   \text{$\cF_\infty$-stably under \ $\PP_{G \cap \{\exists\,\Beta^{-1}\}}$ \ as \ $n \to \infty$,}
 \end{equation}
 where \ $(\bZ_j)_{j\in\ZZ_+}$ \ denotes a sequence of \ $\PP$-independent and identically distributed \ $\RR^d$-valued random vectors
 being \ $\PP$-independent of \ $\cF_\infty$ \ with \ $\PP(\bZ_0 \in B) = \mu(B)$ \ for all \ $B \in \cB(\RR^d)$.
\end{Thm}

\begin{Rem}
(i) The series \ $\sum_{j=0}^\infty \bP^j \bZ_j$ \ in \eqref{conv_BU} and in \eqref{conv_QU} is absolutely convergent $\PP$-almost surely,
 since $\bP$ is invertible, $\varrho(\bP)<1$, $\EE_\PP(\log^+(\|\bZ_0\|))<\infty$ (by condition (iv) of Theorem \ref{MSLTES})
 and one can apply Lemma 1.3 in Barczy and Pap \cite{BarPap1} (that, for completeness, recalled below, see Lemma \ref{series_convergence}).

(ii) The random variable \ $\Beta$ \ and the sequence \ $(\bZ_j)_{j\in\ZZ_+}$ \ are \ $\PP$-independent in Theorem \ref{MSLTES},
 since \ $\Beta$ \ is \ $\cF_\infty$-measurable and the sequence \ $(\bZ_j)_{j\in\ZZ_+}$ \ is \ $\PP$-independent of \ $\cF_\infty$.
 \proofend
\end{Rem}

\begin{Lem}[Barczy and Pap {\cite[Lemma 1.3]{BarPap1}}]\label{series_convergence}
Let \ $(\bZ_j)_{j\in\ZZ_+}$ \ be a \ $\PP$-independent and identically distributed sequence of \ $\RR^d$-valued random vectors.
Let \ $\bP \in \RR^{d\times d}$ \ be an invertible matrix with \ $\varrho(\bP) < 1$.
\ Then the following assertions are equivalent:
 \begin{enumerate}
  \item[\textup{(i)}]
   $\EE_\PP(\log^+(\|\bZ_0\|)) < \infty$.
  \item[\textup{(ii)}]
   $\sum_{j=0}^\infty \|\bP^j \bZ_j\| < \infty$ \ $\PP$-almost surely.
  \item[\textup{(iii)}]
   $\sum_{j=0}^\infty \bP^j \bZ_j$ \ converges \ $\PP$-almost surely in \ $\RR^d$.
  \item[\textup{(iv)}]
   $\bP^j \bZ_j \to \bzero$ \ as \ $j \to \infty$ \ $\PP$-almost surely.
 \end{enumerate}
\end{Lem}

We note that from the proof of Lemma \ref{series_convergence} it turns out that
 for the implications \ (i) $\Rightarrow$ (ii) $\Rightarrow$  (iii) $\Rightarrow$ (iv),
 we do not need the invertibility of \ $\bP$, \ we only need it for \ (iv) $\Rightarrow$ (i).

\section{SDE representation and moments of CBI processes}
\label{section_SDE_moments}

Recall that a CBI process \ $(X_t)_{t\in\RR_+}$ \ with parameters
\ $(c, a, b, \nu, \mu)$ \ satisfying \ $\EE(X_0)<\infty$ \ and the moment condition
 \eqref{moment_condition_m_nu} can be represented as a
 pathwise unique strong solution of the SDE \eqref{SDE_atirasa_dim1}.
The SDE \eqref{SDE_atirasa_dim1_mod} is a rewriting of the SDE \eqref{SDE_atirasa_dim1}
 in a form which does not contain integrals with respect to non-compensated Poisson random measures.
Further, one can perform a transformation in order to remove the randomness
 from the drift of the SDE \eqref{SDE_atirasa_dim1_mod}, see Lemma 4.1 in Barczy et al.\ \cite{BarLiPap3}
 (which we recall below as well), and then this transformed SDE can be well-used for deriving a representation
 of the martingale differences \ $(M_k)_{k\in\NN}$ \ defined in \eqref{Mk} in terms of stochastic integrals
 with respect to a standard Wiener process and compensated Poisson integrals.

\begin{Lem}\label{SDE_transform_sol}
Let \ $(c, a, b, \nu, \mu)$ \ be a set of admissible parameters
 such that the moment condition \eqref{moment_condition_m_nu} holds.
\ Let \ $(X_t)_{t\in\RR_+}$ \ be a pathwise unique non-negative strong
 solution to the SDE \eqref{SDE_atirasa_dim1} such that \ $\EE(X_0)<\infty$.
\ Then
 \begin{align}\label{SDE_driftless}
  \begin{split}
  X_t
  &= \ee^{B(t-s)} X_s + \int_s^t \ee^{B(t-u)} A \, \dd u
     + \int_s^t \ee^{B(t-u)} \sqrt{2 c X_u} \, \dd W_u \\
  &\quad
     + \int_s^t \int_0^\infty \int_0^\infty
        \ee^{B(t-u)} z \bbone_{\{v\leq X_{s-}\}} \, \tN(\dd u, \dd z, \dd v)
     + \int_s^t \int_0^\infty \ee^{B(t-u)} r \, \tM(\dd u, \dd r)
  \end{split}
 \end{align}
 for all \ $s, t \in \RR_+$ \ with \ $s \leq t$.
\ Consequently, for the martingale differences \ $(M_k)_{k\in\NN}$ \ defined in \eqref{Mk}, we have
 \begin{align*}
   M_k & = X_k - \EE(X_k \mid \cF^X_{k-1})\\
       & = \int_{k-1}^k
           \ee^{B(k-u)} \sqrt{2 c X_u} \, \dd W_u
            + \int_{k-1}^k \int_0^\infty \int_0^\infty
           \ee^{B(k-u)} z \bbone_{\{v\leq X_{s-}\}} \, \tN(\dd u, \dd z, \dd v) \\
   &\quad
     + \int_{k-1}^k \int_0^\infty \ee^{B(k-u)} r \, \tM(\dd u, \dd r) ,
   \qquad k \in \NN .
 \end{align*}
\end{Lem}

\noindent
\textbf{Proof.}
The first statement \eqref{SDE_driftless} follows from Lemma 4.1 in Barczy et al.\ \cite{BarLiPap3}
 and the fact that \ $(X_t)_{t\in\RR_+}$ \ is a time-homogeneous Markov process.
The second statement follows from \eqref{Mk} and \eqref{SDE_driftless}.
Namely, by \eqref{Mk} and \eqref{SDE_driftless}, for each $k\in\NN$ we have
 \begin{align*}
  M_k &=  X_k - \EE(X_k \mid \cF^X_{k-1}) = X_k - \varrho X_{k-1} - \cA
       =  X_k - \ee^B X_{k-1} - A \int_0^1 \ee^{Bs}\,\dd s\\
      &=  \ee^{Bk} X_0 + \int_0^k \ee^{B(k-u)} A\,\dd u
         + \int_0^k \ee^{B(k-u)} \sqrt{2c X_u}\,\dd W_u\\
      &\quad   + \int_0^k \int_0^\infty \int_0^\infty
               \ee^{B(k-u)} z \bbone_{\{v\leq X_{s-}\}} \, \tN(\dd u, \dd z, \dd v)
         + \int_0^k \int_0^\infty \ee^{B(k-u)} r \, \tM(\dd u, \dd r)   \\
      &\quad  -\ee^B\Bigg[
              \ee^{B(k-1)} X_0 + \int_0^{k-1} \ee^{B(k-1-u)} A\,\dd u
         + \int_0^{k-1} \ee^{B(k-1-u)} \sqrt{2c X_u}\,\dd W_u \\
      &\phantom{\quad  -\ee^B\Bigg[\, }
           + \int_0^{k-1} \int_0^\infty \int_0^\infty
              \ee^{B(k-1-u)} z \bbone_{\{v\leq X_{s-}\}} \, \tN(\dd u, \dd z, \dd v) \\
      &\phantom{\quad  -\ee^B\Bigg[\, }
           + \int_0^{k-1} \int_0^\infty \ee^{B(k-1-u)} r \, \tM(\dd u, \dd r)
               \Bigg]\\
      &\phantom{=\quad}- A\int_0^1 \ee^{Bs}\,\dd s\\
      &=  \int_{k-1}^k \ee^{B(k-u)}A \,\dd u
          + \int_{k-1}^k \ee^{B(k-u)} \sqrt{2c X_u}\,\dd W_u\\
      &\quad
          + \int_{k-1}^k \int_0^\infty \int_0^\infty
              \ee^{B(k-u)} z \bbone_{\{v\leq X_{s-}\}} \, \tN(\dd u, \dd z, \dd v) \\
      &\quad + \int_{k-1}^k \int_0^\infty \ee^{B(k-u)} r \, \tM(\dd u, \dd r)
         - A \int_0^1 \ee^{Bs}\,\dd s,
\end{align*}
 yielding the assertion, since \ $\int_{k-1}^k \ee^{B(k-u)} \, \dd u = \int_0^1 \ee^{Bs} \,\dd s$.
\proofend

Note that the formula for \ $(M_k)_{k\in\NN}$ \ in Lemma \ref{SDE_transform_sol} can also be found as the first displayed
 formula in the proof of Lemma 2.1 in Huang et al.\ \cite{HuaMaZhu}.

In the next lemma, in case of a special set of admissible parameters, we can prove
 that \ $(M_k)_{k\in\NN}$ \ consists of independent and identically distributed random variables,
 and we determine the characteristic function of $M_1$ as well.

\begin{Lem}\label{SDE_transform_sol_1t}
Let \ $(X_t)_{t\in\RR_+}$ \ be a CBI process with parameters
 \ $(c, a, b, \nu, \mu)$ \ such that \ $\EE(X_0) < \infty$, \ $c = 0$, \ $\mu = 0$,
 \ and the moment condition \eqref{moment_condition_m_nu} hold.
Then
 \[
   M_k = X_k - \EE(X_k \mid \cF_{k-1}^X)
       = \int_{k-1}^k \int_0^\infty \ee^{B(k-u)} r \, \tM(\dd u, \dd r) ,
   \qquad k \in \NN ,
 \]
 and the sequence \ $(M_k)_{k\in\NN}$ \ consists of independent and
 identically distributed random variables such that the characteristic function of \ $M_1$ \ takes the form
 \[
 \EE(\ee^{\ii\theta M_1})
   = \exp\biggl\{\int_0^1 \int_0^\infty
                    \Big(\ee^{\ii\theta r \ee^{Bu} }- 1  - \ii \theta r \ee^{Bu} \Big)
                    \, \dd u \, \nu(\dd r)\biggr\},
                    \qquad \theta\in\RR.
 \]
\end{Lem}

\noindent
\textbf{Proof.}
By Lemma \ref{SDE_transform_sol}, we obtain the formula for \ $M_k$, \ $k \in \NN$.

Under the moment condition \eqref{moment_condition_m_nu}, we have
 \[
   \int_{k-1}^k \int_0^\infty \ee^{B(k-u)} r \,\dd u\, \nu(\dd r)
     = \int_{k-1}^k \ee^{B(k-u)} \,\dd u \int_0^\infty r \, \nu(\dd r)
     <\infty, \qquad k\in\NN,
 \]
 and hence, for each $k \in \NN$, the characteristic function of the random variable $M_k$ has the form
 \begin{align}\label{help_SDE_transform_sol_1t_1}
  \begin{split}
  \EE(\ee^{\ii\theta M_k})
  &= \EE\biggl(\exp\biggl\{\ii\theta\int_{k-1}^k \int_0^\infty \ee^{B(k-u)} r \, M(\dd u, \dd r) - \ii\theta\int_{k-1}^k \int_0^\infty
                    \ee^{B(k-u)} r
                    \, \dd u \, \nu(\dd r)\biggr\} \biggr) \\
  &= \exp\biggl\{\int_{k-1}^k \int_0^\infty
                    \Big(\ee^{\ii\theta r \ee^{B(k-u)} } - 1 - \ii \theta r \ee^{B(k-u)} \Big)
                    \, \dd u \, \nu(\dd r)\biggr\} \\
  &= \exp\biggl\{\int_0^1 \int_0^\infty
                    \Big(\ee^{\ii\theta r\ee^{B(1-u)}} -1 - \ii \theta r \ee^{B(1-u)} \Big)
                    \, \dd u \, \nu(\dd r)\biggr\} \\
  &= \exp\biggl\{\int_0^1 \int_0^\infty
                    \Big(\ee^{\ii\theta r\ee^{Bu}} -1 - \ii \theta r \ee^{Bu} \Big)
                    \, \dd u \, \nu(\dd r)\biggr\}
   = \EE(\ee^{\ii\theta M_1}), \qquad \theta \in \RR,
  \end{split}
 \end{align}
 yielding that $M_k$, $k\in\NN$, are identically distributed,
 where for the second equality one can use Campbell's theorem (see, e.g., Kyprianou \cite[Theorem 2.7 (ii)]{Kyp} or Kingman \cite[Section 3.2]{Kin}),
 which can be applied, since
 \[
   \int_{k-1}^k \int_0^\infty  (  \vert \ee^{B(k-u)} r \vert \wedge 1) \,\dd u\,\nu(\dd r)
     \leq \int_{k-1}^k \int_0^\infty \ee^{B(k-u)} r  \, \dd u \, \nu(\dd r) <\infty,
 \]
 as we have seen before.
It remains to check that $M_k$, $k\in\NN$, are independent.
Using Corollary 3.1 in Kingman \cite{Kin} and \eqref{help_SDE_transform_sol_1t_1},
 under the moment condition \eqref{moment_condition_m_nu}, for each $k\in\NN$ and $\theta_1,\ldots,\theta_k\in\RR$, we have
 \begin{align*}
  &\EE\Big(\exp\left\{\ii(\theta_1 M_1+ \cdots+\theta_k M_k)\right\}\Big)
     = \EE\left(\exp\left\{\ii \sum_{j=1}^k \theta_j \int_{j-1}^j \int_0^\infty \ee^{B(j-u)} r \, \tM(\dd u, \dd r) \right\}\right)\\
   & = \EE\left(\exp\left\{\ii \sum_{j=1}^k \theta_j \int_0^\infty  \int_0^\infty \ee^{B(j-u)}r\bbone_{(j-1,j]}(u) \, \tM(\dd u, \dd r) \right\}\right)\\
   & = \EE\left(\exp\left\{\ii \sum_{j=1}^k \theta_j \int_0^\infty  \int_0^\infty \ee^{B(j-u)}r\bbone_{(j-1,j]}(u) \, M(\dd u, \dd r) \right\}\right)\\
   &\phantom{=\,}
        \times \exp\left\{ -\ii \sum_{j=1}^k \theta_j \int_0^\infty  \int_0^\infty \ee^{B(j-u)}r\bbone_{(j-1,j]}(u) \, \dd u \,\nu(\dd r)  \right\}
  \end{align*}
 \begin{align*}
   & = \exp\left\{ \int_0^\infty  \int_0^\infty \left( \exp\left\{\ii \sum_{j=1}^k \theta_j \ee^{B(j-u)}r\bbone_{(j-1,j]}(u) \right\}
                     - 1 \right) \, \dd u \, \nu(\dd r) \right\}\\
   &\phantom{=\,}
        \times \exp\left\{ -\ii \sum_{j=1}^k \theta_j \int_{j-1}^j  \int_0^\infty \ee^{B(j-u)}r \, \dd u \,\nu(\dd r)  \right\}\\
   &= \exp\left\{ \sum_{j=1}^k \int_{j-1}^j  \int_0^\infty \Big( \exp\left\{\ii  \theta_j \ee^{B(j-u)}r \right\}
                     - 1 -\ii \theta_j \ee^{B(j-u)}r \Big) \, \dd u \, \nu(\dd r) \right\}\\
   &=\prod_{j=1}^k \exp\left\{ \int_{j-1}^j  \int_0^\infty \Big( \exp\left\{\ii  \theta_j r\ee^{B(j-u)} \right\}
                     - 1 -\ii \theta_j r\ee^{B(j-u)} \Big) \, \dd u \, \nu(\dd r) \right\} \\
   &=\prod_{j=1}^k \exp\left\{ \int_0^1 \int_0^\infty \Big( \exp\left\{\ii  \theta_j r\ee^{Bu} \right\}
                                - 1 -\ii \theta_j r\ee^{Bu} \Big) \, \dd u \, \nu(\dd r) \right\} \\
   &=\prod_{j=1}^k \EE(\ee^{\ii\theta_j M_1})
    = \prod_{j=1}^k \EE(\ee^{\ii\theta_j M_j}),
 \end{align*}
 yielding that for each \ $k\in\NN$, \ the random variables \ $M_1,\ldots,M_k$ \ are independent.
This implies that $M_k$, $k\in\NN$, are independent, as desired.
\proofend

Note that, using the definition of $M_k$ given in \eqref{Mk} and
 that $(X_t)_{t\in\RR_+}$ is a time-homogeneous Markov process,
 under the moment conditions $\EE(X_0^2) < \infty$ and \eqref{moment_condition_2},
 we have that $\EE(X_t^2)<\infty$, $t\in\RR_+$ (see Theorem \ref{Thm_L2}), and
 \[
    \var(M_k \mid \cF_{k-1}^X)
        = \EE( M_k^2 \mid \cF_{k-1}^X)
        = \var(X_k \mid \cF_{k-1}^X)
        = \var(X_k \mid X_{k-1}),\qquad k\in\NN,
 \]
     and
 \ $\var(X_k \mid X_{k-1} = x) = \var(X_1 \mid X_0 = x)$ \ for each \ $k\in\NN$ \ and
 \ $x \in \RR_+$.
Hence Proposition 4.8 in Barczy et al. \cite{BarLiPap3} implies the following
 formula for \ $\var(M_k \mid \cF_{k-1}^X)$, $k\in\NN$.

\begin{Pro}\label{moment_formula_2}
Let \ $(X_t)_{t\in\RR_+}$ \ be a CBI process with parameters
 \ $(c, a, b, \nu, \mu)$ \ such that \ $\EE(X_0^2) < \infty$
 \ and the moment condition \eqref{moment_condition_2} hold.
Then for each \ $k \in \NN$, \ we have
 \[
   \var(M_k \mid \cF_{k-1}^X) = V X_{k-1} + V_0 ,
 \]
 where
 \begin{align*}
  V &:= C \int_0^1 \ee^{B(1+u)} \, \dd u , \\
  V_0 &:= \int_0^\infty r^2 \, \nu(\dd r)
          \int_0^1 \ee^{2Bu} \, \dd u
          + A C
            \int_0^1
             \left( \int_0^{1-u} \ee^{Bv} \, \dd v \right)
             \ee^{2Bu} \, \dd u
 \end{align*}
 with \ $A,B$ \ and \ $C$ \ given in \eqref{help1_BarKorPap}, \eqref{tBbeta} and \eqref{help_C}.
\end{Pro}

\begin{Rem}\label{REMARK_par}
Note that \ $V_0 = \var(X_1 \mid X_0 = 0)$; \ and $V = 0$ if and only if $c = 0$ and $\mu = 0$.
Further, $V = \var(Y_1 \mid Y_0 = 1)$, where $(Y_t)_{t\in\RR_+}$ is a CBI process with parameters $(c, 0, b, 0, \mu)$
 (i.e., it is a pure branching process).
Clearly, $V$ depends only on the branching mechanism.
\proofend
\end{Rem}

\section{Lenglart's inequality}
\label{section_Lenglart}

The following form of the Lenglart's inequality can be found, e.g., in H\"ausler and Luschgy \cite[Theorem A.8]{HauLus}.

\begin{Thm}\label{Lenglart_X}
Let \ $(\xi_n)_{n\in \ZZ_+}$ \ be a nonnegative submartingale with respect to a filtration \ $(\cF_n)_{n\in\ZZ_+}$ \ and with compensator
 \ $A_n := \sum_{k=1}^n \EE(\xi_k - \xi_{k-1} \mid \cF_{k-1})$, \ $n \in \NN$, \ $A_0:=0$.
\ Then for all \ $a, b \in \RR_{++}$ \ and \ $n\in\NN$, \ we have
 \[
   \PP\Big(\max_{k\in\{0,1,\ldots,n\}} \xi_k \geq a\Big)
   \leq \frac{b}{a}
        + \PP(\xi_0 + A_n > b) .
 \]
\end{Thm}

Applying Theorem \ref{Lenglart_X} for the square of a square integrable martingale, we obtain the following corollary.

\begin{Cor}\label{Lenglart_Y}
Let \ $(\eta_n)_{n\in \ZZ_+}$ \ be a square integrable martingale with respect to a filtration \ $(\cF_n)_{n\in\ZZ_+}$ \ and
 with quadratic characteristic \ $\langle \eta\rangle_n := \sum_{k=1}^n \EE((\eta_k - \eta_{k-1})^2 \mid \cF_{k-1})$, \ $n \in \NN$,
 \ $\langle \eta\rangle_0:=0$.
\ Then for all \ $a, b \in \RR_{++}$ \ and \ $n\in\NN$, \ we have
 \[
   \PP\Big(\max_{k\in\{0,1,\ldots,n\}} \eta_k^2 \geq a\Big)
   \leq \frac{b}{a}
        + \PP(\eta_0^2 + \langle \eta\rangle_n > b) .
 \]
\end{Cor}

Note that, under the conditions of Corollary \ref{Lenglart_Y}, the quadratic characteristic process
 $(\langle \eta\rangle_n)_{n\in\ZZ_+}$ of $(\eta_n)_{n\in\ZZ_+}$ coincides with the compensator of $(\eta_n^2)_{n\in\ZZ_+}$.

\section{A version of the continuous mapping theorem for stable convergence}
\label{CMT}

If \ $\xi_n$, \ $n \in \NN$, \ and \ $\xi$ \ are random elements on a probability space $(\Omega,\cF,\PP)$
 with values in a metric space \ $(E, d)$, \ then we denote by \ $\xi_n \distrP \xi$, \ the weak convergence of the distributions of \ $\xi_n$ \
 on the space \ $(E, \cB(E))$ \ under the probability measure \ $\PP$ \ towards the distribution of \ $\xi$ \
 on the space \ $(E, \cB(E))$ \ under the probability measure \ $\PP$ \ as
 \ $n \to \infty$, \ where \ $\cB(E)$ \ denotes the Borel \ $\sigma$-algebra on \ $E$ \ induced by the given metric \ $d$.

First, we recall a version of the continuous mapping theorem (for convergence in distribution)
 which can be found, for example, in Kallenberg \cite[Theorem 4.27]{K}.

\begin{Lem}\label{Lem_Kallenberg}
Let \ $(S, d_S)$ \ and \ $(T, d_T)$ \ be metric spaces and
 \ $(\xi_n)_{n \in \NN}$ \ and \ $\xi$ \ be random elements with values in \ $S$ \ on a  probability space
 \ $(\Omega,\cF,\PP)$ \ such that \ $\xi_n \distrP \xi$ \ as \ $n \to \infty$.
\ Let \ $f : S \to T$ \ and \ $f_n : S \to T$, \ $n \in \NN$, \ be measurable
 mappings (with respect to the Borel $\sigma$-algebras
 induced by the metrics $d_S$ and $d_T$) and \ $\cC \in \cB(S)$ \ such that \ $\PP(\xi \in \cC) = 1$ \ and
 \ $\lim_{n \to \infty} d_T(f_n(s_n), f(s)) = 0$ \ if
 \ $\lim_{n \to \infty} d_S(s_n,s) = 0$ \ and \ $s \in \cC$.
\ Then \ $f_n(\xi_n) \distrP f(\xi)$ \ as \ $n \to \infty$.
\end{Lem}

Next, we formulate a version of the continuous mapping theorem for stable convergence.

\begin{Lem}\label{Lem_Kallenberg_stable}
Let \ $d\in\NN$, \ $(\Omega,\cF,\PP)$ \ be a probability space, \ $\cG\subset \cF$ \ be a sub-$\sigma$-field,
 and \ $(\xi_n)_{n\in\NN}$ \ and \ $\xi$ \ be \ $\RR^d$-valued random variables on \ $(\Omega,\cF,\PP)$ \ such that \ $\xi_n\to\xi$ \
 $\cG$-stably as \ $n\to\infty$.
\ Let \ $f:\RR^d\to\RR^d$ \ and \ $f_n:\RR^d\to\RR^d$, \ $n\in\NN$, \ be Borel measurable mappings
 such that \ $\lim_{n\to\infty}f_n(s_n)=f(s)$ \ if \ $\lim_{n\to\infty}s_n=s \in \RR^d$.
\ Then \ $f_n(\xi_n) \to f(\xi)$ \ $\cG$-stably \ as \ $n \to \infty$.
\end{Lem}

\noindent{\bf Proof.}
Since \ $\xi_n\to \xi$ \ $\cG$-stably as \ $n\to\infty$, \ by the equivalence between parts (i) and (iii)
 of Theorem 3.17 in H\"ausler and Luschgy \cite{HauLus}, for each probability measure
 \ $\QQ$ \ on \ $(\Omega,\cF)$ \ such that \ $\QQ$ \ is absolutely continuous with respect to \ $\PP$
 \ and the Radon-Nikodym derivative \ $\frac{\dd\QQ}{\dd\PP}$ \ is \ $\cG$-measurable, we have that
 \ $\QQ^{\xi_n}$ \  converges weakly to \ $\QQ^\xi$ \ as \ $n\to\infty$, \ where
 \ $\QQ^{\xi_n}$ \ and \ $\QQ^\xi$ \ denote the distribution of \ $\xi_n$ \ and \ $\xi$ \ on \ $(\RR^d,\cB(\RR^d))$
 \ under the probability measure \ $\QQ$, \ respectively.
\ Using that \ $\QQ^{\xi_n}$ \ converges weakly to \ $\QQ^\xi$ \ as \ $n\to\infty$ \ if and only if
 \ $\xi_n\distrQ \xi$ \ as \ $n\to\infty$, \ Lemma \ref{Lem_Kallenberg} implies that \
 \ $f_n(\xi_n)\distrQ f(\xi)$ \ as \ $n\to\infty$ \ or equivalently
 \ $\QQ^{f_n(\xi_n)}$ \ converges weakly to \ $\QQ^{f(\xi)}$ \ as \ $n\to\infty$ \ for each probability measure
 \ $\QQ$ \ on \ $(\Omega,\cF)$ \ such that \ $\QQ$ \ is absolutely continuous with respect to \ $\PP$
 \ and the Radon-Nikodym derivative \ $\frac{\dd\QQ}{\dd\PP}$ \ is \ $\cG$-measurable.
Using again the equivalence between parts (i) and (iii) of Theorem 3.17 in H\"ausler and Luschgy \cite{HauLus},
 we get that \ $f_n(\xi_n) \to f(\xi)$ \ $\cG$-stably \ as \ $n \to \infty$, \ as desired.
\proofend

\section*{Acknowledgement}
We acknowledge the valuable suggestions from the two referees that helped us to improve the paper.

\section*{Declaration of competing interest}

The author declares that he has no known competing financial interests or personal relationships
 that could have appeared to influence the work presented in this paper.

\end{document}